\documentclass[13pt,tikz,border=8pt,multi]{amsart}
\usepackage{microtype,amsmath,hyperref,amsthm,tikz,forest,array}
\usetikzlibrary{shapes,arrows,trees,calc}
\usetikzlibrary{shadows}
\newtheorem{theorem}{Theorem}[section]
\newtheorem{lemma}[theorem]{Lemma}
\newtheorem{conjecture}[theorem]{Conjecture}
\newtheorem{original conjecture}[theorem]{Original conjecture}
\DisableLigatures{encoding=*,family=*}
\theoremstyle{definition}
\newtheorem{definition}[theorem]{Definition}
\newtheorem{example}[theorem]{Example}

\newtheorem{cor}[theorem]{Corollary}
\theoremstyle{remark}
\newtheorem{remark}[theorem]{Remark}

\DeclareMathOperator{\lcm}{lcm}

\numberwithin{equation}{section}

\usepackage[square,numbers]{natbib}

\usepackage{listings}
\usepackage{color}
\definecolor{mygreen}{RGB}{28,172,0}
\definecolor{mylilas}{RGB}{170,55,241}

\begin{document}
\title{Well-defined Erd\"{o}s-Straus equations and L.C.M }
\author{Mohammad Arab}
\address{Department of Mathematics, Lorestan University, Khoramabad, Iran}
\email{arab.mohammad20000@gmail.com}
\curraddr{}
\date{}
\dedicatory{}

\begin{abstract}
	The Erd\"{o}s–Straus conjecture states that the equation
		$\frac{4}{n}=\frac{1}{x}+\frac{1}{y}+\frac{1}{z}$ 
		has positive integer solutions
		$x, y, z$
		for every postive integers 
		$ n \ge 2 $. We generalize the Erd\"{o}s–Straus equation, state several methods for obtaining well-defined equations, and some features of well-defined equations. Using well-defined equations, we explain the inefficiency of the solutions that did not lead to a complete proof of the problem. Finally, we prove that this conjecture cannot be completely proven by the methods expressed in various articles. For this reason, we express conjectures equivalent to the Erd\"{o}s–Straus conjecture and make the concept of conjecture clearer. The well-defined equations of Erd\"{o}s–Straus have a lot to do with the concept of the least common multiple and the greatest common divisor, and in the last section, we will independently express the functions and features based on the subject of the least common multiple.
\end{abstract}

\maketitle

\section{Introduction}
One important topic in number theory is the study of Diophantine equations, equations in which only integer solutions are permitted. The type of Diophantine equation discussed in this paper concerns Egyptian fractions \footnote{An Egyptian fraction is a representation of an irreducible fraction as a sum of unit fractions, as e.g. 7/10 = 1/2 + 1/5}, which deal with the representation of rational numbers as the sum of unit fractions\cite{A}:
\begin{equation}
\frac{a}{b}=\frac{1}{x_1}+\frac{1}{x_2}+\cdots+\frac{1}{x_r} \label{saraismylife}
\end{equation}
Specifically, in 1948 Paul Erd\"{o}s  and Ernst G. Straus formulated the following:
\begin{conjecture}
For every positive integer $  n \ge 2 $ there exist positive integers $ x, y, z $ such that:
\end{conjecture}
\begin{equation}
 \frac{4}{n}=\frac{1}{x}+\frac{1}{y}+\frac{1}{z} . \label{sara2}
\end{equation}
As stated in article\cite{U},this conjecture has attracted a lot of attention not only
among reseachers in Number Theory but also among many
people involved among the different areas in Mathematics, such as L. Bernstein\cite{B},
 M. B. Crawford\cite{C}, M. Di Giovanni, S. Gallipoli , M. Gionfriddo\cite{D}, C. Elsholtz and T. Tao\cite{E}, J. Guanouchi\cite{F,G}, L. J. Mordell\cite{H}, D. J. Negash\cite{I}, R. Obl\'{a}th\cite{J}, L. A. Rosati \cite{K}, J. W. Sander\cite{L}, R. C. Vaughan\cite{O}, K. Yamamoto\cite{P}, just to cite some of them. For example, Swett in\cite{Q} has established validity of the conjecture for all $n\le10^{14}$ and found S. E. Salez\cite{Z} for all  found $ n \le 10^{17}$. 

It is not clear if conjecture \ref{saraismylife} is true or not since there are many papers in which some authors adfirm to have proved that the conjecture is true and others in which the authors proved that it is false (see\cite{N}). L. J. Mordell\cite{H}has proven that the conjecture is true for all $n$ except possibly cases in which $n$ is congruent to $12, 112, 132, 172, 192 $ or $232$ mod $840$.

\subsection*{Schinzel method}
 A.Schinzel\cite{X}has demonstrated that if $r$ and $q$ are relatively prime integers and one can express:
\begin{equation} \label{sara 0021}
\frac{4}{rt+q}=\frac{1}{X(t)}+\frac{1}{Y(t)}+\frac{1}{Z(t)}
\end{equation}
With $X(t)$, $Y(t)$ and $Z(t)$ being integer polynomials in $t$ with positive leading coefficients, then $q$ cannot be a quadratic residue mod $r$.\\ In all of these papers, the authors prove that there exist solutions for the Erd\"{o}s-Straus Conjecture, but noone of them found explicit solutions.\\

In this article, considering the equation\eqref{saraismylife} as a partition of positive rational numbers we will address the following:
\begin{enumerate} 
\item The reason for the inefficiency of the methods used to solve the equation will be stated.
\begin{enumerate}
\item Schinzel method.
\item Identity $ 1=\frac{a-1}{a}+\frac{1}{a}$.
\end{enumerate}

\item Methods for obtaining well-defined equations for the Erd\"{o}s-Straus equation which can be generalized to the equation\eqref {saraismylife} will be presented.
\begin{enumerate}
\item Algebraic formulas, partitions of positive integers.\eqref{sara 10},\eqref{saraismylife.very}
\item Using any partition any positive integer.\eqref{sara31}
\item By identity $a = a$ and auxiliary variable.\eqref{sara arab}
\item Use of distributive property.\eqref{sara 097}
\item Using the property, the least common multiple.\eqref{sara.aras}
\end{enumerate}

\item Some features of well-defined Erd\"{o}s-Straus equations will be described.
\begin{enumerate}
\item Well-defined equations based on the partition of numbers in the form $kb-1$.\eqref{sara 500}
\item Different forms of a well-defined equation.\eqref{sara..fat}
\item Well-defined equations triangle.\eqref{sara.fat}
\item Well-defined equation based on positive integer partitions with positive integer sizes.\eqref{sarabor}
\item Well-defined equations based on numbers in the form $a^{b} $.\eqref{sara....f},\eqref{SARA}
\item The expression of a pair of twin equations.\eqref{sarababa}
\end{enumerate}

\item The domain of definition $n$ in the Erd\"{o}s-Straus equation .
\item Conjectures equivalent to the Erd\"{o}s-Straus conjecture and its generalization.
\item Arithmetic functions based on the least common multiple and their properties.
\end{enumerate}

\begin{original conjecture}
For every positive integer $ \frac{k}{n} \le 3 $ there exist positive integers $ \alpha,\beta,\gamma $ such that:
\begin{equation} 
 \frac{k}{n}=\frac{1}{\alpha}+\frac{1}{\beta}+\frac{1}{\gamma}. \label{sara3}
\end{equation}
\end{original conjecture}

\section{Definitions of partitions}

In this section, we will explain the definitions that we are dealing with in the process of the article.

\subsection*{Partitions of Integers}

In the book principles and techniques of combinations\cite{T}, the definition of integer partitions is as follows:
\begin{definition}
A partition of a positive integer $n$ is a collection of positive integers whose sum is $n$. Since the ordering is immaterial, we may regard 
a partition of $n$ as a finite nonincreasing sequence $ n_1 \ge n_2 \ge \cdots \ge n_l$ of positive integers such that $\sum_{i=1}^l
n_i= n$. So If $n =n_1 +n_2 +\cdots+n_l$ is a partition of $n$, we say that $n$ is partitioned into $l$ parts of sizes $n_1,n_2,\cdots,n_l$ respectively.
\end{definition}
We express two other definitions similar to this one.
\subsection*{Generalized partitions of Integers} 
\begin{definition}
A partition of a positive integer $n$ is a collection of positive rational numbers whose sum is $n$. Since the ordering is immaterial, we may regard 
a partition of $n$ as a finite nonincreasing sequence $ n_1 \ge n_2 \ge \cdots \ge n_l$ of positive rational numbers such that $\sum_{i=1}^l
n_i=n$. So If $ n=n_1+n_2+\cdots+n_l $ is a partition of $n$, we say that $n$ is partitioned into $l$ parts of sizes $n_1,n_2,\cdots,n_l$ respectively.
\end{definition}

\subsection*{Partitions of rational numbers}

\begin{definition}
A partition of a positive rational numbers $n$ is a collection of positive rational numbers whose sum is $n$ . Since the ordering is immaterial, we may regard a partition of $n$ as a finite nonincreasing sequence $ n_1 \ge n_2 \ge \cdots \ge n_l$ of positive rational numbers such that $\sum_{i=1}^l
n_i=n$. So If $ n=n_1+n_2+\cdots+n_l $ is a partition of $n$, we say that $n$ is partitioned into $l$ parts of sizes $n_1,n_2,\cdots,n_l$ respectively.
\end{definition}
\begin{remark}
	In the process of the article, wherever it is expressed, the partitioning of positive integers with rational sizes means that there are a number of pairs of rational numbers in the sizes whose face and denominator are relatively prime.
	\[
	\forall n \in \mathbb{N}  \quad,\quad \forall l \in \mathbb{E}=2\mathbb{N}  \quad,\quad \exists n_1, n_2,\cdots, n_l,m_1, m_2,\cdots, m_{\frac{l}{2}} \in \mathbb{N} \]
	\[
	\quad  m_{\frac{i}{2}}\big{|} n_{i-1} + n_{i} \quad,\quad 2\le i \le l \quad,\quad i \in \mathbb{E}
	\] 
\[
n=\frac{n_1 + n_2}{m_{1}}+\frac{n_3 + n_4}{m_{2}}+\cdots+\frac{n_{l-1} + n_{l}}{m_{\frac{l}{2}}}
\]
\end{remark}
\begin{example}
	Partition of numbers in form $(4z+1) $ with rational sizes:
\[(4z+1)=4z+\frac{1}{2}+\frac{1}{2}\]
Partition of numbers in form $(4z+2) $ with positive integer sizes:
\[(4z+2)=4z+1+1\]
\end{example}
\section{Formulation partitions of positive rational numbers}

The partitions of positive rational numbers  can be obtained from the partitions of positive integers, which will be shown in the following theorems:
\begin{theorem} \label{sara 10}
Multiplying the partition $ n = n_1 + n_2 + \cdots + n_l $ in $ \frac{1}{\lcm\left(n_1, n_2,\cdots, n_l,r\right)} $  the set of partitions of positive rational numbers  will be obtained:
\[ \forall n \in \mathbb{N} \quad,\quad \forall l \in \mathbb{N}_{n} \quad,\quad \exists n_1, n_2,\cdots, n_l \in \mathbb{N} \quad,\quad  \forall r \in \mathbb{N}  \]
\begin{equation*}
\frac{n}{\lcm\left(n_1, n_2,\cdots, n_l,r\right)}=\frac{1}{\frac{\lcm\left(n_1, n_2,\cdots, n_l,r\right)}{n_1}}+\frac{1}{\frac{\lcm(n_1, n_2,\cdots, n_l,r)}{n_2}}+\cdots+\frac{1}{\frac{\lcm(n_1, n_2,\cdots, n_l,r)}{n_l}}.
\end{equation*}
\end{theorem}
\begin{proof}
Given the definition of the least common multiple\cite{sameh}, the correctness of the theorem is obvious.
\end{proof}

\begin{remark}
The above theorem can be generalized from partitions of integers with integer sizes to partitions of integers with rational sizes as follows:
\[
\forall n \in \mathbb{N}  \quad,\quad \forall l \in \mathbb{E}=2\mathbb{N}  \quad,\quad \exists n_1, n_2,\cdots, n_l,m_1, m_2,\cdots, m_{\frac{l}{2}} \in \mathbb{N} \]
\[
 m_{\frac{i}{2}}\big{|} n_{i-1} + n_{i} \quad,\quad 2\le i \le l \quad,\quad i \in \mathbb{E}
\] 
\[
n=\frac{n_1 + n_2}{m_{1}}+\frac{n_3 + n_4}{m_{2}}+\cdots+\frac{n_{l-1} + n_{l}}{m_{\frac{l}{2}}}
\]
In this case:
\begin{align*}
\frac{n}{\lcm\left(n_1, n_2,\cdots, n_l,r\right)}&=\frac{1}{\frac{\lcm(n_1, n_2,\cdots, n_l,r)m_{1}}{n_1}}+\frac{1}{\frac{\lcm(n_1, n_2,\cdots, n_l,r)m_{1}}{n_2}}+\\
&+\cdots +\frac{1}{\frac{\lcm(n_1, n_2,\cdots, n_l,r) m_\frac{l}{2}}{n_l}}.
\end{align*}
\end{remark}
\begin{theorem}\label{sara31}
From each partition of positive integers, a subset of Partitions of rational numbers can be obtained. The number of subsets is equal to the number of different sizes.
\end{theorem}
\begin{proof}
Consider the partition of positive integers as follows:
\[ \forall n \in \mathbb{N} \quad,\quad \forall l \in \mathbb{N}_{n} \quad,\quad \exists n_1, n_2,\cdots, n_l \in \mathbb{N} \]
\[ n = n_1 + n_2 +\cdots+ n_l   \]
We will get the fraction $\frac{n}{\lcm(n_1, n_2,\cdots, n_l)}$ in the form of an irreducible fraction: \[\frac{n}{\lcm(n_1, n_2,\cdots, n_l)}=\frac{r}{s}\] 
In the event that $(r,s)=1 $ and we multiply the positive integer multiples of $ s $ by $ n $:
\[ n(sg) \ne n_1 + n_2 + \cdots + n_l. \]
Assuming $ n_1 \ne n_2 \ne \cdots \ne n_l $ the subsets and their number are obtained as follows:
\[subset(1): n(sg) = n_1u_{1} + n_2 + \cdots + n_l \]
\[subset(2): n(sg) = n_1+ n_2u_{2} +\cdots + n_l \]
\[\vdots\]
\[subset(l): n(sg) = n_1+ n_2 + \cdots + n_lu_{l} \]
The number of subsets is equal to the number of different sizes.\\
To get each subset you have to multiply it by $\frac{1}{\lcm(n_1, n_2,\cdots, n_l,s)gu_{i}} $ for $1\le i \le l$.\\
For the subset(1),we will have multiply the partition by $\frac{1}{\lcm(n_1, n_2,\cdots, n_l,s)gu_{1}}$:
\begin{align*}
\frac{n(sg)}{\lcm(n_1, n_2,\cdots, n_l,s)gu_{1}}&=\frac{n_1u_{1}}{\lcm(n_1, n_2,\cdots, n_l,s)gu_{1}}+\frac{n_2}{\lcm(n_1, n_2,\cdots, n_l,s)gu_{1}}\\
&+\cdots+\frac{n_l}{\lcm(n_1, n_2,\cdots, n_l,s)gu_{1}} 
\end{align*}
Which can be written as follows:
\begin{align*} \label{mohammad}
\frac{n}{\left(\frac{\lcm(n_1, n_2,\cdots, n_l,s)}{s}\right)u_{1}}&=\frac{1}{\left(\frac{\lcm(n_1, n_2,\cdots, n_l,s)}{n_1}\right) g}+\frac{1}{\left(\frac{\lcm(n_1, n_2,\cdots, n_l,s)}{n_2}\right)gu_{1}}\\
&+\cdots+\frac{1}{\left(\frac{\lcm(n_1, n_2,\cdots, n_l,s)}{n_l}\right)gu_{1}} 
\end{align*}
In the same way, the whole subset will be obtained and the structure of the equation is based on the following equation:
\[
\frac{r}{s}=\frac{1}{\frac{\lcm(n_1, n_2,\cdots, n_l,s)}{n_1}}+\frac{1}{\frac{\lcm(n_1, n_2,\cdots, n_l,s)}{n_2}}+\cdots +\frac{1}{\frac{\lcm(n_1, n_2,\cdots, n_l,s)}{n_l}}
\]
And we can do the same process for other subsets.
\end{proof}

\begin{remark}
In the $ u_{j}$ for $j\in \mathbb{N} $  position there is a condition that $u$ is obtained independently:
\[ n(sg) = n_1+ n_2 + \cdots+ n_{j}u_{j}+\cdots+n_{l} \quad,\quad 1\le j \le l \]
\[n_{j}u_{j}=n(sg)-\left(\sum_{j\ne i, i=1}^l n_{i}\right) \]
\[u_{j} = \left(\frac{ns}{n_{j}}\right)g-\left(\frac{\sum_{j\ne i, i=1}^l n_{i}}{n_{j}}\right) \]
In this case, the following two divisibility must be established:
\begin{itemize}
\item $n_{j} \big{|} \left(\sum_{j\ne i, i=1}^l n_{i}\right) $\\
\item $n_{j} \big{|} ns$.
\end{itemize}
\end{remark}
\begin{remark}
Theorem \ref{sara31} can also be expressed by generalizing partitions of positive integers .\\
In theorem \ref{sara31}, forall $m_{\frac{i}{2}}$ equals one.
\end{remark}
\begin{example}
To further illustrate theorem \ref{sara31} we will give an example below.
Consider partition number $24$:
\[
24=1+2+21
\]
Now based on what we stated in the proof of the theorem \ref{sara31}, for the partition number $24$ will be obtained as follows:
\[
   subset(1):24(7g)=1u_{1}+2+21
\]
\[
subset(2):24(7g)=1+2u_{2}+21
\]
\[
subset(3):24(7g)=1+2+ 21u_{3}
\]
For example, subset(2):
\[
\frac{24(7g)}{42gu_{2}}=\frac{1}{42gu_{2}}+\frac{2u_{2}}{42gu_{2}}+\frac{21}{42gu_{2}} \\ 
\]
Which can be written as follows:
\[
\frac{4}{u_{2}}=\frac{1}{42gu_{2}}+\frac{1}{21g}+\frac{1}{2gu_{2}}\\
\]
According to subset(2) we now enter $u_{2}$:
\[
   u_{2}=84g-11
\]
\[
\frac{4}{84g-11}=\frac{1}{42g(84g-11)}+\frac{1}{21g}+\frac{1}{2g(84g-11)} \quad,\quad  \forall g \in \mathbb{N}  
\]
The structure of the equation is based on the following equation:
\[
\frac{4}{7}=\frac{1}{42}+\frac{1}{21}+\frac{1}{2}
\]
And we can do the same process for other subsets.
\end{example}
\begin{theorem}\label{saraismylife.very}
Equation\eqref{saraismylife}can be well defined by generalizing positive integers partitions.
\end{theorem}
\begin{proof}
First, consider the following set of positive integer partitions:\\
\[ \forall k\in \mathbb{N} \quad,\quad \forall r \in \mathbb{O}=2\mathbb{N}+1 \quad,\quad 1\le i \le \frac{r-1}{2} 
\]
\[
\exists n,c_{i},b_ {2i-1},b_ {2i} \in \mathbb{N} \quad,\quad c_i\mid b_ {2i-1} +b_ {2i} \quad,\quad \exists q \in \mathbb{N} 
\]
\begin{equation}\label{sara 5}
k\left( \frac{b_1b_2\cdots b_{r-1}}{\gcd\left(b_1,b_2,\cdots,b_{r-1}\right)}\right)q=n+\sum_{i=1}^\frac{r-1}{2}
	  \frac{b_ {2i-1} +b_ {2i}}{c_i}
\end{equation}
By multiplication $\frac{1}{\left( \frac{b_1b_2\cdots b_{s-1}}{\gcd(b_1,b_2,\cdots,b_{s-1})}\right)nq}$ in partition \eqref{sara 5} diophantine equation \eqref{saraismylife} is generally well defined:
	\begin{equation}
	\frac {k}{n}=\frac {1}{\left(\frac{b_1b_2\cdots b_{r-1}}{\gcd(b_1,b_2,\cdots,b_{r-1})}\right)q}+\sum_{i=1}^\frac{r-1}{2}
	\frac{1}{\frac{\left(\frac{b_1b_2\cdots b_{r-1}}{\gcd(b_1,b_2,\cdots,b_{r-1})}\right)nqc_i}{b_{ 2i-1} +b_ {2i}}}
	\end{equation}
Thus equation \eqref{saraismylife} was expressed in the form of well-defined.

\end{proof}

Erd\"{o}s–Straus conjecture is actually a partition of positive rational numbers in the form of $\frac{4}{n}$ with sizes in the form of unit fractions.\\
In this case, according to theorem \ref{saraismylife.very} for two-part and three-part partitions, can be expressed.

\begin{cor}
	Two-part partitions and Three-part partitions (Original conjecture) of positive rational numbers are well-defined by Two-part and Three-part partitions of positive integers.
\end{cor}
\begin{proof}
Two-part partitions:
\[ \forall k \in \mathbb{N} \quad,\quad \forall b \in \mathbb{N} \quad,\quad \exists a,q \in \mathbb{N} \quad,\quad kq>1 \quad,\quad a+b=kbq \] 
Other forms of this partition are as follows:
\[  a+b=kbq  \rightarrow  a=kbq-b=b(kq-1) \rightarrow a \equiv -b \pmod{kb}\]
Multiplying Two-part partitions of positive integers by $ \frac{1}{abq} $ will give Two-part rational partitions:
\[ \frac{k}{a}=\frac{1}{bq}+\frac{1}{aq} \rightarrow \frac{k}{(kq-1)b}=\frac{1}{bq}+\frac{1}{(kq-1)bq} \]\\

Three-part partitions:
\[\forall k \in \mathbb{N} \quad,\quad \exists a \in \mathbb{N} \quad,\quad \exists b,c,d \in \mathbb{N} \quad,\quad c\mid b+d \quad,\quad \exists q \in \mathbb{N}  \]
\begin{equation} \label{sara 3000}
k\left(\frac{bd}{\gcd(b,d)}\right)q=a+\frac{b}{c}+\frac{d}{c} 
\end{equation}
Three-part partitions of positive rational numbers for all $ k \in \mathbb{N} $ are obtained using partitions\eqref{sara 3000}:
\[a=k\left(\frac{bd}{\gcd(b,d)}\right)q-\left(\frac{b+d}{c}\right)\]
According to the definition of the greatest common divisor can be written:
\[ \gcd(b,d)=D \rightarrow \gcd\left(b_{1}=\frac{b}{D},d_{1}=\frac{d}{D}\right)=1 \]
\[a=kbd_{1}q-\frac{b_{1}D-d_{1}D}{c} \]
Therefore:
\[a \equiv -\frac{D}{c}(b_{1}+d_{1}) \pmod{kbd_1} \]
So the Three-part partitions of rational numbers are expressed as follows.\\
The partition \eqref{sara 3000} will be multiplied by $ \frac{1}{\left(\frac{bd}{\gcd(b,d)}\right)qa}$:
\[\frac{k}{a}=\frac{1}{\left(\frac{bd}{\gcd(b,d)}\right)q}+\frac{1}{\frac{\left(\frac{bd}{\gcd(b,d)}\right)qac}{b}}+\frac{1}{\frac{\left(\frac{bd}{\gcd(b,d)}\right)qac}{d}}.\]
\end{proof}
The difference of unit fractions from each other or their equality with each other in this corollary can be checked in terms of the expressed partition because it includes both cases.

\section{Well-defined equations}

There are several ways to obtain well-defined equations for equation \eqref{saraismylife}. The solution for obtaining well-defined equations from partitions was stated in the previous section.
In this section, the solution to obtain well-defined equations using identity and auxiliary variables will be expressed.\\
We now generalize equation\eqref{sara 0021} and what is obtained is the well-defined equations for equation\eqref{saraismylife}. \\
The domain of definition $n$ in any well-defined equation in the form of equation\eqref{sara 0021} is a subset of natural numbers. In general we will generalize the domain $n$ to a family of subsets of natural numbers, for this reason in equation\eqref{sara 0021}, we will generalize $r,q$ to the relations.\\
Assuming $ S $ and $ G$ as multivariate relations, families belonging to the domain $n$ and the generalization of the equation\eqref{sara 0021} can be expressed as follows:
\[
M(t)\equiv \mp G \pmod{S}
\]
\[
M(t)=St \mp G
\]
\[
\frac{k}{M(t)}=\frac{1}{\beta_{1}(t)}+\frac{1}{\beta_{2}(t)}+\cdots \pm\frac{1}{\beta_{s}(t)}
\]
$t$ is the positive leading coefficient.\\
The two relations that underlie well-defined equations in many theorems are as follows:
\begin{enumerate}
\item The polynomial relation underlies many of the well-defined equations expressed in the theorems presented:
\[
S(w)=aw+ b \quad,\quad G(w^{'})=cw^{'}+ d
\]
$w$ and $w^{'}$ positive leading coefficients.\\
And families belonging to the scope of $n$ is expressed as follows:
\[
M(t)=(aw+ b)t\pm (cw^{'}+ d)
\]
\[
M(t) \equiv \pm (cw^{'}+ d) \pmod{(aw+ b)} 
\]
One of the main points of the article is to present well-defined equations in the form of this type of polynomial, which is presented in this section and the following sections.
\item
$S$ and $ G $ are expressed in the form of two rational relations: 
\[
S(a,b)=\frac{a}{b} \quad,\quad G(c,d)=\frac{c}{d}  
\]
\[ M(t)=S(a,b)t\pm G(c,d) \]
\[
M(t) \equiv \pm G(c,d) \pmod{S(a,b)} 
\]
In equation\eqref{sara mohammad} we have given a well-defined example of this equation.
\end{enumerate}
As we have said, a family is a set of sets. If we call a family a $2$-set, then in general we can express it in the form of a mathematical expression:
\[
z\left(w,a_{1},a_{2},\cdots,a_{s-1}\right):=(\cdots(((kw-b_{1})a_{1}-b_{2})a_{2}-b_{3})a_{3}-\cdots-b_{s-1})a_{s-1}-b_{s}.
\]
$z\left(w,a_{1},a_{2},\cdots,a_{s-1}\right)$ is called $s$-set.
\begin{theorem}
An equation based on $s$-set for partitions of positive rational numbers can be expressed as follows:
\[
\forall k\in\mathbb{N} \quad,\quad  \forall w\in\mathbb{N} \quad,\quad \forall s\in\mathbb{N} \quad,\quad \exists b_{1},b_{2},\cdots,b_{s} \in \mathbb{N}
\]
\[
w=\lcm\left(b_{1},b_{2},\cdots,b_{s}\right) \quad,\quad v=(\cdots(((kw-b_{1})a_{1}-b_{2})a_{2}-b_{3})a_{3}-\cdots-b_{s-1})a_{s-1}-b_{s} 
\]
\[
\frac{k}{v}=\frac{1}{wa_{1}a_{2}\cdots a_{s-1}}+\frac{1}{\frac{vw}{b_{1}}}+\frac{1}{\frac{vwa_{1}}{b_{2}}}+\frac{1}{\frac{vwa_{1}a_{2}}{b_{3}}}+\frac{1}{\frac{vwa_{1}a_{2}a_{3}}{b_{4}}}+\cdots+\frac{1}{\frac{vwa_{1}a_{2}\cdots a_{s-1}}{b_{s}}}
\]
For all $a_{1},a_{2},\cdots,a_{s-1}\in\mathbb{N}$ s.t $v>0$.
\end{theorem}
\begin{proof}
Consider partitions like this:
\[
kwa_{1}a_{2}\cdots a_{s-1}=v+b_{1}a_{1}a_{2}\cdots a_{s-1}+b_{2}a_{2}a_{3}\cdots a_{s-1}+\cdots+b_{s}
\]
Multiplying these types of partitions by the following fraction will give the equation in question:
\[
\frac{1}{vwa_{1}a_{2}\cdots a_{s-1}}.
\]
\end{proof}
According to the theorems of the previous section, the well-defined equations of equation\eqref{saraismylife} can be obtained using the identity $ a = a $. \\
Well-defined equations for the original conjecture \ref{sara3} will now be expressed with respect to the identity $ a = a $. \\
In the following theorem we explain how to obtain well-defined equations using the auxiliary variable.
\begin{theorem} \label{sara arab}
Well-defined equations the original conjecture is obtained by using the auxiliary variable in the identity a=a:
\begin{enumerate}
\item  \label{sara 401}
 \[\forall k\in \mathbb{N}-\{1\} 
\]
\[
 \forall m\in\{1 \le m\le k-1 \mid m\in\mathbb{N} \} \quad,\quad \exists q\in \mathbb{N} \quad,\quad mq>1
\]
\[
 \exists v\in \mathbb{N} \quad,\quad k-m\mid v 
\]
\[
 \frac {k}{mt-v}=\frac{1}{\frac{1}{k-m}(mt-v)}+\frac{1}{t}+\frac{1}{\frac{t}{v}(mt-v)} 
\]
\[
\forall  t\in\{ t \mid t=vqr \:,\: r\in \mathbb{N}  \} 
\]

\item
\[
\forall k\in \mathbb{N} \quad,\quad \forall  v\in\{2\le v \le k \mid v\in \mathbb{N} \}
\]
\[
\exists z\in\{1 \le z\le v-1 \mid z\in \mathbb{N} \} \quad,\quad
 v>z
\]
\[
 \frac{k}{kt-v}=\frac{1}{t}+\frac{1}{\frac{t(kt-v)}{z}}+\frac{1}{\frac{t(kt-v)}{v-z}}
\]
\[ 
\forall t \in\left\{t\mid t=lx \:,\: \lcm(z,v-z)=l\:,\:klx>v\:,\:l,x\in\mathbb{N}\right\}
\]
\item 
\[\forall m\in \mathbb{N}  \quad,\quad \exists v \in \mathbb{N} \quad,\quad  v\mid m
\]
\[
 \frac{\frac{m}{v}}{mt-v}=\frac{1}{tv}+\frac{1}{t(mt-v)} 
\]
\[
\forall t\in
\begin{cases}
\mathbb{N}-\{1\}  & \mbox{if } v=m\\
\mathbb{N}  & \mbox{if } v\ne m
\end{cases}
\]

\item
\[\forall m\in \mathbb{N} \quad,\quad \exists r\in \mathbb{N} \quad,\quad r\mid m \quad,\quad r=ld \quad,\quad l,d \in \mathbb{N} \]
\[
 \frac{\frac {m}{r}}{mt-(l+d)}=\frac{1}{lt(mt-(l+d))}+\frac{1}{rt}+\frac{1}{dt(mt-(l+d))} 
\]
\[
\forall t\in
\begin{cases}
\mathbb{N}-\{1\}  & \mbox{if } r=m \quad ,\quad l+d=m\\
\mathbb{N}  & \mbox{if } r\ne m
\end{cases}
\]

\item
\[ 
\forall k \in \mathbb{N} \quad,\quad \exists w \in \mathbb{N} \quad ,\quad (2k-1)w\ne1  
\]
\begin{align*}
&\frac{k}{(kw-1)t+((k-1)w-1)}=\frac{1}{w(t+1)}+\\
&+\frac{1}{w((kw-1)t+((k-1)w-1))}+\frac{1}{(t+1)((kw-1)t+((k-1)w-1))}  
\end{align*}
\[
\forall t\in
\begin{cases}
\mathbb{N}-\{1\}  & \mbox{if } (2k-1)w=2\\
\mathbb{N}  & \mbox{if } (2k-1)w>2 \quad,\quad k\in\{1,2\}
\\
\mathbb{N} \cup \{0\} & \mbox{if }(2k-1)w>2 \quad,\quad k>2

\end{cases}
\]
\textrm{In another form}:
\[ 
\forall w \in \mathbb{N} \quad,\quad S\equiv -1 \pmod{w} \quad,\quad G\equiv -1 \pmod{w} \]
\[
\frac{S+1}{k}=\frac{G+1}{k-1}=w \quad,\quad k \in \mathbb{N} \quad,\quad M\equiv G \pmod{S}
\]
\begin{equation*}
\frac{k}{M}=\frac{k}{St+G}=\frac{1}{w(t+1)}+\frac{1}{w(St+G)}+\frac{1}{(t+1)(St+G)}  
\end{equation*}

\item
\[ \forall k \in \mathbb{N} \quad,\quad \exists w \in \mathbb{N}\quad ,\quad kw>1\]
\[ \frac{k}{(kw-1)(t-w)}=\frac{1}{w(t-w)}+\frac{1}{wt(kw-1)}+\frac{1}{t(kw-1)(t-w)} \]
\[
 \forall t \in\{t > w\mid t \in \mathbb{N}\}    
\]

\item
\[
\forall k \in \mathbb{N} \quad,\quad \exists w \in \mathbb{N} \quad ,\quad kw>1
\]
\[
\frac{k}{(kw-1)(v-1)}=\frac{1}{w(v-1)}+\frac{1}{vw(kw-1)}+\frac{1}{vw(kw-1)(v-1)}
\]
\[
\forall v \in \mathbb{N}- \{1\}
\]
\item\label{SARA}
\[
\forall A\in \mathbb{N}-\{1\}
\]
We consider the decomposition of $A$ into the prime factors:
\[
 A=p_{1}^{a_{1}}p_{2}^{a_{2}}\cdots p_{i}^{a_{i}} 
\]
\[
A_{1}=p_{1}p_{2}\cdots p_{i}
\]
\[
\forall b\in \mathbb{N}  \quad,\quad \exists c \in \mathbb{N}
\begin{cases}
	b>c \quad,\quad b\ne 1  & \mbox{if } A=A_{1}\\
	b\ge c  & \mbox{if } A\ne A_{1}
\end{cases}
\]
\[
\exists n \in \mathbb{N}  \quad,\quad A^{b}-A_{1}^{c}>n
\]
\begin{align*}
&\frac{A^b}{p_{1}^{a_{1}}p_{2}^{a_{2}}\cdots p_{i}^{a_{i}}\left(A^{b}-A_{1}^{c}-n\right)n }&=\\&=\frac{1}{\left(p_{1}^{a_{1}}p_{2}^{a_{2}}\cdots p_{i}^{a_{i}}\right)n}+\frac{1}{p_{1}^{a_{1}+c}p_{2}^{a_{2}+c}\cdots p_{i}^{a_{i}+c}\left(A^{b}-A_{1}^{c}-n\right)n}&+\\&+\frac{1}{p_{1}^{a_{1}}p_{2}^{a_{2}}\cdots p_{i}^{a_{i}}\left(A^{b}-A_{1}^{c}-n\right)}
\end{align*}
\item
\[
\forall k \in \mathbb{N} \quad,\quad \exists r,v,w,s,l,q \in \mathbb{N} \quad,\quad (krv-l)s\ge rw \]
\[
(krv-l)s-rw\mid rvwls
\]
\[
\left(\frac{k}{wl}=\frac{1}{rvw}+\frac{1}{vls}+\frac{1}{\frac{rvwls}{(krv-l)s-rw}}\right) \equiv \left(\frac{1}{s}=\frac{krv-l}{rw}-\frac{vl}{q}\right)
\]

\item
\[
\forall k\in \mathbb{N}\quad,\quad\forall m\in\mathbb{N} \quad,\quad \exists v,b\in\mathbb{N}\quad,\quad kmv>\mp b
\]
\[
 \exists s\in\mathbb{N}\quad,\quad \mp sb\pm v \mid mvs(kmv\pm b)
\]
\[
\frac{k}{kmv\pm b}=\frac{1}{mv}\mp\frac{1}{ms(kmv\pm b)}+\frac{1}{\frac{mvs(kmv\pm b)}{\mp sb\pm v}}
\]

\end{enumerate}
\end{theorem}
\begin{proof}
The correctness of the equations can be proved in two ways:\\
 The first method using the sum of unit fractions and reaching the other side of the equation and the second method using the identity and the auxiliary variable.\\
The first method is simple, so we use it for the first equation of the theorem and the rest of the equations will be proved as the first equation but the second method will be stated for all equations.
\begin{enumerate}
\item Consider the following identity:
\[
k=k  
\]
Adding auxiliary variables $m$ and $\frac{v}{t}$:
\[
k=k-m+m-\frac{v}{t}+\frac{v}{t}=k-m+\frac{mt-v}{t}+\frac{v}{t}
\]  
We multiply the partition by $ \frac {1}{mt-v} $:
\[
 \frac {k}{mt-v}=\frac{1}{\frac{mt-v}{k-m}}+\frac{1}{t}+\frac{1}{\frac{t}{v}(mt-v)}.
\]
Another proof is as follows:
\begin{align*}
	 &\frac{1}{\frac{mt-v}{k-m}}+\frac{1}{t}+\frac{1}{\frac{t}{v}(mt-v)}&=\\
&=\frac{\frac{(k-m)t}{v}+\frac{mt-v}{v}+1}{\frac{t}{v}(mt-v)}&=\\ 
&=\frac{(k-m)t+mt-v+v}{t(mt-v)}&=\\
 &=\frac{k}{(mt-v)}.
	 \end{align*}

\item 
Consider the following identity:
\[
kt=kt 
\]
Adding auxiliary variables $v$ and $z$:
\[
 kt=kt-v+z+v-z 
\]
We multiply the partition by $ \frac {1}{t(kt-v)} $:
\[
 \frac{k}{kt-v}=\frac{1}{t}+\frac{1}{\frac{t(kt-v)}{z}}+\frac{1}{\frac{t(kt-v)}{v-z}}.
\]
\item 
Consider the following identity:
\[
 mt=mt
 \]
Add auxiliary variable $ v $:
\[
 mt=mt-v+v 
\]
We multiply the partition by  $ \frac {1}{vt(mt-v)} $:
\[
 \frac{\frac{m}{v}}{mt-v}=\frac{1}{tv}+\frac{1}{t(mt-v)}.
\]
\item
 Consider the following identity:
\[
mt=mt
\]
Adding auxiliary variables $d$ and $l$:
\[
 mt=d+mt-(l+d)+l 
\]
We multiply the partition by  $ \frac {1}{rt(mt-(l+d))} $:
\[
 \frac{\frac {m}{r}}{mt-(l+d)}=\frac{1}{lt(mt-(l+d))}+\frac{1}{rt}+\frac{1}{dt(mt-(l+d))}.
\] 
\item 
Consider the following identity:
\[
kw(t+1)=kw(t+1)
\]
Adding auxiliary variables $w$ , $t$ and numbers $\{-1,+1\}$:
\[
 kw(t+1)=(kw-1)t+((k-1)w-1)+t+1+w
\]
We multiply the partition by $ \frac {1}{w(t+1)((kw-1)t+((k-1)w-1))} $:
\begin{align*}
&\frac{k}{((kw-1)t+((k-1)w-1))}=\frac{1}{w(t+1)}+\\
&\frac{1}{w((kw-1)t+((k-1)w-1))}+\frac{1}{(t+1)((kw-1)t+((k-1)w-1))}. 
\end{align*}
\item
Consider the following identity:
\[
kwt=kwt
\]
Adding auxiliary variables $w$ and $t$:
\[
 kwt=t(kw-1)+t-w+w
\]
We multiply the partition by $ \frac {1}{wt((kw-1)(t-w))} $:
\[ \frac{k}{(kw-1)(t-w)}=\frac{1}{w(t-w)}+\frac{1}{wt(kw-1)}+\frac{1}{t(kw-1)(t-w)}. \]

\item
Consider the following identity:
\[ 1=\frac{v-1}{v}+\frac{1}{v} \]
We multiply the partition by $\frac{1}{w}$:
\[\frac{1}{w}=\frac{v-1}{vw}+\frac{1}{vw} \]
Add auxiliary variable $ k$:
\[\frac{1}{w}+k=\frac{v-1}{vw}+\frac{1}{vw}+k \]
The equation can be rewritten as follows:
\[k=\frac{v-1}{vw}+\frac{1}{vw}+\frac{kw-1}{w} \]
We multiply the partition by $ \frac{1}{(kw-1)(v-1)}$:
\[
\frac{k}{(kw-1)(v-1)}=\frac{1}{vw(kw-1)}+\frac{1}{vw(kw-1)(v-1)}+\frac{1}{w(v-1)}.
\]

\item
Consider the following identity:
\[
A^{b}=A^{b}
\]
Adding auxiliary variables $A_{1}^{c}$ and $n$:
\[
A^{b}=A^{b}-A_{1}^{c}-n+A_{1}^{c}+n
\]
We multiply the partition by $ \frac {1}{A\left(A^{b}-A_{1}^{c}-n\right)n} $:
\begin{align*}
	&\frac{A^b}{p_{1}^{a_{1}}p_{2}^{a_{2}}\cdots p_{i}^{a_{i}}\left(A^{b}-A_{1}^{c}-n\right)n }&=\\&=\frac{1}{\left(p_{1}^{a_{1}}p_{2}^{a_{2}}\cdots p_{i}^{a_{i}}\right)n}+\frac{1}{p_{1}^{a_{1}+c}p_{2}^{a_{2}+c}\cdots p_{i}^{a_{i}+c}\left(A^{b}-A_{1}^{c}-n\right)n}&+\\&+\frac{1}{p_{1}^{a_{1}}p_{2}^{a_{2}}\cdots p_{i}^{a_{i}}\left(A^{b}-A_{1}^{c}-n\right)}.
\end{align*}

\item
Consider the following identity:
\[
krvs=krvs
\]
Adding auxiliary variables $rw$ and $ls$:
\[
krvs=krvs+ls+rw-ls-rw=ls+rw+(krv-l)s-rw
\]
We multiply the partition by $ \frac {1}{rvwls} $:
\[
\frac{k}{wl}=\frac{1}{rvw}+\frac{1}{vls}+\frac{1}{\frac{rvwls}{(krv-l)s-rw}}.
\]

\item
Consider the following identity:
\[
kmvs=kmvs
\]
Adding auxiliary variables $sb$ and $v$:
\[
kmvs=s(kmv\pm b)\mp v\mp sb\pm v
\]
We multiply the partition by $ \frac {1}{mvs(kmv\pm b)} $:
\[
\frac{k}{kmv\pm b}=\frac{1}{mv}\mp\frac{1}{ms(kmv\pm b)}+\frac{1}{\frac{mvs(kmv\pm b)}{\mp sb\pm v}}.
\]

The second proof for the equations is the same as the second proof of the first item.
\end{enumerate}
\end{proof}
Another equation that can be expressed from the tenth equation of theorem \ref{sara arab} is as follows:
\begin{cor}
	\[
	\forall k\in\mathbb{N}\quad,\quad \forall x\in\mathbb{N}\quad,\quad \exists y\in\mathbb{N} \quad,\quad kx>\mp y
	\]
	\[
	 \exists g\in\mathbb{N} \quad,\quad f+g=ry\quad,\quad f\mid xrg(kx\pm y)
	\]
\[
\frac{\frac{k}{g}}{kx\pm y}=\frac{1}{xg}\mp\frac{1}{xr(kx\pm y)}\mp\frac{1}{\frac{xrg(kx\pm y)}{f}}.
\]
\end{cor}
\begin{proof}
	\[
	\frac{1}{xg}\mp\frac{1}{xr(kx\pm y)}\mp\frac{1}{\frac{xrg(kx\pm y)}{f}}=\frac{r(kx\pm y)\mp f\mp g}{xrg(kx\pm y)}=\frac{\frac{k}{g}}{kx\pm y}.
	\]
\end{proof}
Another way to get well-defined equations is to use the distributive property that we use to get the next equations.

\begin{theorem}\label{sara 097}
The well-defined equation of the original conjecture is obtained using the distributive property:
\[
 \forall k\in \mathbb{N} \quad,\quad \exists z \in \mathbb{N} \quad,\quad \exists n,v \in \mathbb{N} \quad,\quad z\mid k \quad,\quad k\mid n+v   \]
\[
 \exists d \in \mathbb{N} \quad,\quad \exists \alpha , \beta \in \mathbb{N} \quad,\quad \alpha+\beta=\frac{k}{z}dv 
\]
\[
\exists r \in \mathbb{N} \quad,\quad\lcm(\alpha,\beta)=\frac{dn\left(\frac{n+v}{z}\right)}{r}
\]
\begin{equation} \label{sara 590}
\frac{k}{n}=\frac{1}{\frac{n+v}{k}}+\frac{1}{\frac{dn(n+v)}{z\alpha}}+\frac{1}{\frac{dn(n+v)}{z\beta}}.
\end{equation}
\end{theorem}
\begin{proof}
Using the following distribution, equation\eqref{sara 590} will be obtained:
\[
 kd\left(\frac{n+v}{z}\right)=\frac{k}{z}dn+\frac{k}{z}dv 
  \]
$k$ will be obtained in the form of the following fraction:
\[
 k=\frac{\frac{k}{z}dn+\frac{k}{z}dv}{d\left(\frac{n+v}{z}\right)}
\]
We multiply the partition by  $ \frac {1}{n} $:
\[
\frac{k}{n}=\frac{\frac{k}{z}dn+\frac{k}{z}dv}{dn\left(\frac{n+v}{z}\right)}
\]
We express it in the form of unit fractions:
\begin{equation} \label{sara 591}
\frac{k}{n}=\frac{1}{\frac{n+v}{k}}+\frac{1}{\frac{n(n+v)}{kv}}
\end{equation}
We will now expand the equation into three unit fractions:
\[
\frac{k}{n}=\frac{1}{\frac{n+v}{k}}+\frac{1}{\frac{dn(n+v)}{z\alpha}}+\frac{1}{\frac{dn(n+v)}{z\beta}}.
\]
\end{proof}
\begin{cor} Equation \eqref{sara 590} can be expressed in the following form:
\[
  \forall k \in \mathbb{N}  \quad,\quad \exists z \in \mathbb{N} \quad,\quad z\mid k \quad,\quad \exists d \in \mathbb{N}
\]
\[
 \exists \alpha,\beta \in \mathbb{N} \quad,\quad \frac{\alpha+\beta}{\frac{k}{z}d}=v \quad,\quad \frac{k\alpha\beta}{\gcd(k,\alpha,\beta)}=m \quad,\quad v,m \in \mathbb{N}
  \]
 \begin{equation}\label{sara mohammad}
  \frac{k}{mt-v}=\frac{1}{\frac{\alpha\beta t}{\gcd(k,\alpha,\beta)}}+\frac{1}{(mt-v)\left(\frac{dk\beta t}{z\gcd(k,\alpha,\beta)}\right)}+\frac{1}{(mt-v)\left(\frac{dk\alpha t}{z\gcd(k,\alpha,\beta)}\right)},
 \end{equation}
\[
\forall t\in
\begin{cases}
\quad\quad\quad\quad\mathbb{N}  & \mbox{if } m>v\\
\quad\quad\quad\mathbb{N}-\{1\}  & \mbox{if } m=v\\
\{t>s\mid v-m=s \:,\: t,s\in\mathbb{N}\}  & \mbox{if } m< v. 

\end{cases}
\]
\end{cor}
\begin{remark}
The values covered by $n$ in equation \eqref{sara mohammad} can be expressed as follows:
  \[
 n \equiv -\left(\frac{\alpha+\beta}{\frac{k}{z}d}\right)\ \left(\textrm{mod}\  \frac{k\alpha\beta}{\gcd(k,\alpha,\beta)}\right).
  \]
\end{remark}
\begin{example}

To make theorem \ref{sara 097} clearer, we present the following example.\\
We consider the equation in the following form:
\[
\frac{k}{n}=\frac{\frac{k}{z}dn+\frac{k}{z}dv}{dn\left(\frac{n+v}{z}\right)}
\]
Now for the rational number $\frac{k}{n}=\frac{4}{5569}$ , we express its two partition according to the presented equation as follows:
 \begin{enumerate}
\item
We will get the first partition by substituting $v=71,d=z=1$ values:
\[
\frac{4}{5569}=\frac{4\times5569+4\times71}{5569\times(5569+71)}=\frac{4\times5569+282+2}{5640\times5569}
\]
\[
\frac{4}{5569}=\frac{1}{1410}+\frac{1}{111380}+\frac{1}{15704580}.
\]
\item
We will get the second partition by substituting $v=7,d=3,z=2$ values:
\[
\frac{4}{5569}=\frac{6\times5569+6\times7}{3\times5569\times\left(\frac{5569+7}{2}\right)}=\frac{33414+41+1}{16707\times2788}
\]
\[
\frac{4}{5569}=\frac{1}{1394}+\frac{1}{1136076}+\frac{1}{46579116}.
\]
\end{enumerate}
Thus it can be stated:
\[
\frac{4}{5569}=\frac{1}{1410}+\frac{1}{111380}+\frac{1}{15704580}=\frac{1}{1394}+\frac{1}{1136076}+\frac{1}{46579116}.
\]
\end{example}
\begin{cor}\label{SARAAA}
	Equation\eqref{sara 591} for $k = 1$ can be expressed as follows:
	\[\forall n \in \mathbb{N} \quad,\quad \exists v\in \mathbb{N} \quad,\quad v\mid n(n+v)\]
	\[ \frac{1}{n}=\frac{n+v}{n(n+v)}=\frac{1}{n+v}+\frac{1}{\frac{n(n+v)}{v}}.             
	\]
\end{cor}
\begin{proof}
	Consider the following equation:
	\[ n+v=n+v
	\]\[ 1=\frac{n+v}{n+v}
	\]
	We multiply the partition by $ \frac {1}{n} $:
	\begin{equation}\label{ssara}
		\frac{1}{n}=\frac{n+v}{n(n+v)}=\frac{1}{n+v}+
		\frac{1}{\frac{n(n+v)}{v}}.             
	\end{equation}
\end{proof}
\begin{definition}
	$\tau(n) $ is the arithmetic function of counting the number of divisors the natural number $n$.\cite{SARAISMYLIFE}\\
We consider the decomposition of $w$ into the prime factors:
\[
w=p_{1}^{\alpha_{1}}p_{2}^{\alpha_{2}}\cdots p_{q}^{\alpha_{q}}
\]
Therefore:
\[
\tau(w)=\prod_{i=1}^{q} \left(\alpha_{i}+1\right).
\]
\end{definition}
\begin{definition}
		$A_{n}(v)$ is the arithmetic function of counting numbers such as $v$ that occur in equation \ref{ssara} .
\end{definition}
	\begin{definition}
			$P_{k}(n) $ is the arithmetic function of counting the number of partitions of the unit fraction $\frac{1}{n} $ in $k$-part.
	\end{definition}
\begin{remark}
	\[
	\forall n\in \mathbb{N} \quad,\quad \exists v\in \mathbb{N}  \quad,\quad v\mid n^{2}
	\]
	We consider the decomposition of $n$ into the prime factors:
	\[
	n=p_{1}^{\alpha_{1}}p_{2}^{\alpha_{2}}\cdots p_{q}^{\alpha_{q}}
	\]
	Hence it can be stated: 
	\[
	A_{n}(v)=\prod_{i=1}^{q} (2\alpha_{i}+1)=\tau(n^2)
	\]
	\[
	P_{2}(n)=\left \lceil \frac{A_{n}(v)}{2} \right \rceil=\left \lceil \frac{\tau(n^2)}{2} \right \rceil
	\]
	Now consider the following values:
	\[
	P_{3}(n)=\left \lceil \frac{A_{n+v}(v_{2})}{2} \right \rceil+\left \lceil \frac{A_{\frac{n(n+v)}{v}=s}(v_{3})}{2} \right \rceil
	\]
	\[
	P_{4}(n)=\left \lceil \frac{A_{n+v}(v_{2})}{2} \right \rceil\times\left \lceil \frac{A_{\frac{n(n+v)}{v}=s}(v_{3})}{2} \right \rceil
	\]
	\[
	P_{5}(n)=\left \lceil \frac{A_{n+v+v_{2}}(v_{4})}{2} \right \rceil+\left \lceil \frac{A_{\frac{(n+v)(n+v+v_{2})}{v_{2}}}(v_{6})}{2} \right \rceil+\left \lceil \frac{A_{s+v_{3}}(v_{5})}{2} \right \rceil+\left \lceil \frac{A_{\frac{s(s+v_{3})}{v_3}}(v_{7})}{2} \right \rceil
	\]
	\[
	P_{8}(n)=\left \lceil \frac{A_{n+v+v_{2}}(v_{4})}{2} \right \rceil\times\left \lceil \frac{A_{\frac{(n+v)(n+v+v_{2})}{v_{2}}}(v_{6})}{2} \right \rceil\times\left \lceil \frac{A_{s+v_{3}}(v_{5})}{2} \right \rceil\times\left \lceil \frac{A_{\frac{s(s+v_{3})}{v_3}}(v_{7})}{2} \right \rceil
	\]
	In this case, it is clear that in the general case for any $a$ belonging to natural numbers, we can say that:
	\[
	P_{2^{a-1}+1}(n)=\sum_{i=1}^{2^{a-1}}\left \lceil \frac{A_{c_i}(b_{i})}{2} \right \rceil
	\]
	\[
	P_{2^a}(n)=\prod_{i=1}^{2^{a-1}}\left \lceil \frac{A_{c_i}(b_{i})}{2} \right \rceil
	\]
\end{remark}
So far we have described methods for obtaining well-defined equations of equation \ref{saraismylife}, such as using identities and partitions. Now in the following theorem we present another method for obtaining well-defined equations.
\begin{theorem}\label{keysara}
	A well-defined equation of equation \ref{saraismylife} will be obtained by applying algebraic operations on the $\frac{k}{n}$ fraction as follows:
	\[
	\forall k\in\mathbb{N}\quad,\quad \exists n\in\mathbb{N} \quad,\quad  \forall z\in\mathbb{N} \quad,\quad\exists r\in\mathbb{N} \quad,\quad \exists \alpha_{1},\alpha_{2},\cdots,\alpha_{r}\in\mathbb{Q^+}
	\]
	\[
	z=\alpha_{1}+\alpha_{2}+\cdots+\alpha_{r}
	\]
	\[
	 \forall a\in\mathbb{N}\quad,\quad\exists l\in\mathbb{N}\quad,\quad\exists n_{1},n_{2},\cdots,n_{l},m_{1},m_{2},\cdots,m_{l}\in\mathbb{N} 
\]
\[
 \frac{a}{n}=\frac{m_{1}}{n_{1}}+\frac{m_{2}}{n_{2}}+\cdots+\frac{m_{l}}{n_{l}}
	\]
	\[
	km_{i}\mid an_{i} \quad,\quad 1\le i \le l-1 \quad,\quad i\in\mathbb{N}
	\]
\[an_{l}(\alpha_{1}+\alpha_{2}+\cdots+\alpha_{r})=km_{l}
  \lcm(\alpha_{1},\alpha_{2},\cdots,\alpha_{r})q \quad,\quad q\in\mathbb{N}
\]
	\[
	\frac{k}{n}=\frac{1}{\frac{an_{1}}{km_{1}}}+\frac{1}{\frac{an_{2}}{km_{2}}}+\cdots+\frac{1}{\frac{an_{l}(\alpha_{1}+\alpha_{2}+\cdots+\alpha_{r})}{km_{l}\alpha_{1}}}+\frac{1}{\frac{an_{l}(\alpha_{1}+\alpha_{2}+\cdots+\alpha_{r})}{km_{l}\alpha_{2}}}+\cdots+\frac{1}{\frac{an_{l}(\alpha_{1}+\alpha_{2}+\cdots+\alpha_{r})}{km_{l}\alpha_{r}}}.
	\]
\end{theorem}
\begin{proof}
	We consider fraction $\frac{k}{n}$ as follows:
	\[
	\frac{k}{a}\times\frac{a}{n} 
	\]
	We consider the fraction $\frac{a}{n}$ partition according to the theorem as follows:
	\[
	\frac{a}{n}=\frac{m_{1}}{n_{1}}+\frac{m_{2}}{n_{2}}+\cdots+\frac{m_{l}z}{n_{l}z}
	\]
	By inserting the partition of the number $z$, we rewrite the partition as follows:
	\[
     \frac{a}{n}=\frac{m_{1}}{n_{1}}+\frac{m_{2}}{n_{2}}+\cdots+\frac{m_{l}(\alpha_{1}+\alpha_{2}+\cdots+\alpha_{r})}{n_{l}(\alpha_{1}+\alpha_{2}+\cdots+\alpha_{r})}
	\]
In this case we will have:
	\[
	\frac{k}{a}\left(\frac{m_{1}}{n_{1}}+\frac{m_{2}}{n_{2}}+\cdots+\frac{m_{l}\alpha_{1}}{n_{l}(\alpha_{1}+\alpha_{2}+\cdots+\alpha_{r})}+\cdots+\frac{m_{l}\alpha_{r}}{n_{l}(\alpha_{1}+\alpha_{2}+\cdots+\alpha_{r})}\right)
	\]
	Hence it can be easily seen that the stated equation is correct.
\end{proof}
In the following corollary, we express the two equations that are expressed in \cite{E} and are a special case of theorem \ref{keysara} .
\begin{cor}
	\[
	\frac{4}{m}=\frac{1}{\frac{m+c}{4}}+\frac{1}{\frac{m+c}{4a}\frac{a+b}{c}m}+\frac{1}{\frac{m+c}{4b}\frac{a+b}{c}m}
	\]
	\[
	\frac{4}{m}=\frac{1}{\frac{mc+1}{4}m}+\frac{1}{\frac{mc+1}{4a}\frac{a+b}{c}}+\frac{1}{\frac{m+c}{4b}\frac{a+b}{c}}
	\]
\end{cor}
The following is the C++ programming code that is effective in obtaining various results from the equations of this corollary:
	\begin{lstlisting}[language=C]
	#include <iostream>
	
	using namespace std;
	
	int main()
	{
	int msumcdiv4, mmulcsum1div4;
	int checkifpworks;
	for(int m=1;m<=594000;m++){	
	checkifpworks=(m)%840;
	if(checkifpworks==1){		
	for(int c=1;c<100;c++){
	if ((c+1)%4==0 ){
	msumcdiv4=(m+c)/4;
	mmulcsum1div4=(m*c+1)/4;
	for (int a=1;a<msumcdiv4|| a<mmulcsum1div4;a++){
	for(int b=1;b<11;b++){
	if(((m+c)/4)%a==0 && ((m+c)/4)%b==0 && (a+b)%c==0 && ((a+b)/c)<100 || ((m*c+1)/4)%a==0 && ((m*c+1)/4)%b==0&&(a+b)%c==0 && ((a+b)/c)<100 ){
	cout<<"m is"<< m <<endl;
	cout<<"a is"<< a <<endl;				
	cout<<"b is"<< b <<endl;					
	cout<<"c is"<< c <<endl;
}	}	}}}}}}
\end{lstlisting}
There are several domains that can be considered; the domain that the program expresses is:
\[
m\in\{1\le m\le 594000|m \equiv 1 \pmod{840}\}\quad,\quad c\in\{1\le c<100|c\in\mathbb{N}\}
\]
\[
b\in\{1\le b<11|b\in\mathbb{N}\}
\]

In the following theorems we use different modes of positive integer partitions.

\begin{theorem}\label{sara 098}
The well-defined equation of the original conjecture is obtained by using the special mode of partitions of positive integers.
\[\forall k\in \mathbb{N}  \quad,\quad  \forall w \in \mathbb{N} \quad,\quad \exists c,d \in \mathbb{Q}^+\]
Such that:
\[
c=\left\{\frac{e}{f}\bigg{|} e,f\in\mathbb{N}\:\:,\:\:f\ne 0\right\} \quad,\quad d=\left\{\frac{u}{f}\bigg{|} u,f\in\mathbb{N}\:\:,\:\: f\ne 0\right\}
 \]

\[
\quad f\mid e+u \quad,\quad w=\lcm(e,u) 
\]
In this case, we can say:
\begin{equation} \label{sara a}
\frac{k}{kw-(c+d)}=\frac{1}{w}+\frac{1}{\frac{w}{c}(kw-(c+d))}+\frac{1}{\frac{w}{d}(kw-(c+d))}.
\end{equation}
\end{theorem}
\begin{proof}
Suppose that $v$ belongs to positive integers now consider the partition of positive integers in the following form: 
\[
 kw=v+c+d
\]
We multiply the partition by $ \frac{1}{vw} $:
\begin{equation} \label {sara 021}
\frac{k}{v}=\frac{1}{w}+\frac{1}{\frac{w}{c}v}+\frac{1}{\frac{w}{d}v}
\end{equation}
By inserting the relation $v$ the equation can be rewritten as follows:
\[
\frac{k}{kw-(c+d)}=\frac{1}{\frac{w}{c}(kw-(c+d))}+\frac{1}{w}+\frac{1}{\frac{w}{d}(kw-(c+d))}.
\]
\end{proof}
Another proof using the distribution property:
\begin{proof}
Equation\eqref{sara a} is also obtained with the property of distributability as follows:
\[
k(c+v+d)=kc+kv+kd 
\]
The distribution property can be rewritten as follows:
\[
k=\frac{kc+kv+kd}{(c+v+d)}
\]
We multiply the partition by $ \frac{1}{v} $ :
\[
\frac{k}{v}=\frac{kc+kv+kd}{v(c+v+d)}
\]
Now the equation is of the form:
\[
\frac{k}{v}=\frac{1}{\frac{c+v+d}{k}}+\frac{1}{\frac{v(c+v+d)}{kc}}+\frac{1}{\frac{v(c+v+d)}{kd}}.
\]
\end{proof}
\begin{remark}
We will state the reason for such a structure for equation\eqref{sara a} with respect to the following equation:
\[
\forall k \in \mathbb{N} \quad,\quad \exists x,y \in \mathbb{N} \quad,\quad \exists r,i \in \mathbb{N} \quad,\quad r=\frac{\lcm(x,y)}{x} \quad,\quad i=\frac{\lcm(x,y)}{y}
\]
\begin{equation}\label{sara o}
k=\frac{1}{x}+\frac{1}{y}+\frac{1}{\frac{\lcm(x,y)}{k\lcm(x,y)-r-i}}
\end{equation}
To get the partitions of rational numbers, we multiply the equation by $\frac{1}{n}$:
\[
\frac{k}{n}=\frac{1}{xn}+\frac{1}{yn}+\frac{1}{\frac{\lcm(x,y)n}{k\lcm(x,y)-r-i}}
\]
In this case, for the equation to hold, we consider the following equality:
\[
\exists n \in\mathbb{N} \quad,\quad n=k\lcm(x,y)-(r+i)
\]
And the equation will be obtained as follows:
\[
\frac{k}{n}=\frac{1}{xn}+\frac{1}{yn}+\frac{1}{\lcm(x,y)}
\]
And equation\eqref{sara o} is obtained using the following identity:
\[
k\lcm(x,y)=k\lcm(x,y)
\]
Adding auxiliary variables $r$ and $i$:
\[
k\lcm(x,y)=k\lcm(x,y)-r-i+r+i
\]
 We multiply the equation by $\frac{1}{\lcm(x,y)}$:
\[
k=\frac{1}{x}+\frac{1}{y}+\frac{1}{\frac{\lcm(x,y)}{k\lcm(x,y)-r-i}}.
\]
\end{remark}

\begin{cor}\label{sara -1}
Equation\eqref{sara a} can be expressed in another form:
\[
\forall k \in \mathbb{N} \quad,\quad \forall e,u \in \mathbb{N} \quad,\quad  \exists w,f \in \mathbb{N} \quad,\quad w=\lcm(e,u) \quad,\quad f \big{|} e+u
\]
\begin{equation}\label{sara rf}
\frac{k}{kwt-\left(\frac{e+u}{f}\right)}=\frac{1}{wt}+\frac{1}{\frac{wtf\left(kwt-\left(\frac{e+u}{f}\right)\right)}{e}}+\frac{1}{\frac{wtf\left(kwt-\left(\frac{e+u}{f}\right)\right)}{u}}
\end{equation}
\[
\forall t\in
\begin{cases}
\quad\quad\quad\quad\mathbb{N}  & \mbox{if } kw>\left(\frac{e+u}{f}\right)\\
\quad\quad\quad\mathbb{N}-\{1\}  & \mbox{if } kw=\left(\frac{e+u}{f}\right)\\
\left\{t>s\bigg{|} \left(\frac{e+u}{f}\right)-kw=s \:,\: t,s\in\mathbb{N}\right\}  & \mbox{if } kw<\left(\frac{e+u}{f}\right) .

\end{cases}
\]
\end{cor}
\begin{proof}
Given the proof of equation theorem\eqref{sara 098} consider the following equation:
\[
\frac{k}{v}=\frac{1}{a}+\frac{1}{\frac{afv}{e}}+\frac{1}{\frac{afv}{u}}
\]
We express the sum of unit fractions as follows:
\[
\frac{1}{a}+\frac{1}{\frac{afv}{e}}+\frac{1}{\frac{afv}{u}}=\frac{fv+e+u}{afv}
\]
In this case, this fraction must be valid in the following equation:
\[
\frac{fv+e+u}{afv}=\frac{ky}{vy}
\]
Now this equation can be rewritten in the following form:
\[
f=\frac{e+u}{ka-v}
\]
With placements $a=wt$ and $v=kwt-\left(\frac{e+u}{f}\right)$ equation\eqref{sara rf} is obtained correctly.
\end{proof}
\begin{cor}\label{sara sara}
Another form of expression of equation\eqref{sara a} is as follows:

\[ 
\forall k \in \mathbb{N} \quad,\quad \forall c,d \in \mathbb{N} \quad,\quad \exists w \in \mathbb{N}  \quad,\quad w=\lcm(c,d) 
\]
\begin{equation} \label{sara 6}
\frac{k}{(kw-c)t-d}=\frac{1}{\frac{w}{c}((kw-c)t-d)}+\frac{1}{wt}+\frac{1}{\frac{wt}{d}((kw-c)t-d)}
\end{equation}
\[
\forall t\in
\begin{cases}
\quad\quad\quad\quad\mathbb{N}  & \mbox{if } kw>c+d\\
\quad\quad\quad\mathbb{N}-\{1\}  & \mbox{if } kw=c+d\\
\left\{t>s\mid(c+d)-kw=s \:,\: t,s\in\mathbb{N}\right\}  & \mbox{if } kw<c+d 

\end{cases}
\]
According to the variable $t$, it is difficult to obtain $v$ as a positive integer, so the field of values $c$ and $d$ is limited to positive integers.
\end{cor}
\begin{proof}
Consider partitions of positive integers in the following form:
\[
 kwt=v+ct+d
\]
$v$ will be obtained in the following form:
\[
v=(kw-c)t-d
\]
Multiply the partition by $\frac{1}{vwt}$ to get the desired equation:
\[
\frac{k}{v}=\frac{1}{wt}+\frac{1}{\frac{wv}{c}}+\frac{1}{\frac{vwt}{d}}.
\]
\end{proof}
Equation\eqref{sara 6} is an equation for partitions of positive integers  with positive integer sizes.
\begin{remark} \label{sara spe}
The values covered by $n$ in equation\eqref{sara a} can be expressed as follows:
\[
v\equiv -(c+d) \pmod{k(\lcm(e,u))} \rightarrow \exists t \in \mathbb{N} \quad v=k(\lcm(e,u))t-(c+d).
\]
The values covered by $n$ in equation\eqref{sara 6} can be expressed as follows:
\[
v \equiv -d \pmod{k(\lcm(c,d)-c)} \rightarrow \exists t \in \mathbb{N} \quad v=k(\lcm(c,d)-c)t-d. 
\]
\end{remark}
\begin{remark}
In equation\eqref{sara 6}, for each $w$ ,$ (\tau(w))^2 $ well-defined equations can be obtained from this equation.
\end{remark}
\begin{cor}
We will now state equation\eqref{sara 021} in general terms:
\[ \forall k \in \mathbb{N} \quad,\quad \forall w \in \mathbb{N} \quad,\quad \exists s \in\mathbb{N} \quad,\quad  \exists c_{1},c_{2},\cdots,c_{s}\in\mathbb{N}
\]
\[ w=\lcm\left(c_{1},c_{2},\cdots,c_{s}\right) \quad,\quad v=kw-\left(\sum_{i=1}^s c_i \right) \quad,\quad v\in \mathbb{N}   \]
\begin{equation} 
\frac{k}{v}=\frac{1}{w}+\frac{1}{\frac{vw}{c_1}}+\frac{1}{\frac{vw}{c_2}}+\cdots+\frac{1}{\frac{vw}{c_s}}. \label{sara 91}
\end{equation}
\end{cor}
\begin{proof}
\begin{align*}
&\frac{1}{w}+\frac{1}{\frac{vw}{c_1}}+\frac{1}{\frac{vw}{c_2}}+\cdots+\frac{1}{\frac{vw}{c_s}}&=\\
&=\frac{v+\sum_{i=1}^s c_i}{vw}&=\\
&=\frac{k}{v}.
\end{align*}
\end{proof}
\begin{cor}
A well-defined equation in the form of equation\eqref{sara 021}:
\[ \forall k \in \mathbb{N}  \quad,\quad \exists m,l,e \in \mathbb{N} \quad,\quad \forall d \in \mathbb{N}   
\]

\[
 k\mid m+2l \quad,\quad e+1\bigg{|} \frac{m+2l}{k} \quad,\quad  e+2l\bigg{|} \frac{m+2l}{k} 
\]

\[
M=k\left(\frac{m+2l}{k} m^{d-1}\right)t-\left( (2lt+e)m^{d-1}+(e+1)\right)
\]
\begin{align*}
\frac{k}{M}=\frac{k}{m^{d-1}(mt-e)-(e+1)}=&\frac{1}{\left(\frac{m+2l}{k}\right)m^{d-1}}+\frac{1}{\frac{\left(\frac{m+2l}{k}\right)\left(m^{d-1}(mt-e)-(e+1)\right)}{2l+e}}\\
&+\frac{1}{\frac{\left(\frac{m+2l}{k} m^{d-1}\right)\left(m^{d-1}(mt-e)-(e+1)\right)}{e+1}}
\end{align*}
\[
\forall t\in
\begin{cases}
\quad\quad\quad\quad\mathbb{N}  & \mbox{if } m^{d}-1>e\left(m^{d-1}+1\right)\\
\left\{m^{d}t>e\left(m^{d-1}+1\right)+1\mid t\in\mathbb{N}\right\}  & \mbox{if } m^{d}-1\not>e\left(m^{d-1}+1\right). 

\end{cases}
\]
\end{cor}
\begin{proof}
\begin{align*}
&\frac{1}{\left(\frac{m+2l}{k}\right)m^{d-1}}+\frac{1}{\frac{\left(\frac{m+2l}{k}\right)\left(m^{d-1}(m-e)-(e+1)\right)}{2l+e}}+\frac{1}{\frac{\left(\frac{m+2l}{k} m^{d-1}\right)\left(m^{d-1}(m-e)-(e+1)\right)}{e+1}}&=\\
&=\frac{m^{d-1}(m-e)-(e+1)+(2l+e)m^{d-1}+(e+1)}{\left(\frac{m+2l}{k}\right)\left(m^{d-1}(mt-e)-(e+1)\right)}&=\\
&=\frac{m^{d-1}(2l+m)}{\left(\frac{m+2l}{k}\right)\left(m^{d-1}(mt-e)-(e+1)\right)}&=\\
&=\frac{k}{m^{d-1}(m-e)-(e+1)}.
\end{align*}
\end{proof}
Therefore, according to the expressed structure, different equations can be expressed.
\begin{remark}
The set of values covered by the equation is shown according to the following congruence:
\[
 M\equiv-\left((2lt+e)m^{d-1}+(e+1)\right)\ \left(\textrm {mod}\ k\left(\frac{m+2l}{k} m^{d-1}\right)\right).
\]
\end{remark}

\begin{theorem} \label{sara 013}
Consider the following equation which is obtained from a particular case of positive integer partitions:
 \[
\forall k \in \mathbb{N} \quad,\quad \forall c \in \mathbb{N} \quad,\quad \exists v \in \mathbb{N} \quad,\quad \exists q \in \mathbb{N}-\{1\}  
\]
\[
c\mid q-1 \quad,\quad  \frac{qv+c}{k}=b \quad,\quad \frac{b}{q-1}=w \quad,\quad b,w \in \mathbb{N}
\]
\[
\frac{k}{v}=\frac{1}{w}+\frac{1}{b}+\frac{1}{\frac{vb}{c}}.
\]
\end{theorem}
\begin{proof}
Consider partitions of positive integers in the following form:
\[
kwq=vq+c+kw
\]
The partition can be written as:
\[
	  	 kw(q-1)=vq+c
\]
The equation using the auxiliary variable is as follows:
\[
	  	 kw(q-1)=vq+c+v-v=v+v(q-1)+c 
\]
We multiply the partition by $ \frac {1}{vw(q-1)} $:
\[
\frac{k}{v}=\frac{1}{w(q-1)}+\frac{1}{w}+\frac{1}{\frac{vw(q-1)}{c}}
\]
Now it will be shown that the equation has answers or in other words there are partitions in the expressed form.
\[
\forall k \in \mathbb{N} \quad,\quad \forall c \in \mathbb{N} \quad,\quad \exists q \in \mathbb{N}-\{1\}  \quad,\quad c\mid q-1  
\]
Consider the following linear congruence equation:
\[
qv \equiv -c\ (\textrm{mod}\ k(q-1)) \quad,\quad d=\gcd(k(q-1),q) \quad,\quad d \mid c
\]
$v$ is obtained in the following form:
\[
v=v_{1}+\frac{k(q-1)}{d}t \quad,\quad t\in\mathbb{N}
\]
$v_{1}$ an answer that holds in the linear congruence equation.\\
However it can be stated:
\[
\exists w\in \mathbb{N}
\]
\[
\frac{qv+c}{k}=w(q-1)
\]
The equation with these conditions is expressed in the following form:
\[
\frac{k}{v}=\frac{1}{w(q-1)}+\frac{1}{w}+\frac{1}{\frac{vw(q-1)}{c}}.
\]
\end{proof}
\begin{theorem} \label{sara 530}
\begin{enumerate}
The following equations are well-defined equations that are obtained by adding a natural number to the answer value of a linear congruence equation:
\item
\[
\forall k \in \mathbb{N} 
\]
\[
\forall \omega \in \mathbb{N}  \quad,\quad \exists v,\eta \in \mathbb{N} \quad,\quad k\eta  \equiv -1\pmod{v}\quad,\quad r=\eta+\omega 
\]
\[
 \forall m\in\{kt-2 \mid t \in \mathbb{N}\} \quad,\quad \exists c \in \mathbb{N} \quad,\quad c\mid m 
\]
\[
 \frac{k}{v}=\frac{1}{r}+\frac{1}{rm}+\frac{1}{\frac{vrm}{c}}
\]
\item
\[ 
\forall k \in \mathbb{N} 
\]
\[
\forall \omega \in \mathbb{N}  \quad,\quad \exists v,\eta \in \mathbb{N} \quad,\quad k\eta  \equiv -1\pmod{v}\quad,\quad r=\eta+\omega 
\]
\[
 \forall m\in \mathbb{N} \quad,\quad \exists c \in \mathbb{N} \quad,\quad c\mid m
\]
\[
\frac{k}{v}=\frac{1}{r}+\frac{1}{vm}+\frac{1}{\frac{vrm}{c}}.
\]
\end{enumerate}
\end{theorem}
\begin{proof}
\begin{enumerate}
\item 
Consider partitions of positive integers in the following form:
\[
 krm=vm+v+c
 \]
We multiply the partition by $ \frac {1}{vrm} $:
\[
 \frac{k}{v}=\frac{1}{r}+\frac{1}{rm}+\frac{1}{\frac{vrm}{c}}
\]
The equation is expressed in the following well-defined form:
\[
\forall k\in \mathbb{N}-\{1,2\} \quad,\quad \forall \omega \in \mathbb{N} \quad,\quad \forall z\in \mathbb{N}_{(k-2)\omega}
\]
\begin{eqnarray*}
     \left\lbrace \begin{array}{lc}
   v=k(k\omega-1)t-(2(k\omega-1)+(kz-(k-1)))\\
 \eta=(k\omega-1)t-(z+2\omega-1)\\
    c=kz-(k-1)\\
     m=kt-2 \\
    \end{array}\right.
\end{eqnarray*}
\begin{align*}
&\frac{k}{k(k\omega-1)t-(2(k\omega-1)+(kz-(k-1)))}&=\\
&=\frac{1}{(k\omega-1)t-(z+\omega-1)}+\frac{1}{((k\omega-1)t-(z+\omega-1))(kt-2)}&+\\
&+\frac{1}{\frac{((k\omega-1)t-(z+\omega-1))(kt-2)(k(k\omega-1)t-(2(k\omega-1)+(kz-(k-1))))}{kz-(k-1)}}
\end{align*}
\begin{eqnarray*}
  \forall t \in   \left\lbrace \begin{array}{lc}
   \quad\quad\quad\quad\quad\quad\mathbb{N} \quad &if \quad z=1 \\
    \left\{b(z-1)+\frac{b+2}{k} \big{|}  b=ky-2 , y\in\mathbb{N}\right\} \quad &if \quad z\ne1
    \end{array}\right.
\end{eqnarray*}
Now the proof of this well-defined equation is expressed in the form of the sum of the unit fractions and reaching the other side of the equation:
\begin{align*}
&\frac{1}{(k\omega-1)t-(z+\omega-1)}+\frac{1}{((k\omega-1)t-(z+\omega-1))(kt-2)}&+\\
&+\frac{1}{\frac{((k\omega-1)t-(z+\omega-1))(kt-2)(k(k\omega-1)t-(2(k\omega-1)+(kz-(k-1))))}{kz-(k-1)}}&=\\
&=\frac{1}{\frac{(k\omega-1)t-(z+\omega-1))(kt-2)(k(k\omega-1)t-(2(k\omega-1)+(kz-(k-1))))}{(kt-2)(k(k\omega-1)t-(2(k\omega-1)+(kz-(k-1)))+k\omega-1)}}&=\\
&=\frac{1}{\frac{(k\omega-1)t-(z+\omega-1))(kt-2)(k(k\omega-1)t-(2(k\omega-1)+(kz-(k-1))))}{(kt-2)(k((k\omega-1)t-(z+\omega-1)))}}&=\\
&=\frac{k}{k(k\omega-1)t-(2(k\omega-1)+(kz-(k-1)))}
\end{align*}
To express the equation as a set of unit fractions, we consider the following equation.\\
Suppose $b$ belongs to the set of natural numbers, in which case we can express:
\[
\frac{kt-2}{kz-k+1}=b 
\]
Therefore, $t$ will be obtained as follows:
\[
 t=b(z-1)+\frac{b+2}{k}
\]
In order for t to belong to the set of natural numbers, it must $k\mid b+2$.\\
For $k=1$, $v$ will be obtained as follows:
\[
v=(\omega-1)t-(2(\omega-1)+z)
\]
The values of $t$ are obtained as follows:
\[
t>2+\frac{z}{\omega-1}
\]
Therefore, we choose the values of $z$ and $\omega$ in such a way that $\omega-1\mid z$.\\
For $k=2$, $v$ will be obtained as follows:
\[
v=2(2\omega-1)t-(2(2\omega-1)+2z-1)
\]
The values of $t$ are obtained as follows:
\[
t>z+1
\]
For each $(2\omega-1)$ consecutive set, is $t>z+1$; Because if we consider the following inequality:
\[
2(2\omega-1)<2z+4\omega-3<4(2\omega-1)
\]
In this case z will be obtained:
\[
\frac{1}{2}<z<\frac{4\omega-1}{2}
\]
Since $z$ belongs to natural numbers, we will have:
\[
\left\lfloor{\frac{4\omega-1}{2}}\right\rfloor=2\omega-1.
\]

\item
Consider partitions of positive integers in the following form:
\[
 krw=vw+r+c
\]
We multiply the partition by $ \frac {1}{vrw} $:
\[
\frac{k}{v}=\frac{1}{r}+\frac{1}{vm}+\frac{1}{\frac{vrm}{c}}
\]
The equation is expressed in the following well-defined form:
\[
\forall k \in \mathbb{N}-\{1\} 
\]
\[
\forall \omega\in
\begin{cases}
\quad\quad\mathbb{N}  & \mbox{if } k \in \mathbb{N}-\{1,2,3\}\\
\{\omega>1 \mid \omega\in\mathbb{N}\}  & \mbox{if }  k=3\\
\{\omega>2 \mid \omega\in\mathbb{N}\}  & \mbox{if }  k=2
\end{cases}
\]
\[
 \exists l,c \in \mathbb{N} \quad ,\quad  \omega(k-1)>l+c
\]
\begin{eqnarray*}
    \left\lbrace \begin{array}{lc}
   v=&k(k\omega-l)t-k(\omega+c)+l\\
   m=&t \\
   \eta =&(k\omega-l)t-(c+\omega)
    \end{array}\right.
\end{eqnarray*}
\begin{align*}
&\frac{k}{k(k\omega-l)t-k(\omega+c)+l}&=\\
&=\frac{1}{(k\omega-l)t-c}+\frac{1}{(k(k\omega-l)t-k(\omega+c)+l)t}&+\\
&+\frac{1}{\frac{((k\omega-l)t-c)(k(k\omega-l)t-k(\omega+c)+l)t}{c}}
\end{align*}
\[
\forall t\in
\begin{cases}
\quad\quad\mathbb{N}  & \mbox{if } c=1\\
\{cx \mid x\in\mathbb{N}\}  & \mbox{if }  c\ne1

\end{cases}
\]
To prove the inequality $\omega(k-1)>l+c$, $t=1$ we consider: 
\[
k(k\omega-l)-k(\omega+c)+l=k(\omega(k-1)-(l+c))+l
\]
Now for the product to belong to natural numbers, you have to consider:
\[
\omega(k-1)-(l+c)>0
\]
Therefore:
\[
\omega(k-1)>(l+c)
\]
Now the proof of this well-defined equation is expressed in the form of the sum of the unit fractions and reaching the other side of the equation:
\begin{align*}
&\frac{1}{(k\omega-l)t-c}+\frac{1}{(k(k\omega-l)t-k(\omega+c)+l)t}&+\\
&+\frac{c}{((k\omega-l)t-c)(k(k\omega-l)t-k(\omega+c)+l)t}&=\\
&=\frac{(k(k\omega-l)t-k(\omega+c)+l)t+(k\omega-l)t-c+c}{((k\omega-l)t-c)(k(k\omega-l)t-k(\omega+c)+l)t}&=\\
&=\frac{t(k(k\omega-l)t-c)}{((k\omega-l)t-c)(k(k\omega-l)t-k(\omega+c)+l)t}&=\\
&=\frac{k}{k(k\omega-l)t-k(\omega+c)+l}
\end{align*}
\end{enumerate}
We now express the equation for $k=1$ :
\[
(\omega-l)t-(\omega+c)+l=\omega(t-1)-l(t-1)-c>0
\]
Therefore it can be stated:
\[
\forall \omega \in \mathbb{N}-\{1,2\} \quad,\quad \exists l,c \in \mathbb{N} \quad,\quad \omega>l+c
\]
And the values of $t$ is as follows:
\[
\forall t \in \mathbb{N}-\{1\}.
\]
\end{proof}
Let us now give an example to make the proof of the first part of theorem \ref{sara 530} clearer.
\begin{example}
	To clarify the matter of $k=2$, we will state it for $\omega\in\{1,2,3\}$:
	\begin{enumerate}
		\item 
		$\omega=1$
		\begin{enumerate}
			\item
			$ z=1 \rightarrow \{2t-3\mid t\ge2\}$
		\end{enumerate}
		\item
		$\omega=2$
		\begin{enumerate}
			\item
			$ z=1 \rightarrow \{6t-7\mid t\ge2\}$
			\item
			$ z=2 \rightarrow\{6t-9 \mid t\ge2\}$
			\item
			$ z=3 \rightarrow\{6t-11 \mid t\ge2\}$
		\end{enumerate}
		\item
		$\omega=3$
		\begin{enumerate}
			\item
			$ z=1 \rightarrow \{10t-11\mid t\ge2\}$
			\item
			$ z=2 \rightarrow\{10t-13 \mid t\ge2\}$
			\item
			$ z=3 \rightarrow\{10t-15 \mid t\ge2\}$
			\item
			$ z=4 \rightarrow\{10t-17 \mid t\ge2\}$
			\item
			$ z=5 \rightarrow\{10t-19 \mid t\ge2\}$.
		\end{enumerate}
	\end{enumerate}
\end{example}

\begin{theorem} \label{sssaraaa}
A well-defined equation for the original conjecture:
\[
\forall k \in \mathbb{N}\quad,\quad\forall v \in \mathbb{N} \quad,\quad \exists b,w \in \mathbb{N} 
\]
\[
 b\mid k   \quad,\quad v\mid w \quad,\quad kw>b \quad,\quad w >v
\]
\[
\frac{\frac{k}{b}}{(kw-b)t-kv}=\frac{1}{b\left(wt-v\right)}+\frac{1}{w((kw-b)t-kv)}+\frac{1}{\frac{w\left(wt-v\right)\left((kw-b)t-kv\right)}{v}}
\]
\[
\forall t\in
\begin{cases}
\quad\quad\quad\mathbb{N}  & \mbox{if } k(w-v)>b\\
\{(kw-b)t>kv \mid t\in\mathbb{N}\}  & \mbox{if }  k(w-v)<b.

\end{cases}
\]

\end{theorem}
\begin{proof}
\begin{align*}
&\frac{1}{b\left(wt-v\right)}+\frac{1}{w((kw-b)t-kv)}+\frac{1}{\frac{w\left(wt-v\right)((kw-b)t-kv)}{v}}&=\\
&=\frac{w((kw-b)t-kv)+b\left(wt-v\right)+bv}{wb\left(wt-v\right)\left((kw-b)t-kv\right)}&=\\
&=\frac{w\left(k\left(wt-v\right)\right)}{wb\left(wt-v\right)\left((kw-b)t-kv\right)}&=\\
&=\frac{\frac{k}{b}}{(kw-b)t-kv}.
\end{align*}
\end{proof}
In this way many well-defined equations can be expressed for equation\eqref{saraismylife} using the identity $a=a$, the property of distributive ,the auxiliary variable and different modes of partitions of positive integers.\\
The original conjecture has two states that the unit fractions are different from each other or not. The methods expressed are independent of these modes in this way any mode that can be considered can be obtained.\\

\section{Some features of well-defined equations}

\begin{theorem}\label{sara....f}
The well-defined equation for each $k$ in the form $a^b$ can be expressed as follows:
\[
\forall a,b \in \mathbb{N} \quad,\quad\exists c_{1},c_{2},c_{3},\cdots,c_{s} \in \mathbb{N}
\]
\[
 \exists v \in \mathbb{N} \quad,\quad v>b+c_{i} \quad,\quad 1\le i\le s
\]
\begin{align*}\label{sara...fat}
&\frac{a^b}{a^{v}t-\left(\sum_{i=1}^s a^{c_{i}}+1\right)}=\frac{1}{a^{v-b}t}+\frac{1}{a^{v-(b+c_{1})}t(a^{v}t-\left(\sum_{i=1}^s
a^{c_{i}}+1\right))}&+\\
&+\frac{1}{a^{v-(b+c_{2})}t(a^{v}t-\left(\sum_{i=1}^s
a^{c_{i}}+1\right))}+\cdots+\frac{1}{a^{v-(b+c_{s})}t(a^{v}t-\left(\sum_{i=1}^s
a^{c_{i}}+1\right))}&+\\
&+\frac{1}{a^{v-b}t(a^{v}t-\left(\sum_{i=1}^s
a^{c_{i}}+1\right))}
\end{align*}
\[
\forall t\in
\begin{cases}
\quad\quad\quad\mathbb{N}  & \mbox{if } a^{v}>\left(\sum_{i=1}^s a^{c_{i}}+1\right) \\
\left\{a^{v}t>\left(\sum_{i=1}^s a^{c_{i}}+1\right)  \bigg{|} t\in\mathbb{N}\right\}  & \mbox{if }  a^{v}<\left(\sum_{i=1}^s a^{c_{i}}+1\right). 

\end{cases}
\]
\end{theorem}
\begin{proof}
Correctness of the equation is proved by the method of total unit fractions and reaching the other side of the equation in the following form:

\begin{align*}
&\frac{1}{a^{v-b}t}+\frac{1}{a^{v-(b+c_{1})}t(a^{v}t-(\sum_{i=1}^s 
a^{c_{i}}+1))}+\frac{1}{a^{v-(b+c_{2})}t(a^{v}t-(\sum_{i=1}^s
a^{c_{i}}+1))}&+\\
&+\cdots+\frac{1}{a^{v-(b+c_{s})}t(a^{v}t-(\sum_{i=1}^s
a^{c_{i}}+1))}+\frac{1}{a^{v-b}t(a^{v}t-(\sum_{i=1}^s
a^{c_{i}}+1))}&=\\
&=\frac{a^{v}t-(\sum_{i=1}^s
a^{c_{i}}+1)+a^{c_{1}}+a^{c_{2}}+\cdots+a^{c_{s}}+1}{a^{v-b}t(a^{v}t-(\sum_{i=1}^s
a^{c_{i}}+1))}&=\\
&=\frac{a^{v}t}{a^{v-b}t(a^{v}t-(\sum_{i=1}^s
a^{c_{i}}+1))}=\frac{a^b}{a^{v}t-(\sum_{i=1}^s
a^{c_{i}}+1)}.
\end{align*}
\end{proof}

\subsection{Well-defined equations based on the partition of numbers in the form $kb-1$}\label{sara 500}
\begin{cor}
Consider numbers in the form $kb-1$ as the sum of two rational numbers:
\[\forall k \in \mathbb{N} \quad,\quad \exists b \in \mathbb{N} \quad,\quad kb>1 \quad,\quad \exists \lambda_1,\lambda_2,\gamma \in \mathbb{N}\]
\[
\frac{k}{\frac{\lambda_1\lambda_2}{\gcd(\lambda_1,\lambda_2)}}=\frac{1}{\frac{\left(\frac{\lambda_1\lambda_2}{\gcd(\lambda_1,\lambda_2)}\right)\gamma b}{\lambda_1}}+\frac{1}{\frac{\left(\frac{\lambda_1\lambda_2}{\gcd(\lambda_1,\lambda_2)}\right)\gamma b}{\lambda_2}}+\frac{1}{\left(\frac{\lambda_1\lambda_2}{\gcd(\lambda_1,\lambda_2)}\right)b}.
\]
\end{cor}
\begin{proof}
Consider generalizing Two-part partitions of positive integers to the following form:
\[
kb-1=\frac{\lambda_1}{\gamma}+\frac{\lambda_2}{\gamma} \rightarrow kb=\frac{\lambda_1}{\gamma}+\frac{\lambda_2}{\gamma}+1 
\]
According to theorem \ref{sara 10} now multiply the identity by  $ \frac{1}{\left(\frac{\lambda_1\lambda_2}{\gcd(\lambda_1,\lambda_2)}\right)b}$:
\[
\frac{k}{\frac{\lambda_1\lambda_2}{\gcd(\lambda_1,\lambda_2)}}=\frac{1}{\frac{\left(\frac{\lambda_1\lambda_2}{\gcd(\lambda_1,\lambda_2)}\right)\gamma b}{\lambda_1}}+\frac{1}{\frac{\left(\frac{\lambda_1\lambda_2}{\gcd(\lambda_1,\lambda_2)}\right)\gamma b}{\lambda_2}}+\frac{1}{\left(\frac{\lambda_1\lambda_2}{\gcd(\lambda_1,\lambda_2)}\right)b}.
\]
\end{proof}
\begin{cor} \label{sara 300}
Now consider the numbers in the form $kb-1$ as a multiplication of two positive integers:
\[ \forall k \in \mathbb{N} \quad,\quad \exists b \in \mathbb{N} \quad,\quad kb>1 \quad,\quad \exists  q,v \in \mathbb{N} \quad,\quad qv=kb-1\quad,\quad q-1\mid b  \]
\begin{equation} \label{sara 20 }
 \frac{k}{v}=\frac{1}{\frac{b}{q-1}}+\frac{1}{b}+\frac{1}{bv}
\end{equation}
Equation\eqref{sara 20 } is a special mode of equation\eqref{sara 013}.

\end{cor}
\begin{proof}
According to the condition of the corollary can be expressed:
\[
qv=kb-1 \rightarrow qv+1=kb \rightarrow \frac{qv+1}{b}=k  
\]
We multiply the equation by $ \frac{1}{v}$:
\[ \frac{qv+1}{bv}=\frac{k}{v} \]
Now we express it using auxiliary variable as a well-defined equation:
\begin{align*}
\frac{qv+1}{bv}&=\frac{qv+1+v-v}{bv}=\frac{v(q-1)+v+1}{bv}\\&=\frac{v(q-1)}{bv}+\frac{v+1}{bv}=\frac{1}{\frac{b}{q-1}}+\frac{1}{b}+\frac{1}{bv}=\frac{k}{v}.
\end{align*}
\end{proof}

\subsection{Different forms of a well-defined equation}\label{sara..fat}
\begin{cor}
Another form of corollary\eqref{sara 300} is as follows:
\[  \forall k \in \mathbb{N}  \quad,\quad \forall v \in \mathbb{O}-\{1\} \quad,\quad \exists a,b,u,u{'} \in \mathbb{N} \]

\[v-1\mid u\quad,\quad v-1\mid u{'} \quad,\quad av=ku  \quad,\quad  1-bv=ku{'} \quad  \] 

\[ 
\frac{k}{at-b}=\frac{1}{\frac{(at-b)v+1}{k(v-1)}}+\frac{1}{\frac{(at-b)v+1}{k}}+\frac{1}{(at-b)\left(\frac{(at-b)v+1}{k}\right)}
\]
$ t $ with positive leading coefficients.
\end{cor}
\begin{proof}
According to corollary \ref{sara 300}:
\[
\frac{(at-b)v+1}{k}
\]
It can now be rewritten as follows:
\[
\frac{atv-bv+1}{k}=\frac{av}{k}t+\frac{1-bv}{k}
\]
According to the constraints of the corollary it can be stated:
\[
ut+u{'}
\]
According to corollary \ref{sara 300}:
\[
\frac{ut+u{'}}{v-1}=\frac{u}{v-1}t+\frac{u{'}}{v-1}
\]
According to the constraints of the corollary $ v-1\mid u$ and  $v-1\mid u{'} $.
\end{proof}
\begin{cor}
According to subsection \ref{sara 500} and corollary \ref{sara 300} a well-defined equation for numbers in the form $kb-1$ is expressed:
\[ \forall k \in \mathbb{N}\]
\begin{equation} \label {sara 71}
\frac{k}{kb-1}=\frac{1}{b}+\frac{1}{kb^2}+\frac{1}{(kb-1)kb^2} 
\end{equation}
\[
\forall b\in
\begin{cases}
\mathbb{N}  & \mbox{if } k\in\mathbb{N}-\{1\} \\
\mathbb{N}-\{1\}  & \mbox{if } k=1. 

\end{cases}
\]

\end{cor}
\begin{proof}
According to corollary \ref{sara 300}:
\[
(kb-1)(kb+1)+1=kb^2
\]
The equation is expressed as follows:
\[ \frac{k}{kb-1}=\frac{1}{\frac{(kb-1)(kb+1)+1}{k(kb+1-1)}}+\frac{1}{\frac{(kb-1)(kb+1)+1}{k}}+\frac{1}{(kb-1)\left(\frac{(kb-1)(kb+1)+1}{k} \right)} \] 
After simplification, equation \eqref{sara 71} will be obtained.
\end{proof}
Another proof according to the auxiliary variable:
\begin{proof}
Consider the following identity:
\[
kb=kb
\]
Adding auxiliary variable $\frac{1}{kb}$ and numbers $\{-1,+1\}$:
\[
 kb=1-\frac{1}{kb}+kb-1+\frac{1}{kb}
\]
Now multiply the identity by $\frac{1}{(kb-1)b}$:
\[
\frac{k}{kb-1}=\frac{1}{kb^2}+\frac{1}{b}+\frac{1}{kb^2(kb-1)}. 
\]
\end{proof}
Well-defined equations can be described in different ways, now in this section we will give an example of well-defined Erd\"{o}s-Straus equations in different ways.\\
\begin{cor} \label{sara 91}
According to corollary \ref{sara 300} we express a well-defined equation of Erd\"{o}s-Straus as follows.\\
We express the equation in two logical ways:
\begin{enumerate}
\item
\[ \forall q \in \mathbb{O}-\{1\}\quad,\quad \exists v \in \mathbb{O} \quad,\quad \exists  b\in \mathbb{E} \quad,\quad qv+1=4b \quad,\quad q-1 \mid b \]
\[
 \frac{4}{v}=\frac{1}{\frac{b}{q-1}}+\frac{1}{b}+\frac{1}{bv} 
\]
\item
\[ \forall w \in \mathbb{N} \quad,\quad \exists  q \in \mathbb{O}-\{1\}\quad,\quad \exists  b \in \mathbb{E} \quad ,\quad b=(q-1)w
\]
\[
 \exists v \in \mathbb{O} \quad,\quad qv+1=4b    
\]
\end{enumerate}
\begin{equation} \label{sara 21}
 \frac{4}{v}=\frac{1}{\frac{b}{q-1}}+\frac{1}{b}+\frac{1}{bv}
\end{equation}
\end{cor}
A new way to prove equation\eqref{sara 21}:
\begin{proof}
\begin{enumerate}
\item
Equation\eqref{sara 21} can be expressed in the form of well-defined equations\eqref{sara 0021} according to lemma \ref{sara 430} the corollary \ref{sara 91} can be proved as follows:
\begin{align*}
\frac{4}{8t-1}&=\frac{1}{\frac{(8t-1)(8t+1)+1}{4}}+\frac{1}{\frac{(8t-1)(8t+1)+1}{32t}}+\frac{1}{(8t-1)\left(\frac{(8t-1)(8t+1)+1}{4}\right)}\\
\frac{4}{8t-5}&=\frac{1}{\frac{(8t-5)(8t-3)+1}{4}}+\frac{1}{\frac{(8t-5)(8t-3)+1}{4(8t-4)}}+\frac{1}{(8t-5)\left(\frac{(8t-5)(8t-3)+1}{4}\right)}\\
 \frac{4}{24t-7}&=\frac{1}{\frac{(24t-7)(8t-1)+1}{4}}+\frac{1}{\frac{(24t-7)(8t-1)+1}{4(8t-2)}}+\frac{1}{(24t-7)\left(\frac{(24t-7)(8t-1)+1}{4}\right)}\\
\frac{4}{24t-7}&=\frac{1}{\frac{(24t-7)(7)+1}{4}}+\frac{1}{\frac{(24t-7)(7)+1}{24}}+\frac{1}{(24t-7)\left(\frac{(24t-7)(7)+1}{4}\right)}\\
\frac{4}{24t-19}&=\frac{1}{\frac{(24t-19)(8t-5)+1}{4}}+\frac{1}{\frac{(24t-19)(8t-5)+1}{4(8t-6)}}+\frac{1}{(24t-19)\left(\frac{(24t-19)(8t-5)+1}{4}\right)} \\
\frac{4}{8t-3}&=\frac{1}{\frac{(8t-3)(3)+1}{4}}+\frac{1}{\frac{(8t-3)(3)+1}{8}}+\frac{1}{(8t-3)\left(\frac{(8t-3)(3)+1}{4}\right)}
\end{align*}
\[
\forall t \in \mathbb{N}.
\]
\item
We now prove equation\eqref{sara 21} using the logical method in the second case.
Consider constraints:
\[
qv+1=4b \quad,\quad q-1 \mid b
\]
Considering $b=w(q-1)$, can be expressed:
\[
qv+1=4w(q-1)
\]
The following two relations can now be concluded:
\begin{equation} \label{sara 39}
\frac{1+4w}{4w-v}=q 
\end{equation}
\[
\frac{qv+1}{4(q-1)}=w 
\]
According to relation \eqref{sara 39}for each $w \in \mathbb{N}$ there is at least one pair of $(q,v)$ which is true in equation \eqref{sara 21}.\\
Because, the result of the fraction is equal to $q$ and according to the set of numbers $q$ for each $w \in \mathbb{N} $ there is $\tau(4w+1)-1$  pair$(q, v)$ and the subtraction of the number one is because there is no number one in the set of numbers $q$.
\end{enumerate}
\end{proof}
We now present another proof in the form of the generalization of the schinzel method for equation\eqref{sara 21}:
\begin{cor} \label{sara 731}
Equation \eqref{sara 21} can be expressed as a generalization of the schinsel method as follows:
\[ \forall \omega \in \mathbb{N}\] 
\begin{enumerate}
\item\label{sara dv}
\begin{align*} 
&\frac{4}{8(2\omega-1)t-(4\omega-1)}=\frac{1}{\frac{\left(8(2\omega-1)t-(4\omega-1)\right)(4\omega-1)+1}{4}}+\\&+\frac{1}{\frac{(8(2\omega-1)t-(4\omega-1))(4\omega-1)+1}{4(4\omega-2)}}+\frac{1}{\frac{\left((8(2\omega-1)t-(4\omega-1))(4\omega-1)+1\right)(8(2\omega-1)t-(4\omega-1))}{4}}
\end{align*}
\item\label{sara dg}
\begin{align*}
&\frac{4}{8(2\omega)t-(12\omega+1)}=\frac{1}{\frac{(8(2\omega)t-(12\omega+1))(4\omega+1)+1}{4}}+\\
&+\frac{1}{\frac{(8(2\omega)t-(12\omega+1))(4\omega+1)+1}{16\omega}}+\frac{1}{\frac{(8(2\omega)t-(12\omega+1))\left((8(2\omega)t-(12\omega+1))(4\omega+1)+1\right)}{4}} 
\end{align*}
In both equations $ t \in \mathbb{N} $.
\end{enumerate}
\begin{proof}
\begin{enumerate}
\item Using the sum of the unit fractions and reaching the other side of the equation, the equation is proved:
\begin{align*}
&\frac{4}{\left(8(2\omega-1)t-(4\omega-1)\right)(4\omega-1)+1}&+\\&+\frac{4(4\omega-2)}{(8(2\omega-1)t-(4\omega-1))(4\omega-1)+1}&+\\
&+\frac{4}{\left((8(2\omega-1)t-(4\omega-1))(4\omega-1)+1\right)(8(2\omega-1)t-(4\omega-1))}&=\\
&=\frac{4(8(2\omega-1)t-(4\omega-1)+(4w-2)(8(2\omega-1)t-(4\omega-1))+1)}{\left((8(2\omega-1)t-(4\omega-1))(4\omega-1)+1\right)(8(2\omega-1)t-(4\omega-1))}&=\\
&=\frac{4((8(2\omega-1)t-(4\omega-1))(4\omega-1)+1)}{\left((8(2\omega-1)t-(4\omega-1))(4\omega-1)+1\right)(8(2\omega-1)t-(4\omega-1))}&=\\
&=\frac{4}{8(2\omega-1)t-(4\omega-1)}.
\end{align*}
\item We prove the second equation similar to the first equation:
\begin{align*}
&\frac{4}{(8(2\omega)t-(12\omega+1))(4\omega+1)+1}&+\\&+ \frac{16\omega}{(8(2\omega)t-(12\omega+1))(4\omega+1)+1}&+\\
&+\frac{4}{(8(2\omega)t-(12\omega+1))\left((8(2\omega)t-(12\omega+1))(4\omega+1)+1\right)}&=\\
&=\frac{4(8(2\omega)t-(12\omega+1)+4\omega(8(2\omega)t-(12\omega+1))+1)}{(8(2\omega)t-(12\omega+1))\left((8(2\omega)t-(12\omega+1))(4\omega+1)+1\right)}&=\\
&=\frac{4((8(2\omega)t-(12\omega+1))(4\omega+1)+1)}{(8(2\omega)t-(12\omega+1))\left((8(2\omega)t-(12\omega+1))(4\omega+1)+1\right)}&=\\
&=\frac{4}{8(2\omega)t-(12\omega+1)}.
\end{align*}
\end{enumerate}
\end{proof}
This method shows that equation\eqref{sara 21} for two family of sets of rational numbers expresses two distinct equations.\\
\begin{remark}
The numbers that each equation of the first proof method (Schinzel method) covers are expressed as a set while the numbers that each equation of the second proof method (Schinzel generalization method) covers are expressed as a family of sets.\\\\
The scope of definition of equation\eqref{sara dv}:
\[\left\{8(2\omega-1)t-(4\omega-1)\mid\omega \in \mathbb{N}\right\}_{t =1 }^ {\infty}=\left\{\left\{8t-3\right\},\left\{24t-7\right\},\cdots\right\}_{t =1 }^ {\infty}\]
The scope of definition of equation\eqref{sara dg}:
\[\left\{8(2\omega)t-(12\omega+1)\mid\omega \in \mathbb{N}\right\}_{t =1 }^ {\infty}=\left\{\left\{16t-13\right\},\left\{32t-25\right\},\cdots\right\}_{t =1 }^ {\infty}\]
\end{remark}
\end{cor}
Thus an equation can be expressed in different ways.
\begin{cor}
The generalizations of equations\eqref{sara dv} and \eqref{sara dg} are as follows:
\begin{enumerate}
\item
\[
\forall k \in \mathbb{N} 
\]
$
\forall \omega\in
\begin{cases}
\quad\quad\quad\mathbb{N}  & \mbox{if } k \in \mathbb{N}-\{1,2\} \\
\left\{\omega>1\mid \omega\in\mathbb{N}\right\}  & \mbox{if }  k=2\\
\left\{\omega>2\mid \omega\in\mathbb{N}\right\}  & \mbox{if }  k=1

\end{cases}
$
\begin{align*}
&\frac{k}{\left(k^{2}\omega-2k\right)t-(k\omega-1)}&=\\
&=\frac{1}{\frac{\left(\left(k^{2}\omega-2k\right)t-(k\omega-1)\right)(k\omega-1)+1}{k}}+\frac{1}{\frac{\left(\left(k^{2}\omega-2k\right)t-(k\omega-1)\right)(k\omega-1)+1}{k(k\omega-2)}}&+\\
&+\frac{1}{\frac{\left(\left(k^{2}\omega-2k\right)t-(k\omega-1)\right)(\left(\left(k^{2}\omega-2k\right)t-(k\omega-1))(k\omega-1)+1\right)}{k}} 
\end{align*}
$
\forall t\in
\begin{cases}
\quad\quad\quad\mathbb{N}  & \mbox{if } k \in \mathbb{N}-\{1,2\} \\
\left\{t>\frac{2\omega-1}{4(\omega-1)}\bigg{|} t\in\mathbb{N}\right\}  & \mbox{if }  k=2\\
\left\{t>\frac{\omega-1}{\omega-2}\mid t\in\mathbb{N}\right\}  & \mbox{if }  k=1

\end{cases}
$
\item
\[
\forall k \in \mathbb{N}  \quad,\quad \forall \omega \in \mathbb{N} 
\]
\begin{align*}
&\frac{k}{k^{2}\omega t-((k-1)k\omega+1)}&=\\
&=\frac{1}{\frac{\left(k^{2}\omega t-((k-1)k\omega+1)\right)(k\omega+1)+1}{k}}+\frac{1}{\frac{\left(k^{2}\omega t-((k-1)k\omega+1)\right)(k\omega+1)+1}{k^2\omega}}&+\\
&+\frac{1}{\frac{\left(k^{2}\omega t-((k-1)k\omega+1)\right)\left(\left(k^{2}\omega t-((k-1)k\omega+1)\right)(k\omega+1)+1)\right)}{k}} 
\end{align*}
\end{enumerate}
$
\forall t\in
\begin{cases}
\quad\quad\quad\mathbb{N}  & \mbox{if } k \in \mathbb{N}-\{1\} \\
\left\{wt>1\mid t\in\mathbb{N}\right\}  & \mbox{if }  k=1.

\end{cases}
$

\end{cor}
\begin{proof}
\begin{enumerate}
\item Using the sum of the unit fractions and reaching the other side of the equation, the equation is proved:
\begin{align*}
&\frac{k}{\left(\left(k^{2}\omega-2k\right)t-(k\omega-1)\right)(k\omega-1)+1}&+\\&+\frac{k(k\omega-2)}{\left(\left(k^{2}\omega-2k\right)t-(k\omega-1)\right)(k\omega-1)+1}&+\\
&+\frac{k}{\left(\left(k^{2}\omega-2k\right)t-(k\omega-1)\right)\left(\left(k^{2}\omega-2k\right)t-(k\omega-1))(k\omega-1)+1\right)}&=\\
&=\frac{k\left(\left(\left(k^{2}\omega-2k\right)t-(k\omega-1)\right)+(k\omega-2)\left(\left(k^{2}\omega-2k\right)t-(k\omega-1)\right)+1\right)}{\left(\left(k^{2}\omega-2k\right)t-(k\omega-1)\right)\left(\left(\left(k^{2}\omega-2k\right)t-(k\omega-1)\right)(k\omega-1)+1\right)}&=\\
&=\frac{k\left(\left(\left(k^{2}\omega-2k\right)t-(k\omega-1)\right)(k\omega-1)+1\right)}{\left(\left(k^{2}\omega-2k\right)t-(k\omega-1)\right)\left(\left(\left(k^{2}\omega-2k\right)t-(k\omega-1)\right)(k\omega-1)+1\right)}&=\\
&=\frac{k}{\left(k^{2}\omega-2k\right)t-(k\omega-1)}
\end{align*}

We will now prove the following divisibility:

\[
k(k\omega-2) \Bigg{|} \left(2k\left(\frac{k}{2}\omega-1\right)t-(k\omega-1)\right)(k\omega-1)+1
\]

We prove divisibility in the form of the following fraction:

\begin{align*}
&\frac{\left(2k\left(\frac{k}{2}\omega-1\right)t-(k\omega-1)\right)(k\omega-1)+1}{k(k\omega-2)}&=\\
&=\frac{\left(2k\left(\frac{k}{2}\omega-1\right)(k\omega-1)\right)t-(k\omega-1)^{2}+1}{k(k\omega-2)}&=\\
&=\frac{\left(2k\left(\frac{k}{2}\omega-1\right)(k\omega-1)\right)t-2\omega\left(\frac{k}{2}\omega-1\right)}{2k\left(\frac{k}{2}\omega-1\right)}&=\\
&=(k\omega-1)t-\omega
\end{align*}
\item We prove the second equation similar to the first equation:
\begin{align*}
&\frac{k}{\left(k^{2}\omega t-((k-1)k\omega+1)\right)(k\omega+1)+1}&+\\&+\frac{k^2\omega}{\left(k^{2}\omega t-((k-1)k\omega+1)\right)(k\omega+1)+1}&+\\
&+\frac{k}{\left(k^{2}\omega t-((k-1)k\omega+1)\right)\left(\left(k^{2}\omega t-((k-1)k\omega+1)\right)(k\omega+1)+1)\right)}&=\\
&=\frac{k\left(\left(k^{2}\omega t-((k-1)k\omega+1)\right)+k\omega\left(k^{2}\omega t-((k-1)k\omega+1)\right)+1\right)}{\left(k^{2}\omega t-((k-1)k\omega+1)\right)\left(\left(k^{2}\omega t-((k-1)k\omega+1)\right)(k\omega+1)+1)\right)}&=\\
&=\frac{k\left(\left(k^{2}\omega t-((k-1)k\omega+1)\right)(k\omega+1)+1\right)}{\left(k^{2}\omega t-((k-1)k\omega+1)\right)\left(\left(k^{2}\omega t-((k-1)k\omega+1)\right)(k\omega+1)+1)\right)}&=\\
&=\frac{k}{k^{2}\omega t-((k-1)k\omega+1)}
\end{align*}

We will now prove the following divisibility:

\[
k^2\omega \Bigg{|} \left(k^{2}\omega t-((k-1)k\omega+1)\right)(k\omega+1)+1
\]

We prove divisibility in the form of the following fraction:

\begin{align*}
&\frac{\left(k^{2}\omega t-((k-1)k\omega+1)\right)(k\omega+1)+1}{k^2\omega}&=\\
&=\frac{\left(k^{2}\omega(k\omega+1)\right)t-((k-1)k\omega+1)(k\omega+1)+1}{k^2\omega}&=\\
&=\frac{\left(k^{2}\omega(k\omega+1)\right)t-((k-1)k\omega)k\omega-(k-1)k\omega-k\omega-1+1}{k^2\omega}&=\\
&=(k\omega+1)t-((k-1)k\omega+1).
\end{align*}
\end{enumerate}
\end{proof}
\begin{lemma} \label {sara 430}
Suppose $ q \in \mathbb{O}-\{1\}, $ then $v \in \mathbb{O}$ exists such that 
$ \frac{qv+1}{4} \in \mathbb{E}$ and $ q-1\big{|}\frac{qv+1}{4} $.
\end{lemma}
\begin{proof}
For $ q \in \mathbb{O}-\{1\} $ when $ \frac{qv+1}{4} $ is a even and $ qv+1 \equiv 0 \pmod{8} $, this lemma is established.\\
So the question becomes the following:
\[
\exists v \in \mathbb{O} \quad s.t \quad qv \equiv 7 \pmod{8} 
\]
 There are four modes for $q $ :
\[ q\equiv +1 \pmod{8}  \quad,\quad q\equiv -1 \pmod{8} \]
\[ q\equiv -3 \pmod{8}  \quad,\quad q\equiv -5 \pmod{8} \]
For congruence can be expressed:
\begin{align*}
 if \quad q&\equiv  +1 \pmod{8}\quad then\quad v=\{8t-1\mid t \in \mathbb{N} \} \\
 if \quad q&\equiv -1 \pmod{8}\quad then\quad v=\{24t-7\mid t \in \mathbb{N} \} \\
 if \quad q&\equiv -3 \pmod{8}\quad then\quad v=\{8t-5\mid t \in \mathbb{N} \} \\
 if  \quad q&\equiv-5\pmod{8}\quad then\quad v=\{24-19\mid t \in \mathbb{N} \}. 
\end{align*}
\end{proof}
\begin{cor}
	The domain determined by the Erd\"{o}s-Straus equation is a set of natural numbers, while the domain of well-defined equations can be generalized.\\
	Now we express a well-defined equation whose domain is a subset of real numbers.
	\begin{equation}\label{saaaraaa}
	\frac{4}{8x-1}=\frac{1}{2x}+\frac{1}{16x^2}+\frac{1}{16x^2(8x-1)} \quad,\quad \forall x\in\mathbb{R}^{+}-\left\{0,\frac{1}{8}\right\}
	\end{equation}
	Using variable variation, a well-defined equation can be obtained from another well-defined equation.\\
	By changing the variable to the following form in equation \ref{saaaraaa}, equation \ref{saaaaraaaa} can be obtained:
	\[
	x=t-\frac{1}{2}
	\]
	\begin{equation}\label{saaaaraaaa}
	\frac{4}{8t-5}=\frac{1}{2t-1}+\frac{1}{4(2t-1)^2}+\frac{1}{4(2t-1)^2(8t-5)} \quad,\quad \forall t \in \mathbb{N}
	\end{equation}
Another clearer form that exists is as follows:
	\[
	\frac{4}{40t-23}=\frac{4}{20(2t-1)-3}
	\]
	\[
	\frac{4}{20(2t-1)-3}=\frac{1}{5(2t-1)}+\frac{1}{2(2t-1)(40t-23)}+\frac{1}{10(2t-1)(40t-23)}
	\]
	\[
	\forall t \in \mathbb{N}.
	\]
\end{cor}

\subsection{Well-defined equations triangle}\label{sara.fat}
Two special states of theorem \ref{sara 530} for the Erd\"{o}s-Straus equation, the numbers covered by these two special states of well-defined equations form two triangles that have the following properties.\\
\begin{cor} \label{sara 732}
Two special modes of the equations of theorem \ref{sara 530} are expressed as follows:

\begin{enumerate}
\item
\[
\forall \omega \in \mathbb{N} \quad,\quad \forall m\in \left\{2(2t-1)\mid t \in \mathbb{N}\right\} \quad,\quad \exists v,\eta \in \mathbb{N} 
\]
\[
 4\eta  \equiv -1\pmod{\nu} \quad,\quad r=\eta+\omega 
\]
\begin{equation}
\frac{4}{v}=\frac{1}{r}+\frac{1}{rm}+\frac{1}{vrm}
\end{equation}
\item
\[
\forall \omega \in \mathbb{N} \quad,\quad \forall t \in \mathbb{N} \quad,\quad \exists v,\eta \in \mathbb{N} \quad,\quad 4\eta \equiv -1 \pmod{v}
 \quad,\quad  r=\eta+\omega 
\]
\begin{equation} 
\frac{4}{v}=\frac{1}{r}+\frac{1}{vt}+\frac{1}{vrt}.
\end{equation}
\end{enumerate}
\end{cor}
\begin {proof}
To prove the two equations, we express them in the form of a generalization of the Schinzel method:
\[
\forall \omega \in \mathbb{N}
\]
\begin{enumerate}
\item 
\begin{align*}
&\frac{4}{4(4\omega-1)t-(2(4\omega-1)+1)}&=\\
&=\frac{1}{(4\omega-1)t-\omega}+\frac{1}{((4\omega-1)t-\omega)(4t-2)}&+\\
&+\frac{1}{4(4\omega-1)t-(2(4\omega-1)+1)((4\omega-1)t-\omega)(4t-2)}
\end{align*}

\item
\begin{align*}
&\frac{4}{4(4\omega-1)t-(4\omega+3)}&=\\
&=\frac{1}{(4\omega-1)t-1}+\frac{1}{(4(4\omega-1)t-(4\omega+3))t}&+\\
&+\frac{1}{(4(4\omega-1)t-(4\omega+3))((4\omega-1)t-1)t}   
\end{align*}

In both equations $ t \in \mathbb{N}$ .
\end{enumerate}
\end{proof}

\begin{cor} \label{sara 222}
The generalization of the well-defined equations of corollary \ref{sara 732} is expressed as follows:
\begin{enumerate}
\item
\[
\forall k \in \mathbb{N} 
\]
\[
\forall \omega\in
\begin{cases}
\quad\quad\quad\mathbb{N}  & \mbox{if } k \in \mathbb{N}-\{1\} \\
\left\{\omega>1\mid \omega\in\mathbb{N}\right\}  & \mbox{if }  k=1

\end{cases}
\]
\begin{align*}
&\frac{k}{k(k\omega-1)t-(2(k\omega-1)+1)}&=\\
&=\frac{1}{(k\omega-1)t-\omega}+\frac{1}{((k\omega-1)t-\omega)(kt-2)}&+\\
&+\frac{1}{(k(k\omega-1)t-(2(k\omega-1)+1))((k\omega-1)t-\omega)(kt-2)}   
\end{align*}
\[
\forall t\in
\begin{cases}
\quad\quad\quad\mathbb{N}  & \mbox{if } k \in \mathbb{N}-\{1,2\} \\
\left\{t>\frac{4\omega-1}{4\omega-2}\mid t\in\mathbb{N}\right\}  & \mbox{if }  k=2\\
\left\{t>\frac{2\omega-1}{\omega-1}\mid t\in\mathbb{N}\right\}  & \mbox{if }  k=1

\end{cases}
\]
\item
\[
\forall k \in \mathbb{N} \quad,\quad \forall a \in \mathbb{N}  \quad,\quad \forall \omega \in \mathbb{N}  \quad,\quad  \exists b,c \in \mathbb{N} \quad,\quad c>b \quad,\quad k=|b-c| 
\]
\begin{align*}
&\frac{k}{k(a\omega+b)t-(a\omega+c)}&=\\
&=\frac{1}{(a\omega+b)t-1}+\frac{1}{(k(a\omega+b)t-(a\omega+c))t}&+\\
&+\frac{1}{(k(a\omega+b)t-(a\omega+c))((a\omega+b)t-1) t}   
\end{align*}
\[
\forall t\in\left\{t>\frac{a\omega+c}{k(a\omega+b)}\bigg{|} t\in\mathbb{N}\right\}.
\]
\end{enumerate}
\end{cor}
\begin{proof}
Using the sum of the unit fractions and reaching the other side of the equation, the equation is proved:
\begin{enumerate}
\item
\begin{align*}
&\frac{1}{(a\omega+b)t-1}+\frac{1}{(k(a\omega+b)t-(a\omega+c))t}&+\\&+\frac{1}{(k(a\omega+b)t-(a\omega+c))((a\omega+b)t-1)t}&=\\
&=\frac{(kt-2)(k(k\omega-1)t-(2(k\omega-1)+1))+k(k\omega-1)t-(2(k\omega-1)+1)+1}{(k(k\omega-1)t-(2(k\omega-1)+1))((k\omega-1)t-\omega)(kt-2)}&=\\
&=\frac{(kt-2)(k((k\omega-1)t-\omega))}{(k(k\omega-1)t-(2(k\omega-1)+1))((k\omega-1)t-\omega)(kt-2)}&=\\
&=\frac{k}{k(k\omega-1)t-(2(k\omega-1)+1}
\end{align*}
\item
\begin{align*}
&\frac{1}{(a\omega+b)t-1}+\frac{1}{(k(a\omega+b)t-(a\omega+c))t}&+\\
&+\frac{1}{(k(a\omega+b)t-(a\omega+c))((a\omega+b)t-1)t}&=\\
&=\frac{t(k(a\omega+b)t-(a\omega+c))+(a\omega+b)t-1+1}{(k(a\omega+b)t-(a\omega+c))((a\omega+b)t-1)t}&=\\
&=\frac{t(k((a\omega+b)t-1))}{(k(a\omega+b)t-(a\omega+c))((a\omega+b)t-1)t}&=\\
&=\frac{k}{k(a\omega+b)t-(a\omega+c)}.
\end{align*}
\end{enumerate}
\end{proof}
\begin{definition} \label{sara f}
$M(k,a,d;r,t) $ represents the pattern of a set of symmetric numerical triangles.
\begin{equation} \label{sara 10000}
M(k,a,d;r,t):=(ra+(r-1)k)t-(d+(r-1)k)
\end{equation}
\[
\forall k \in \mathbb{N}  \quad,\quad \forall r,t \in \mathbb{N}  \quad,\quad \exists a,d \in \mathbb{N}  \quad ,\quad ra > d   
\]
Using polynomial\eqref{sara 10000} a set of numerical triangles is obtained, each of which is represented by the following form:
\[
\forall k \in \mathbb{N}  \quad,\quad \exists a,d \in \mathbb{N}   \quad ,\quad ra > d 
\]
\[
m(r,t):=(ra+(r-1)k)t-(d+(r-1)k)
\]
\[
\forall r,t \in \mathbb{N}
\]

The triangle of this relation is as follows:
\begin{table}[ht] \label{sara 777}
\caption{m(r,t)=(ra+(r-1)k)t-(d+(r-1)k)}\label{eqtable}
\renewcommand\arraystretch{1.5}
\noindent\[
\begin{array}{lcccccccr}
{}&{}&{}&{}&{(1,1)}&{}&{}&{}&{}\\

{}&{}&{}&{(1,2)}&{}&{(2,1)}&{}&{}&{}\\

{}&{}&{(1,3)}&{}&{(2,2)}&{}&{(3,1)}&{}&{}\\

{}&{(1,4)}&{}&{(2,3)}&{}&{(3,2)}&{}&{(4,1)}&{}\\

{(1,5)}&{}&{(2,4)}&{}&{(3,3)}&{}&{(4,2)}&{}&{(5,1)}\\
\end{array}
\]
\end{table}
\end{definition}
\begin{cor} \label{sara j}
Special state of the second equation of corollary \ref{sara 222}, the numbers covered by that special state have a numerical triangle structure which are expressed as follows:
\[\forall k \in \mathbb{N}\]
\[
\forall r\in
\begin{cases}
\mathbb{N}-\{1\} & \mbox{if } k=1\\
\mathbb{N} & \mbox{if } k\in\mathbb{N}- \{1\}\\
\end{cases}
\]

\begin{align*}
&\frac{k}{k(kr-1)t-(kr+(k-1))}&=\\
&=\frac{1}{(kr-1)t-1}+\frac{1}{(k(kr-1)t-(kr+(k-1)))t}&+\\
&+\frac{1}{(k(kr-1)t-(kr+(k-1)))((kr-1)t-1)t}   
\end{align*}
The value of $t$:
\[
\forall t\in
\begin{cases}
\mathbb{N} & \mbox{if } k\in\mathbb{N}-\{1,2\}\\
\mathbb{N} & \mbox{if } k=2\quad,\quad r\in\mathbb{N}-\{1\}\\
\mathbb{N}-\{1\} & \mbox{if } k=r+1=2\\
\left\{t>\frac{r}{r-1}\mid t \in\mathbb{N}\right\} & \mbox{if } k =1
\end{cases}
\]
In this well-defined equation for each $ k\in \mathbb{N}$ the values covered by the equation have a numerical triangle pattern as follows:
\[
m(r,t)=k(kr-1)t-(kr+(k-1))
\]
\[ \forall r,t \in \mathbb{N}. \]
\end{cor}
\begin{proof}
We prove the well-defined equation expressed in the form of the sum of unit fractions:

\begin{align*}
&\frac{1}{(kr-1)t-1}+\frac{1}{(k(kr-1)t-(kr+(k-1)))t}&+\\
&+\frac{1}{(k(kr-1)t-(kr+(k-1)))((kr-1)t-1)t}&=\\
&=\frac{t(k(kr-1)t-(kr+(k-1)))+(kr-1)t-1+1}{(k(kr-1)t-(kr+(k-1)))((kr-1)t-1)t}&=\\
&=\frac{tk((kr-1)t-1)}{(k(kr-1)t-(kr+(k-1)))((kr-1)t-1)t}&=\\
&=\frac{k}{k(kr-1)t-(kr+(k-1))}
\end{align*}

According to definition \ref{sara f} by placing $a=k(k-1)$ and $d=2k-1$ in $m(r,t)$ the values covered by the well-defined equation are obtained.
Its well-defined equation is expressed as follows:

\[
m(r,t)=(r(k(k-1))+(r-1)k)t-(2k-1+(r-1)k)
\]

By performing algebraic operations the following form will be obtained:

\[
m(r,t)=k(kr-1)t-(kr+(k-1)).
\]
\end{proof}

\begin{remark}
For each $k \in\mathbb{N}-\{1,2\}$ all numbers in numerical triangles belong to the set of positive integers.
\end{remark}

\begin{cor}
According to definition \ref{sara f} and corollary \ref{sara j} the numbers belonging to the second well-defined equation of corollary \ref{sara 732} have a triangular structure.

\end{cor}
\begin{proof}
We will present the family of numbers in the form of a table \ref{sssss} and a triangular structure \ref{sara 3A4S}:\\
\[
\mathbf{L}=\left\{4(4r-1)t-(4r+3)\mid r \in \mathbb{N}\right\}_{t =1 }^ {\infty}
\]
\[
\mathbf{L}=\left\{\left\{12t-7\right\},\left\{28t-11\right\},\left\{44t-15\right\},\cdots\right\}_{t =1 }^ {\infty}
\]

\[
\mathbf{Q}=\left\{4(4t-1)r-(4t+3)\mid t \in \mathbb{N}\right\}_{r=1}^ {\infty}
\]
\[
\mathbf{Q}=\left\{\left\{12r-7\right\},\left\{28r-11\right\},\left\{44r-15\right\},\cdots\right\}_{r =1 }^ {\infty}.
\]
\end{proof}
\begin{table}[] 
\caption{}\label{eqtable}
\renewcommand\arraystretch{1.5}
\noindent\[
\begin{array}{c|cccccc} \label{sssss}
&1&2&3&4&...&{t}\\
\hline
{1}&{5}&{17}&{29}&{41}&{...}&{12t-7}\\

{2}&{17}&{45}&{73}&{101}&{...}&{28t-11}\\

{3}&{29}&{73}&{117}&{161}&{...}&{44t-15}\\

{4}&{41}&{101}&{161}&{221}&{...}&{60t-19}\\

{\vdots	}&{\vdots	}&{\vdots	}&{\vdots	}&{\vdots	}&{\vdots	}&{\vdots}\\

{r}&{12r-7}&{28r-11}&{44r-15}&{60r-19}&{...}&{\mathbf{Q}= \mathbf{L}}\\

\end{array}
\]
\end{table}
\begin{table}[] 
\caption{m(r,t)=4(4r-1)t-(4r+3)}\label{eqtable}
\renewcommand\arraystretch{1.5}
\noindent\[
\begin{array}{lcccccccr} \label{sara 3A4S}
{}&{}&{}&{}&{5}&{}&{}&{}&{}\\

{}&{}&{}&{17}&{}&{17}&{}&{}&{}\\

{}&{}&{29}&{}&{45}&{}&{29}&{}&{}\\

{}&{41}&{}&{73}&{}&{73}&{}&{41}&{}\\

{53}&{}&{101}&{}&{117}&{}&{101}&{}&{53}\\
\end{array}
\]
\end{table}
So we can express well-defined equations whose numbers are in the form of numerically symmetric triangles, now we are expressing equations whose numbers are in the form of numerically asymmetric triangles,any numerically asymmetric triangles cover two equations.

\begin{definition}
$A(k,a,d;r,t)$ represents the pattern of a set of numerically asymmetric triangles.

\begin{equation} \label{sara 10001}
A(k,a,d;r,t):=(ra+(r-1)k)t-(rd+(r-1))
\end{equation}
\[
\forall k \in \mathbb{N}  \quad,\quad \forall r,t \in \mathbb{N}  \quad,\quad \exists a,d \in \mathbb{N}  \quad ,\quad a > d   
\]

Using polynomial\eqref{sara 10000}, a set of numerically asymmetric triangles is obtained each of which is represented by the following form:

\[
\forall k \in \mathbb{N}  \quad,\quad \exists a,d \in \mathbb{N} \quad ,\quad a > d    
\]
\[
a(r,t):=(ra+(r-1)k)t-(rd+(r-1))
\]
\[
\forall r,t \in \mathbb{N}
\]

The triangle of this relation is as follows:

\begin{table}[ht] \label{sara 777}
\caption{a(r,t)=(ra+(r-1)k)t-(rd+(r-1))}\label{eqtable}
\renewcommand\arraystretch{1.5}
\noindent\[
\begin{array}{lcccccccr}
{}&{}&{}&{}&{(1,1)}&{}&{}&{}&{}\\

{}&{}&{}&{(1,2)}&{}&{(2,1)}&{}&{}&{}\\

{}&{}&{(1,3)}&{}&{(2,2)}&{}&{(3,1)}&{}&{}\\

{}&{(1,4)}&{}&{(2,3)}&{}&{(3,2)}&{}&{(4,1)}&{}\\

{(1,5)}&{}&{(2,4)}&{}&{(3,3)}&{}&{(4,2)}&{}&{(5,1)}\\
\end{array}
\]
\end{table}
\end{definition}
\begin{remark}
According to the definition of numerical asymmetric triangles, considering the special state of them in the following form, the first equation of the corollary \ref{sara 222} covers such numerical asymmetric triangles:
\[
a=k(k-1) \quad,\quad d=2k-1
\]
\[
a(r,t)=(rk(k-1)+(r-1)k)t-(r(2k-1)+(r-1)).
\]
\end{remark}
\begin{cor} \label{sara 555}
The first well-defined equation in corollary \ref{sara 731} and the first well-defined equation in corollary \ref{sara 732} these two equations cover only one family of numbers.
\end{cor}
\begin{proof}
Table\eqref{sara 1181} shows this corollary. As can be seen from the table, consider the numbers in column form, the family of numbers is shown in the first equation\eqref{sara 731}, and the numbers in the surface form of the family of numbers are shown in the first equation\eqref{sara 732}. \\
\\
The family of numbers is the first well-defined equation in the corollary\eqref {sara 731}:
\[
\mathbf{M}=\left\{8(2t_{2}-1)t_{1}-(4t_{2}-1)\mid t_{2} \in \mathbb{N}\right\}_{t_{1} =1 }^ {\infty}
\]
\[
\mathbf{M}=\left\{\left\{8t_{1}-3\right\},\left\{24t_{1}-7\right\},\left\{40t_{1}-11\right\},\cdots\right\}_{t_{1} =1 }^ {\infty}
\]

The family of numbers is the first well-defined equation in the corollary\eqref{sara 732}:
\[
\mathbf{N}=\left\{4(4t_{1}-1)t_{2}-(2(4t_{1}-1)+1)\mid t_{1} \in \mathbb{N}\right\}_{t_{2} =1 }^ {\infty}
\]
\[
\mathbf{N}=\left\{\left\{12t_{2}-7\right\},\left\{28t_{2}-15\right\},\left\{44t_{2}-23\right\},\cdots\right\}_{t_{2} =1 }^ {\infty}
\]

So it can be stated:
\[
\mathbf{M}=\mathbf{N}.
\]
\begin{table}[ht] \label{sara 888}
\caption{a(r,t)=(12r+4(r-1))t-(7r+(r-1))}\label{eqtable}
\renewcommand\arraystretch{1.5}
\noindent\[
\begin{array}{lcccccccr}
{}&{}&{}&{}&{5}&{}&{}&{}&{}\\
{}&{}&{}&{17}&{}&{13}&{}&{}&{}\\
{}&{}&{29}&{}&{41}&{}&{21}&{}&{}\\
{}&{41}&{}&{69}&{}&{65}&{}&{29}&{}\\
{53}&{}&{97}&{}&{109}&{}&{89}&{}&{37}\\
\end{array}
\]
\end{table}
\begin{table}[ht]
\caption{} \label{sara 1181}
\renewcommand\arraystretch{1.5}
\noindent\[
\begin{array}{c|ccccc}

&1&2&3&...&{t_{2}}\\
\hline
{1}&{5}&{17}&{29}&{...}&{12t_{2}-7}\\

{2}&{13}&{41}&{69}&{...}&{28t_{2}-15}\\

{3}&{21}&{65}&{109}&{...}&{44t_{2}-23}\\

{4}&{29}&{89}&{149}&{...}&{60t_{2}-31}\\

{\vdots	}&{\vdots	}&{\vdots	}&{\vdots	}&{\vdots	}&{\vdots}\\

{t_{1}}&{8t_{1}-3}&{24t_{1}-7}&{40t_{1}-11}&{...}&{\mathbf{M}= \mathbf{N}}\\
\end{array}
\]
\end{table}
\end{proof}
\begin{remark}
According to what was stated in corollary \ref{sara 555}:\\
The family of numbers of the two well-defined equations expressed has a triangular structure in the form of\eqref{sara 888}. The diagonal numbers to the right are sets belonging to the well-defined equation\eqref{sara 732} and the diagonal numbers to the left are sets belonging to the well-defined equation\eqref{sara 731}.
\end{remark}

\begin{cor}
Based on theorem \ref{sara 097}, a well-defined equation for Erd\"{o}s-Straus is expressed based on triangular numbers:
\[
\frac{4}{n}=\frac{1}{\frac{n+1}{4}}+\frac{1}{\binom{n+1}{2}}+\frac{1}{\binom{n+1}{2}} \quad,\quad
n\equiv -1 \pmod{4}.
\]
\end{cor}
\begin{proof}
Consider the specific state of partitions in the following form:
\[ 
4\left(\frac{n+1}{2}\right)=2n+1+1 
\]
The partition will now be multiplied by$\frac{1}{\binom{n+1}{2}}$:
\[
 \frac{4}{n}=\frac{1}{\frac{n+1}{4}}+\frac{1}{\binom{n+1}{2}}+\frac{1}{\binom{n+1}{2}}.
\]
\end{proof}
\subsection{Partitions with natural sizes}\label{sarabor}
From theorem \ref{sara 098}, partitions of rational numbers, which are obtained from partitions of positive integers, with natural sizes can be deduced and relatively prime.
\begin{cor}
	\[
	\forall r\in\mathbb{N} \quad,\quad \forall r^{'}\in\mathbb{N}_{r}
	\]
	\begin{align*}
	&\frac{4}{4r^{'}(2r-(r^{'}-1))t-(2r+1)}&=\\
	&=\frac{1}{r^{'}(2r-(r^{'}-1))t}&+\\
	&+\frac{1}{r^{'}t(4r^{'}(2r-(r^{'}-1))t-(2r+1))}&+\\
	&+\frac{1}{(2r-(r^{'}-1))t(4r^{'}(2r-(r^{'}-1))t-(2r+1))}
	\end{align*}
\[
\forall t\in\mathbb{N}.
\]
\end{cor}
\begin{proof}
	By placing the following constraints in theorem \ref{sara 098}, the equation is proved to be correct.
	\[
		k=4 \quad,\quad w=r^{'}(2r-(r^{'}-1))   \quad,\quad c=2r-(r^{'}-1)  \quad,\quad  d=r^{'}.
	\]
\end{proof}

\subsection{Modular equations}\label{sarababa}
\begin{definition}
Equations in the form of equations \ref{sara..} , \ref{sara...} and \ref{sara....} are called modular equations.
\end{definition}
 We have stated the tenth equation of theorem \ref{sara arab} in the following form:
 \[ 
\forall k\in\mathbb{N}\quad,\quad \forall s\in\mathbb{N}\quad,\quad \exists m,v,b\in\mathbb{N}
\]
\[
 kmv>\mp b \quad,\quad  \mp sb\pm v\mid mvs(kmv\pm b)
 \]
 \begin{equation}\label{saram}
 \frac{k}{kmv\pm b}=\frac{1}{mv}\mp\frac{1}{ms(kmv\pm b)}+\frac{1}{\frac{mvs(kmv\pm b)}{\mp sb\pm v}}
 \end{equation}
 According to the following congruence:
\begin{equation}\label{saramm}
s\equiv b \pmod{a} \rightarrow b\in\{0,1,\cdots,a-1\}
\end{equation}
Equation \ref{saram} can be expressed using the congruence \ref{saramm} in the following form:
\begin{align}
\label{sara..}
&\frac{k}{kmv\pm b}=\frac{1}{mv}\mp\frac{1}{m(az)(kmv\pm b)}+\frac{1}{\frac{mv(az)(kmv\pm b)}{\mp (az)b\pm v}}&\\\label{sara...}
&\frac{k}{kmv\pm b}=\frac{1}{mv}\mp\frac{1}{m(az+1)(kmv\pm b)}+\frac{1}{\frac{mv(az+1)(kmv\pm b)}{\mp (az+1)b\pm v}}&\\
&\quad\quad\vdots\quad\quad\quad\vdots\quad\quad\quad\quad\quad\quad\vdots\quad\quad\quad\quad\quad\quad\quad\quad\quad\quad\quad\vdots\\\label{sara....}
&\frac{k}{kmv\pm b}=\frac{1}{mv}\mp\frac{1}{m(a(z+1)-1)(kmv\pm b)}+\frac{1}{\frac{mv(a(z+1)-1)(kmv\pm b)}{\mp (a(z+1)-1)b\pm v}}.
\end{align}
\begin{example}
	If we consider $a=2$ we will have:
	\[
	s\equiv b \pmod{2} \rightarrow b\in\{0,1\}
	\]
	\begin{align*}
		&\frac{k}{kmv\pm b}=\frac{1}{mv}\mp\frac{1}{m(2z)(kmv\pm b)}+\frac{1}{\frac{mv(2z)(kmv\pm b)}{\mp (2z)b\pm v}}&\\
		&\frac{k}{kmv\pm b}=\frac{1}{mv}\mp\frac{1}{m(2z+1)(kmv\pm b)}+\frac{1}{\frac{mv(2z+1)(kmv\pm b)}{\mp (2z+1)b\pm v}}
	\end{align*}
\end{example}
Based on this, we will consider the following equations:
	\begin{align}\label{sarammm}
	&\frac{4}{4t+1}=\frac{1}{t+1}+\frac{1}{(t+1)(2s)(4t+1)}+\frac{1}{\frac{(t+1)(2s)(4t+1)}{6s-1}}&\\\label{saramnm}
	&\frac{4}{4t+1}=\frac{1}{t+1}+\frac{1}{(t+1)(2s+1)(4t+1)}+\frac{1}{\frac{(t+1)(2s+1)(4t+1)}{6s+2}}
\end{align}
In equation \ref{sarammm} by placing the following constraints:
\[
s=1 \quad,\quad t=\{10r-6\mid r\in\mathbb{N}\}
\]
Will be achieved:
\begin{equation}\label{sarammmm}
\frac{4}{40r-23}=\frac{1}{5(2r-1)}+\frac{1}{2(2r-1)(40r-23)}+\frac{1}{10(2t-1)(40t-23)}
\end{equation}
\[
\forall r\in\mathbb{N}
\]
And using theorem \ref{sara 098}, we express the following equation: 
\begin{equation}\label{sarammmmm}
\frac{4}{40t-3}=\frac{1}{10t}+\frac{1}{5t(40t-3)}+\frac{1}{10t(40t-3)} \quad,\quad \forall t\in\mathbb{N}
\end{equation}
Equations \ref{sarammmm} and \ref{sarammmmm} cover the following set:
\[
S=\{\{40t-3\}\cup\{40r-23\}\}=\{17,37,57,\cdots\}
\]
Equations \ref{sarammm} and \ref{saramnm} give us different partitions of an rational number.\\
The following set can be considered to express a specific case:
	\[
	\left\{792t+217\right\}_{t=0}^{\infty}
	\]
	According to equation \ref{sarammm}, there are at least two different decompositions for each member of this set. In general, in modular equations, sets can be expressed in which each member of the set has more than one different decompositions.
\begin{example}
	\[
	451=11\times 41
	\]
\[
\frac{4}{1801}=\frac{1}{451}+\frac{1}{\frac{451\times4\times1801}{11}=295364}+\frac{1}{3249004}
\]
\[
\frac{4}{1801}=\frac{1}{451}+\frac{1}{\frac{451\times14\times1801}{41}=277354}+\frac{1}{11371514}
\]
\end{example}
\section{Approximation $n$}

In this section, the focus is on the values to which $n$ belongs and the reasons for the inefficiency of problem solving so far.\\
According to what we said in the introduction one of the solutions to the problem is the Schinzel method which has not yet provided a complete answer to the conjecture.\\
This solution has a defect that makes it never easy to give a complete answer to the problem in the form of equation\eqref{sara 0021}, and this defect is created by the number one. \\ 
When you want to express a well-defined equation for the Erd\"{o}s-Straus equation in the form of the Schinzel method for the set $ \left\{2t-1\right\}_{t=1}^{\infty} $, it is impossible because of the number one; Because the number four in Erd\"{o}s-Straus equation has no answer.\\
In this case, the authors of the articles consider the set $ \left\{2t-1\right\}_{t=1}^{\infty} $ in the form  $ \left\{2t+1\right\}_{t=1}^{\infty} $, which is a conceptual error. Since  according to the following equation if we consider $n$ as a collective expression, there is no answer in Erd\"{o}s-Straus equation; Unless there is an answer in $t=0$.
\[
\forall k \in \mathbb{N}-\{1\} \quad,\quad \exists c,d \in \mathbb{N} \quad,\quad k-1\mid d  \quad,\quad d\mid c 
\]
\begin{equation*}
	\frac{k}{c\mp d}=\frac{1}{\frac{c\mp d}{k-1}}+\frac{1}{c}\pm\frac{1}{\frac{c(c\mp d)}{d}}
\end{equation*}
Therefore, no equation can be formed for set $ \left\{2t+1\right\}_{t=1}^{\infty} $.\\
For this reason, Schinzel method cannot express well-defined equations for each set of numbers belonging to the following congruence:
\begin{center}
$n-1 \equiv 0 \pmod{r}$, where $n$ and $r$ belonging to positive integers.
\end{center}
If we solve the problem using the Schinzel method, we must express countless equations in the form of equation\eqref{sara 0021} or prove the equation using the inductive method.\\
To solve the inefficiency of the Schinzel method, we generalized it, because in this case the well-defined equations cover the family of sets of positive integers, which we discussed in the previous sections.\\
We now turn to other solution that show the inefficiency of the inductive method.\\
According to the seventh equation of theorem \ref{sara arab} the determined subsets can be obtained by placing the variables.\\
For set $\left\{8t-1\right\}_{t=1}^{\infty}$:
\[ k=4 \quad v=8t \quad w=2t \] 
For set $\left\{8t-5\right\}_{t=1}^{\infty}$:
\[ k=4 \quad v=8t-4 \quad w=2t-1 \] 
Sets $\left\{8t-1\right\}_{t=1}^{\infty}$, $\left\{8t-5\right\}_{t=1}^{\infty}$, $\left\{8t-3\right\}_{t=1}^{\infty}$ and $\left\{24t-7\right\}_{t=1}^{\infty}$which is a subset of $\left\{8t-7\right\}_{t=1}^{\infty}$, can be obtained using the equation of corollary \ref{sara 91} .

For the set $\left\{8t-3\right\}_{t=1}^{\infty}$, we express two different partitions using corollarys \ref{sara 91} and \ref{sara sara} :\\
According to corollary \ref{sara sara} :
\[
 8t-3 \equiv -(1+2)\pmod{4\lcm(1,2)}
\]
According to corollary \ref{sara 91} :
\[
q=3 \quad,\quad v=8t-3 \quad,\quad \frac{qv+1}{4}=6t-2 \quad,\quad \frac{6t-2}{2}=3t-1
\]
So we will have:
\[
\frac{4}{8t-3}=\frac{1}{2t}+\frac{1}{\frac{2t(8t-3)}{2}}+\frac{1}{2t(8t-3)}=\frac{1}{6t-2}+\frac{1}{3t-1}+\frac{1}{(6t-2)(8t-3)}.
\]
The only remaining subset of odd numbers are  $\left\{8t-7\right\}_{t=1}^{\infty}$.\\
To provide tips for the set $\left\{8t-7\right\}_{t=1}^{\infty}$, we describe the process of obtaining a well-defined equation for this set:
\[ 1=\frac{8t-7}{8t-6}+\frac{1}{8t-6}=\frac{8t-7}{2(4t-3)}+\frac{1}{2(4t-3)} \]
Now we will multiply the equation in $\frac{1}{b}$: 
\[\frac{4t-3}{b}=\frac{8t-7}{2b}+\frac{1}{2b} \]
We add the numbers $\{-4,+4\} $ to this equation:
\[4=\frac{8t-7}{2b}+\frac{1}{2b}+4-\frac{4t-3}{b} \]
Now we will multiply the equation by $\frac{1}{8t-7}$:
\[\frac{4}{8t-7}=\frac{1}{2b}+\frac{1}{2b(8t-7)}+\frac{4b-4t+3}{b(8t-7)} \]
Considering the following relation, the equation is expressed in the form of the Erd\"{o}s-Straus equation:
 \[  4b-4t+3=b \rightarrow b=\frac{4t-3}{3}  \]
We rewrite the equation by inserting this relation:
\[\frac{4}{8t-7}=\frac{1}{\frac{2(4t-3)}{3}}+\frac{1}{\frac{2(4t-3)(8t-7)}{3}}+\frac{1}{(8t-7)}. \]
Clearly this equation does not have the answer according to Erd\"{o}s-Straus equation.\\
Now we consider the set in the form $\left\{8t+1\right\}_{t=1}^{\infty}$ and express a well-defined equation for it:
\[ 1=\frac{8t+1}{8t+2}+\frac{1}{8t+2} \]
Now we will multiply the equation by $\frac{1}{b}$: 
\[ \frac{1}{b}=\frac{8t+1}{b(8t+2)}+\frac{1}{b(8t+2)} \]
We add the numbers $ \{-4,+4\} $ to this equation:
\[4=\frac{8t+1}{b(8t+2)}+\frac{1}{b(8t+2)} +4-\frac{1}{b} \]
The equation will be obtained as follows:
\[\frac{4}{8t+1}=\frac{1}{b(8t+2)}+\frac{1}{\frac{b(8t+1)}{4b-1}}+\frac{1}{b(8t+2)(8t+1)} \]
For the equation to be correct, the fraction $\frac{b(8t+1)}{4b-1}$ must be equal to a positive integer, so we can say:
\[ \frac{b(8t+1)}{4b-1}=2t+h\quad,\quad h\in\mathbb{N}\]
We express the equation as follows:
\[  h=\frac{2t+b}{4b-1} \]
For some values can be expressed as follows:
 \[ b=1 \rightarrow h=\frac{2t+1}{3} \rightarrow t=\left\{3m+1\mid m \in \mathbb{N} \right\} \]
 \[ b=2 \rightarrow h=\frac{2t+2}{7} \rightarrow t=\{7m-1\mid m \in \mathbb{N} \} \]
 \[ b=3 \rightarrow h=\frac{2t+3}{11} \rightarrow t=\{11m+4\mid m \in \mathbb{N} \} \]
In general, we can say:
\[
\forall b \in \mathbb{N} \quad,\quad \exists r \in \mathbb{Z} \quad,\quad 4b-1 \mid 2r+b \quad,\quad t=\{(4b-1)m+r\mid m \in \mathbb{N}\} 
\]

\[\frac{4}{8t+1}=\frac{1}{b(8t+2)}+\frac{1}{\frac{b(8t+1)}{4b-1}}+\frac{1}{b(8t+2)(8t+1)} \]
In the well-defined equations of sets $\left\{8t-7\right\}_{t=1}^{\infty}$ and $\{8t+1\}_{t=1}^{\infty}$, we did not get equations that cover all numbers in sets.\\
Now consider the set of natural numbers in the following form:
\begin{equation*}
	 	\mathbb{N}=\left\{
	 	\begin{array}{llllll}
	 		 2t\cr
	 		 2t-1\left\{
	 	\begin{array}{llllll}
	 		  4t-1\cr
	 		  4t-3\left\{
	 	\begin{array}{llllll}
	 		8t-3\cr
	 		8t-7\left\{
	 	\begin{array}{llllll}
	 		24t-7\cr
	 		24t-15\cr
            24t-23\left\{
	 	\begin{array}{llllll}
	 		120t-23\cr
            120t-47\cr
            120t-71\cr
            120t-95\cr
            120t-119
	 	\end{array}
	 	\right.
	 	\end{array}
	 	\right.
	 	\end{array}
	 	\right.
	 	\end{array}
	 	\right.
	 	\end{array}
	 	\right\}_{t=1}^{\infty}
	 \end{equation*}
For sets $\left\{120t-23\right\}_{t=1}^{\infty}$ and $\left\{120t-47\right \}_{t=1}^{\infty}$, we have expressed well-defined equations based on the topics in the previous sections.\\
According to the corollary\ref{sara -1} for set  $\left\{120t-23\right\}_{t=1}^{\infty}$ a well-defined equation can be expressed as follows:
\[
\frac{4}{120t-23}=\frac{1}{5(6t-1)}+\frac{1}{10t(120t-23)}+\frac{1}{10t(6t-1)(120t-23)}
\]
And for the set $\left\{120t-47\right\}_{t=1}^{\infty}$ we express the following well-defined equation:
\[
\frac{4}{40v-7}=\frac{1}{10v}+\frac{1}{5v(40v-7)}+\frac{1}{2v(40v-7)}
\]
\[
\left\{120t-47\right\}_{t=1}^{\infty} \subset \left\{40v-7\right\}_{v=1}^{\infty}
\]
Set $\left\{120t-95\right\}_{t=1}^{\infty}$ using the well-defined equation of set $\left\{24t-19\right\}_{t=1}^{\infty}$ is obtained as a corollary \ref{sara 91}.
Set $\left\{120t-119\right\}_{t=1}^{\infty}$ according to what we said, it should be divided into subsets.
But set $\left\{120t-71\right\}_{t=1}^{\infty}$ cannot be expressed in the form of a well-defined equation, that means an equation for all numbers in the set; Because of the point that is in this set.\\
The point of the set $\left\{120t-71\right\}_{t=1}^{\infty}$ is related to positive integers partitions and the generalization of these types of partitions. Because partitions of positive rational numbers  are obtained from both types partitions of positive integers.\\
There are sets that can only be obtained from partitions of positive integers with rational sizes; For example, we present a subset of Mersenne numbers:\\
Partitions a subset of rational numbers in the form of $\left\{\frac{4}{2^{t}-1}\right\}_{t=1}^{\infty}$ are obtained only from the partition of positive integers with rational sizes.\\
Using Mathematica software, we can investigate the accuracy of this statement.\\
The first few numbers of this subset:
\[
\mathbb{M}_{1}=\{3,7,31,127,4095,8191,16383,...\}\subset\left\{2^{t}-1\right\}_{t=1}^{\infty}
\]
Also, more numbers can be obtained from the set $\mathbb{M}_{1}$ using Mathematica software.\\
And the set $\mathbb{M}_{2}$ is obtained from both types partitions of positive integers:
\[
\mathbb{M}_{2}=\left\{2^{t}-1\right\}_{t=1}^{\infty}-\mathbb{M}_{1}
\]
According to remark \ref{sara spe} it can be obtained that the set $\left\{\frac{4}{2^{t}-1}\right\}_{t=1}^{\infty}$ will be obtained from the partition of positive integers with rational size:
\[
 2^{t}-1 \equiv -\left(\frac{1}{2}+\frac{1}{2}\right)\pmod{4}.
\]
There are three ways to express the Three-part partition of an rational number:
\begin{enumerate}
\item
Generalize its Two-part partition (if there are).
\item
Obtain its partition from partitions of integers with positive integer sizes (if there are).
\item
Obtain its partition from partitions of integers with positive rational sizes (if there are).
\end{enumerate}
Rational numbers are divided into several categories: some are obtained in all three ways, some are obtained in two ways, and some are obtained only in the one way.\\
The set of rational numbers that can only be obtained in the third way is called $\mathbb{A}_{1}$, and the set $\mathbb{M}_{1}$ is a subset of this set.\\
Below we introduce some other $\mathbb{A}_{1}$ members:
\[
\mathbb{A}_{1}=\{\mathbb{M}_{1},409,5569,...\}.
\]
\begin{example}
Now we express the specified numbers from the set $\mathbb{A}_{1}$ using theorem \ref{sara 098} :
\[
 409 \equiv -\left(\frac{1}{2}+\frac{13}{2}\right)\pmod{4\lcm(1,13)}
\]
\[
5569 \equiv -\left(\frac{1}{2}+\frac{141}{2}\right)\pmod{4\lcm(1,141)}.
\]
\end{example}
A set of rational numbers that can be obtained from both the second and third ways and do not have a Two-part partition is called $\mathbb{A}_{2}$.
\begin{example}
To prove that set $\mathbb{A}_{2}$ is not empty, we use the following example using theorem \ref{sara 098} :
\[
3049 \equiv -\left(\frac{5}{2}+\frac{17}{2}\right)\pmod{4\lcm(5,17)}
\]
The partition of number $\frac{4}{3049}$, which is obtained from an integer partition with integer sizes, is shown in table $\ref{eeeeqtable}$ .

\end{example}
Therefore, positive rational partitions are obtained from different modes of positive integers partitions.\\
From what has been said, it can be concluded that there are reasons why it is not possible to express a well-defined equation for a set of numbers, which will be presented as follows:
\begin{enumerate}
\item
Existence of number one in the desired set.
\item
The general sentence of the desired set must be in the form of a collective expression that has no answer in the value of zero.
\item
In the set in question, there should be numbers that are obtained from different modes of partitions of positive integers.
\end{enumerate}
Set $\left\{120t-71\right\}_{t=1}^{\infty}$ has a third item. For example, in this set there is the number 49 whose second root is a Mersenne number belonging to the set $\mathbb{M}_{1}$, and there are numbers from the set $\mathbb{A}_{1}$, in this set, and from numbers that are obtained from partitions of positive integers with integer sizes, therefore an equation for this set that covers all numbers is impossible.\\
Using the stated methods to obtain well-defined equations, it can also be concluded that this is impossible. To make this clearer, consider the following equation:
\[
\frac{k}{mr-c}=\frac{1}{c(vr+l)}+\frac{1}{(vr+l)(mr-c)}+\frac{1}{c\left(\frac{vr+l}{b}\right)(mr-c)}
\]
We will get the sum of the unit fractions:
\[
\frac{1}{c(vr+l)}+\frac{1}{(vr+l)(mr-c)}+\frac{1}{c\left(\frac{vr+l}{b}\right)(mr-c)}=\frac{mr-c+c+b}{c(vr+l)(mr-c)}
\]
In order for the fraction $ \frac{1}{c\left(\frac{vr+l}{b}\right)(mr-c)} $ to be obtained in the form of a unit fraction, we must have:
\[
mr+b=k(vr+l) \rightarrow b=(kv-m)r+kl\quad,\quad\frac{vr+l}{b}=\frac{vr+l}{(kv-m)r+kl}=h \quad,\quad h\in \mathbb{N}
\]
In this case we can conclude that:
\[
\frac{vr+l}{(kv-m)r+kl}=h \rightarrow l=\frac{r((kv-m)h-v)}{kh-l}
\]
Hence, whatever state of the equation we consider, $r$ does not cover all positive integers. So the rational numbers below the well-defined equation are not found for them.
\[
\frac{4}{120t+49}=\frac{4}{10(12t+5)-1}=\frac{4}{10r-1}
\]
In this case, they should be divided into subcategories.\\
One of the modes of the partitions of the set $\left\{120t-119\right\}_{t=1}^{\infty}$ is as follows:
\begin{equation*}
	\left\{120t-119\right\}_{t=1}^{\infty}=\left\{
	\begin{array}{llllll}
		1320t-119\cr
		1320t-239\cr
		1320t-359\cr
		1320t-479\cr
		1320t-599\cr
		1320t-719\cr
		1320t-839\cr
		1320t-959\cr
		1320t-1079\cr
		1320t-1199\cr
		1320t-1319
	\end{array}
	\right\}_{t=1}^{\infty}
\end{equation*}
Partitions of sets $\left\{8t-7\right\}_{t=1}^{\infty}$, $\left\{24t-23\right\}_{t=1}^{\infty}$ and $\left\{120t-119\right\}_{t=1}^{\infty}$ have a relation to the following form in a subset:
\begin{align*}
	\left\{24t-15\right\}&=\left\{3(8t-5)\right\}\\
	\left\{120t-95\right\}&=\left\{5(24t-19)\right\}\\
	\left\{1320t-1199\right\}&=\left\{11(120t-109)\right\}
\end{align*}
This relationship ends here despite the fact that the number $1321$ is prime.\\
This relationship is due to the following equation:
\[
x(at-(a-1))+1=a(xt-b)
\]
To establish this equation we must have:
\[
x+1=a(x-b)
\]
In this case, if we consider $x + 1$ as a prime, we can see that $x+1=a$ and $b=x-1$ are established, so we will reach a multiple of $xt+1 $ when the relation reaches the first number.\\
We use corollary \ref{sara 91} to obtain the well-defined equations of sets $\left\{1320t-1079\right\}_{t=1}^{\infty}$ and $\left\{120t-109\right\}_{t=1}^{\infty}$ .\\
For set $\left\{1320t-1079\right\}_{t=1}^{\infty}$ in corollary \ref{sara 91}, we put $q=23$ and for set $\left\{120t-109\right\}_{t=1}^{\infty}$ we divide it into two subsets:
\[
\left\{120t-109\right\}_{t=1}^{\infty}=\left\{240t-229\right\}_{t=1}^{\infty}\cup 	\left\{240t-109\right\}_{t=1}^{\infty}
\]
To obtain the set $\left\{240t-109\right\}_{t=1}^{\infty}$ in corollary \ref{sara 91}, $q=13$.\\
For the set $\left\{240t-109\right\}_{t=1}^{\infty}$ in corollary \ref{sara 91}, $q=5$.\\
The third case can not be easily examined in sets because the position of different states in some sets is not clear. Therefore, the schnitzel method can not be considered to solve this equation and the best solution is to generalize the schinzel method, which was the main focus of the article.\\
Because in this solution you are dealing with a family of sets that have more properties than other solutions that can be used to obtain acceptable sets in the equation and to obtain properties from the equation.\\
According to theorem \ref{sara  098} , for the two sets $\left\{120t+1\right\}_{t=1}^{\infty}$ and $\left\{120t-71\right\}_{t=1}^{\infty}$, we have created table \ref{eeeeqtable} , which we show that theorem \ref{sara 098} will cover them.\\
In table \ref{eeeeqtable} , our focus is on positive integer partitions using positive integer sizes; therefore, we have omitted numbers that can only be obtained using positive rational sizes in theorem \ref{sara  098} .
\begin{table}[ht] 
	\caption{$\left\{120t^{'}+1\mid t^{'}\in\mathbb{N}\:,\:120t^{'}+1\in\mathbb{P}\right\}$\\\quad$\left\{120t^{'}-71\mid t^{'}\in\mathbb{N}\:,\:120t^{'}-71\in\mathbb{P}\right\}$}\label{eeeeqtable}
	\renewcommand\arraystretch{1.5}
	\noindent\[
	\begin{array}{|c|c|c|c|c|c|c|c|c|}
		\hline
		120t^{'}+1&v&t&120t^{'}-71&v&t&120t^{'}-71&v&t\\
		\hline
		{241}&{11t-1}&{22}&409&-&-&8209&7t-2&1173\\
		\hline
		{601}&{7t-1}&{86}&769&11t-1&70&8329&7t-1&1190\\
		\hline
		{1201}&{31t-8}&{39}&1009&23t-3&44&8689&11t-1&790\\
		\hline
		{1321}&{7t-2}&{189}&1129&87t-2&13&8929&19t-1&470\\
		\hline
		{1801}&{11t-3}&{164}&1249&19t-5&66&9049&7t-2&1293\\
		\hline
		{2161}&{23t-1}&{94}&1489&7t-2&213&9649&1609t-5&6\\
		\hline
		{2281}&{7t-1}&{326}&1609&7t-1&230&9769&23t-6&425\\
		\hline
		{2521}&{87t-2}&{29}&2089&19t-1&110&10009&7t-1&1430\\
		\hline
		{3001}&{19t-1}&{158}&2689&23t-2&117&10369&19t-5&546\\
		\hline
		{3121}&{7t-1}&{446}&3049&3055t-6&1&10729&7t-2&1533\\
		\hline
		{3361}&{99t-5}&{34}&3169&7t-2&453&11329&11t-1&1030\\
		\hline
		{4201}&{11t-1}&{382}&3529&19t-5&186&11689&7t-1&1670\\
		\hline
		{4441}&{19t-5}&{234}&3769&23t-3&164&12049&39t-2&309\\
		\hline
		{4561}&{39t-2}&{117}&3889&59t-5&66&12289&12375t-86&1\\
		\hline
		{4801}&{7t-1}&{686}&4129&7t-1&590&12409&7t-2&1773\\
		\hline
		{5281}&{19t-1}&{278}&4729&11t-1&430&12889&14061t-1172&1\\
		\hline
		{5521}&{11t-1}&{502}&4969&7t-1&710&13009&13013t-4&1\\
		\hline
		{5641}&{7t-1}&{806}&5209&79t-5&66&13249&7t-2&1893\\
		\hline
		{5881}&{111t-1}&{53}&5449&23t-2&237&13729&23t-2&597\\
		\hline
		{6121}&{39t-2}&{157}&5569&-&-&14449&14595t-146&1\\
		\hline
		{6361}&{7t-2}&{909}&{5689}&{7t-2}&{813}&14929&7t-2&2133\\
		\hline
		{6481}&{7t-1}&{926}&6529&7t-2&933&15289&11t-1&1390\\
		\hline
		{6841}&{11t-1}&{622}&7129&23t-1&310&15649&15651t-2&1\\
		\hline
		{6961}&{59t-1}&{118}&7369&11t-1&670&15889&7t-1&2270\\
		\hline
		{7321}&{7t-1}&{1046}&7489&7t-1&1070&16249&16317t-68&1\\
		\hline
		{7561}&{19t-1}&{398}&8089&31t-2&261&16369&16407t-38&1\\
		\hline
	\end{array}
	\]
\end{table}

According to equation\eqref{sara rf} for $k = 2$ and $ k = 3$ well-defined equations can be expressed as follows:\\
For $k=2$:
\begin{align*}
&n=\left\{2w-1\right\}_{w=1}^{\infty} \rightarrow 2w-1 \equiv -\left(\frac{1}{2}+\frac{1}{2}\right)\pmod{2\lcm(1,1)} \\
&n=\left\{2(w-1)\right\}_{w=2}^{\infty} \rightarrow 2(w-1) \equiv -\left(1+1\right)\pmod{2\lcm(1,1)}
\end{align*}
For $k=3$:
\begin{align*}
&n=\left\{3w-1\right\}_{w=1}^{\infty} \rightarrow 3w-1 \equiv -\left(\frac{1}{2}+\frac{1}{2}\right)\pmod{3\lcm(1,1)}\\
&n=\left\{3w-2\right\}_{w=1}^{\infty} \rightarrow 3w-2 \equiv -\left(1+1\right)\pmod{3\lcm(1,1)}\\
&n=\left\{3(w-1)\right\}_{w=2}^{\infty} \rightarrow 3(w-1) \equiv -\left(1+2\right)\pmod{3\lcm(1,2)}\\
&n=\left\{2(3w-2)\right\}_{w=1}^{\infty} \rightarrow 2(3w-2) \equiv -\left(2+2\right)\pmod{3\lcm(2,2)}\\
\end{align*}
According to equation\eqref{sara rf} for $k=4$ we can express:
\begin{align*}
&n=\left\{4w-1\right\}_{w=1}^{\infty} \rightarrow 4w-1 \equiv -\left(\frac{1}{2}+\frac{1}{2}\right)\pmod{4\lcm(1,1)}\\
&n=\left\{8w-3\right\}_{w=1}^{\infty} \rightarrow 8w-3 \equiv -\left(1+2\right)\pmod{4\lcm(1,2)}\\
&n=\left\{12w-3\right\}_{w=1}^{\infty} \rightarrow12w-3 \equiv -\left(\frac{3}{2}+\frac{3}{2}\right)\pmod{4\lcm(3,3)}\\
\end{align*}
For $k=4$ it is not possible to prove the equation with a definite number of equations for $n$ because the undefined number one causes the set to be subdivided into subsets and in this process we encounter other sets that numbers belonging to the set belong to different partitions of positive integers\footnote{This means that there are numbers belonging to partitions of the posittive intrger size and partitions of rational size in the set. Such as set $\left\{120t-71\right\}_{t=1}^{\infty}$.} and we have to break that set into subsets and this pattern has no end.

 \section{Equivalent conjectures}
 
 Now, based on equation 9 of theorem \ref{sara arab}, we express a well-defined equation for the subset of positive rational numbers in the form $\frac{4}{4a+1} $:
 \begin{conjecture}\label{sara tg}
 	For each prime number in the form of $ 4a+1$, a well-defined equation is expressed in the following form:
 	\[
 	\forall w\in\{4a+1\mid a\in\mathbb{N}\:\:,\:\:(4a+1)\in\mathbb{P}\} \quad,\quad \exists r,v,s\in\mathbb{N}
 	\]
 	\begin{equation} \label{sara is my life}
 		4rvs=s+rw+((4rv-1)s-rw)
 	\end{equation}
 	Where in:
 	\[
 	(4rv-1)s-rw \mid rvws
 	\]
 	According to partition\eqref{sara is my life}, it can be stated that:
 	\begin{equation}\label{key}
 	\frac{4}{w}=\frac{1}{rvw}+\frac{1}{vs}+\frac{1}{\frac{rvws}{(4rv-1)s-rw}}
 	\end{equation}
 	Equation$\eqref{key}$ is equivalent to the following equation:
 	\[
 	\frac{4rv-1}{r}=\frac{vw}{q}+\frac{w}{s} \quad,\quad q\in\mathbb{N}.
 	\]
 \end{conjecture}

 \begin{remark}\label{sara eqq}
 	If we consider the multiplication of numbers $r$ and $v$ in the following form:
 	\begin{equation}\label{sara aras}
 		rv=\frac{(4b-1)w+1}{4}\quad,\quad b\in\mathbb{N}
 	\end{equation}	
 	By placing the value $w$ we will have:
 	\begin{equation}\label{sara ras}
 		rv=4ab-a+b=(4b-1)a+b 
 	\end{equation}
 	In this case we will have:
 	\[
 	\frac{4}{w}=\frac{1}{rvw}+\frac{1}{vs}+\frac{1}{\frac{rvws}{(4rv-1)s-rw}}
 	\]
 	By placing\eqref{sara aras}:
 	\[
 	\frac{4}{w}=\frac{1}{rvw}+\frac{1}{vs}+\frac{1}{\frac{rvws}{(4b-1)ws-rw}}=\frac{1}{rvw}+\frac{1}{vs}+\frac{1}{\frac{rvs}{(4b-1)s-r}}
 	\]
 	By placing\eqref{sara ras}:
 	\begin{equation}\label{ssaarraa}
 	\frac{4}{4a+1}=\frac{1}{(4a+1)((4b-1)a+b)}+\frac{1}{vs}+\frac{1}{\frac{((4b-1)a+b)s}{(4b-1)s-r}}
 	\end{equation}
 	When the equation was obtained, we can express:
 	\[
 	\frac{4b-1}{r}=\frac{1}{s}+\frac{v}{q} \quad,\quad q\in\mathbb{N}.
 	\]
 \end{remark}
Equation\ref{ssaarraa} can create distinct partitions of its rational numbers $\frac{4}{4a+1} $.

 \begin{example}
 	Equation\ref{ssaarraa} for the rational number $\frac{4}{13}$ will give us its three partitions:
 	\[\frac{4}{13}=\frac{1}{130}+\frac{1}{10}+\frac{1}{5}\]
 	\[\frac{4}{13}=\frac{1}{130}+\frac{1}{20}+\frac{1}{4}\]
    \[\frac{4}{13}=\frac{1}{468}+\frac{1}{18}+\frac{1}{18}.\]
 \end{example}

 	A set of numbers that do not have an answer in equation \ref{ssaarraa}:
 	\[
 	f=\{4a+1\mid \exists a\in \mathbb{N}\:,\:(4a+1)\not\in\mathbb{P}\}
 	\]
 If $m$ is a set of compound numbers in the form $4a + 1$, then the compound numbers belonging to the following set have at least one answer in the equation \ref{ssaarraa}:
 \[
 z=m-f
 \]
 Hence the numbers that have the answer in equation \ref{ssaarraa}:
 \[
 R=\{\mathbb{P}\cup z \cup \{1\}\}.
 \]
 
	Given conjecture \ref{sara tg} for numbers $4a + 1$ and $4(30a) + 1$ based on equation \ref{ssaarraa}, we created tables \ref{eeeqtable} and \ref{eqtable} and showed them for some numbers.
	
 \begin{table}[pt] 
 	\caption{$\left\{120a^{'}+1\mid a^{'}\in\mathbb{N}\:,\:(120a^{'}+1)\in\mathbb{P}\right\}$}\label{eqtable}
 	\renewcommand\arraystretch{1.5}
 	\noindent\[
 	\begin{array}{|c|c|c|c|c|c|c|c|}
 		\hline
 	120a^{'}+1&(b,r,v,s)&120a^{'}+1&(b,r,v,s)&120a^{'}+1&(b,r,v,s)\\
 	\hline
 	241&(6,66,21,33)&7681&(2,22,611,5)&17041&(2,31,962,15)\\
 	\hline
 	601&(3,87,19,58)&8161&(6,54,869,54)&17401&(3,117,409,26)\\
 	\hline
 	1201&(40,385,124,11)&8521&(1,7,913,6)&17761&(1,77,173,66)\\
 	\hline
 	1321&(234,336,919,21)&8641&(6,78,637,39)&17881&(6,108,952,12)\\
 	\hline
 	1801&(30,190,282,38)&8761&(5,203,205,232)&18121&(3,63,791,16)\\
 	\hline
 	2161&(24,188,273,43)&9001&(6,57,908,190)&18481&(6,267,398,356)\\
 	\hline
 	2281&(4,94,91,47)&9241&(1,29,239,290)&19081&(4,133,538,57)\\
 	\hline
 	2521&(9,817,27,38)&9601&(2,31,542,279)&19441&(4,104,701,52)\\
 	\hline
 	3001&(8,802,29,401)&9721&(1,23,317,184)&19681&(1,29,509,290)\\
 	\hline
 	3121&(4,44,266,11)&10321&(4,41,944,164)&19801&(9,207,837,414)\\
 	\hline
 	3361&(36,124,969,1)&11161&(1,11,761,44)&20161&(2,46,767,14)\\
 	\hline
 	4201&(6,33,732,3)&12241&(3,49,687,98)&20521&(16,401,806,401)\\
 	\hline
 	4441&(6,28,912,4)&12601&(27,559,603,86)&20641&(3,63,901,126)\\
 	\hline
 	4561&(1,11,311,44)&12721&(1,29,329,290)&21001&(4,338,233,169)\\
 	\hline
 	4801&(4,28,643,28)&12841&(6,84,879,6)&21121&(3,57,1019,38)\\
 	\hline
 	5281&(1,17,233,102)&13441&(1,17,593,102)&21481&(3,21,2813,42)\\
 	\hline
 	5521&(1,41,101,574)&13681&(4,53,968,159)&21601&(1,17,953,102)\\
 	\hline
 	5641&(6,34,954,2)&13921&(1,53,197,954)&21841&(28,158,3836,395)\\
 	\hline
 	5881&(1,11,401,44)&14281&(3,53,741,10)&21961&(3,123,491,164)\\
 	\hline
 	6121&(2,13,824,3)&14401&(3,43,921,172)&22441&(9,87,2257,12)\\
 	\hline
 	6361&(3,21,833,42)&15121&(5,85,845,90)&22921&(5,65,1675,10)\\
 	\hline
 	6481&(4,28,868,4)&15241&(1,23,497,184)&23041&(9,513,393,228)\\
 	\hline
 	6841&(2,4,2993,11)&15361&(1,41,281,574)&23761&(1,71,251,1704)\\
 	\hline
 	6961&(1,23,227,184)&15601&(2,34,803,8)&24001&(1,47,383,752)\\
 	\hline
 	7321&(1,17,323,12)&16561&(3,51,893,34)&24121&(2,122,346,549)\\
 	\hline
 	7561&(1,53,107,954)&16921&(8,133,986,342)&24481&(2,62,691,558)\\
 	\hline
 	\end{array}
 	\]
 \end{table}
\begin{table}[ht] 
	\caption{$\left\{4a+1\mid a\in\mathbb{N}\:,\:(4a+1)\in\mathbb{P}\right\}$}\label{eeeqtable}
	\renewcommand\arraystretch{1.5}
	\noindent\[
	\begin{array}{|c|c|c|c|c|c|}
		\hline
		4a+1&(b,r,v,s)&4a+1&(b,r,v,s)&4a+1&(b,r,v,s)\\
		\hline
		{5}&{(1,1,4,1)}&{281}&{(2,6,82,1)}&641&{(2,6,187,4)}\\
		\hline
		{13}&{(1,1,10,1)}&{293}&{(1,4,55,2)}&653&{(1,5,98,2)}\\
		\hline
		{17}&{(2,1,30,1)}&{313}&{(1,5,47,2)}&661&{(1,4,124,2)}\\
		\hline
		{29}&{(1,1,22,1)}&{317}&{(1,7,34,3)}&673&{(1,5,101,10)}\\
		\hline
		{37}&{(1,1,28,1)}&{337}&{(1,11,23,44)}&677&{(1,4,127,2)}\\
		\hline
		{41}&{(2,1,72,1)}&{349}&{(2,13,47,26)}&701&{(5,37,90,37)}\\
		\hline
		{53}&{(1,1,40,1)}&{353}&{(1,5,53,21)}&709&{(1,7,76,3)}\\
		\hline
		{61}&{(1,1,46,1)}&{373}&{(1,4,70,2)}&733&{(1,10,55,5)}\\
		\hline
		{73}&{(1,1,55,2)}&{389}&{(1,4,73,2)}&757&{(1,8,71,4)}\\
		\hline
		{89}&{(2,2,78,2)}&{397}&{(3,12,91,3)}&761&{(2,18,74,6)}\\
		\hline
		{97}&{(2,2,85,1)}&{401}&{(2,9,78,5)}&769&{(4,1,2884,1)}\\
		\hline
		{101}&{(1,1,76,1)}&{409}&{(4,26,59,2)}&773&{(1,1,580,1)}\\
		\hline
		{109}&{(1,1,82,1)}&{421}&{(1,4,79,2)}&797&{(1,1,598,1)}\\
		\hline
		{113}&{(1,1,85,2)}&{433}&{(1,5,65,2)}&809&{(2,3,472,1)}\\
		\hline
		{137}&{(2,3,80,1)}&{449}&{(3,13,95,78)}&821&{(1,2,308,1)}\\
		\hline
		{149}&{(1,2,56,1)}&{457}&{(2,8,100,4)}&829&{(1,2,311,1)}\\
		\hline
		{157}&{(1,2,59,1)}&{461}&{(5,30,73,40)}&853&{(1,2,320,1)}\\
		\hline
		{173}&{(1,2,65,1)}&{509}&{(2,9,99,6)}&857&{(2,3,500,1)}\\
		\hline
		{181}&{(1,2,68,1)}&{521}&{(1,17,23,6)}&877&{(1,2,329,1)}\\
		\hline
		{193}&{(1,5,29,2)}&{541}&{(1,7,58,3)}&881&{(2,6,257,6)}\\
		\hline
		{197}&{(1,2,74,1)}&{557}&{(1,11,38,10)}&929&{(1,1,697,6)}\\
		\hline
		{229}&{(1,2,86,1)}&{569}&{(2,12,83,2)}&937&{(2,2,820,1)}\\
		\hline
		{233}&{(1,5,35,2)}&{577}&{(6,42,79,21)}&941&{(1,1,706,1)}\\
		\hline
		{241}&{(6,14,99,1)}&{593}&{(1,5,89,2)}&953&{(1,1,715,2)}\\
		\hline
		{257}&{(2,5,90,1)}&{601}&{(1,11,41,44)}&977&{(2,2,855,1)}\\
		\hline
		{269}&{(3,10,74,1)}&{613}&{(1,4,115,2)}&997&{(1,1,748,1)}\\
		\hline
		{277}&{(1,4,52,2)}&{617}&{(2,6,180,2)}&1009&{(4,4,946,1)}\\
		\hline
	\end{array}
	\]
\end{table}

According to theorems \ref{saraismylife.very} and \ref{sara 098}, the Erd\"{o}s-Straus conjecture can be considered equivalent to the following conjecture:
\begin{conjecture}
	\[ \forall a \in \mathbb{N}-\{1\} \quad,\quad \exists b,c,d \in \mathbb{N} \quad,\quad c\mid b+d \]
	\begin{equation} \label{sara 30}
	a \equiv-\left(\frac{b}{c}+\frac{d}{c}\right)\pmod{4\lcm(b,d)}
	\end{equation}
\end{conjecture}
Therefore, in different sections, we showed that:\\
To get rational partitions, you have to use positive integer partitions with positive rational sizes. If only positive integers are considered, there are rational numbers for which paritions cannot be expressed.
For example, numbers in the form $\frac{4}{2^{a}-1}$ for some $a$ belong to natural numbers.\\
We now present another method for the Erd\"{o}s-Straus conjecture based on the congruence of module $4$:\\
Two congruence in module $4$ are as follows:
\[
	n \equiv 0 \pmod{4} \quad and \quad  n \equiv -2 \pmod{4}
\]
In this case, we use the partitions of positive rational numbers in the form $\frac{2}{l}$:
\[
\forall l\in\mathbb{N}\quad,\quad\exists i,j\in\mathbb{N}
\]
\[
\frac{4}{n}=\frac{2}{l}=\frac{1}{i}+\frac{1}{j}
\]

If $l$ is an even number, then the correctness of this equation according to \ref{SARAAA} is clear, and if $l$ is an odd number then the equation $l=2x-1$ has answers, so it is $i=x$ and $j=xl$.
We now consider the following congruence mode: 
\[
  n \equiv -1 \pmod{4}
\]
This congruence, can be rewritten in the following form:
\begin{equation}\label{sar a}
n=4z-1 \rightarrow 4z=n+1
\end{equation}
According to the theorem \ref{saraismylife.very}, multiplying equation \ref{sar a} by $\frac{1}{nz}$ will give:
\[
\frac{4}{n}=\frac{1}{z}+\frac{1}{nz}
\]
In equation \ref{sar a}, using the following equation, we get an equation for Three-part partitions:
  \begin{equation}\label{sara aa}
  1=\frac{c-1}{c}+\frac{1}{c}
 \end{equation}
 Therefore:
  \[
  4z=n+\frac{c-1}{c}+\frac{1}{c} 
  \]
  Now we multiply the equation by $\frac{1}{nz}$:
  \[
  \frac{4}{n}=\frac{1}{z}+\frac{c-1}{cnz}+\frac{1}{cnz}.
  \]
  And finally, we will have the last case:
   \[
   n \equiv -3\pmod{4}
   \]
It can be rewritten:
  \begin{equation}\label{sar aa}
  	n=4z-3 \rightarrow 4z=n+3
  \end{equation}
According to the theorem \ref{saraismylife.very}, multiplying equation \ref{sar aa} by $\frac{1}{nz}$ will give:
\[
\frac{4}{n}=\frac{1}{z}+\frac{3}{nz}
\]
In this case, according to the definition of partitions of positive rational numbers, there must be $3\mid nz$ to create the equation.\\
To obtain Three-part partitions, we consider the partitions of number three  in the following form:
 \[
 3=\frac{f}{g}+\frac{h}{g} \quad,\quad f,h,g\in\mathbb{N}
 \]
 Therefore:
 \[
 4z=n+\frac{f}{g}+\frac{h}{g} 
 \]
 Now we multiply the equation by $\frac{1}{nz}$:
 \[
 \frac{4}{n}=\frac{1}{z}+\frac{f}{nzg}+\frac{h}{nzg}
 \]
For partitions number three can be expressed:
  \begin{align}\label{ar aaa}
 	3g&=1+3g-1 &\rightarrow f=1 \quad,\quad h&=3g-1\\\label{ar aaaa}
 	3g&=2+3g-2 &\rightarrow f=2 \quad,\quad h&=3g-2\\\label{ar aaaaaa}
 	3g&=3+3(g-1)&\rightarrow f=3 \quad,\quad h&=3(g-1)
 \end{align}
According to the definition of positive rational numbers partitions, we should have:
 \[
f\mid nzg \quad,\quad h\mid nzg  
\]
Therefore, according to relations \ref{ar aaa}, \ref{ar aaaa} and \ref{ar aaaaaa}, the only set that will not be answered is the following set:
\[
n=4z-3=4(6t+1)-3=24t+1
\]
Therefore, any number in the form $24t+1$ that we consider will remain unanswered sets, so we consider the general case.\\
The partition of numbers in the form $4y+3$, which we consider as:
\[
(4y+3)=\frac{m}{v}+\frac{n}{v} \quad,\quad m,n,v\in\mathbb{N}
\]
In this case it can be stated:
\[
4z=n+(4y+3)=n+\frac{m}{v}+\frac{n}{v}
\]
Now we have for $24t+1$ numbers:
\begin{equation}\label{sara-1000}
4z=(24t+1)+(4y+3)=(24t+1)+\frac{m}{v}+\frac{n}{v}
\end{equation}
According to partition\eqref{sara-1000}, partitions of numbers in the form $24t+1$ will be as follows:
\begin{equation}\label{sara w}
\frac{4}{24t+1}=\frac{1}{z}+\frac{m}{zv(24t+1)}+\frac{n}{zv(24t+1)}
\end{equation}
Where in:
\[
m\mid zv(24t+1) \quad,\quad n\mid zv(24t+1)  
\]
let the following inequality in equation\eqref{sara w}:
\[
\frac{1}{z}<\frac{1}{\frac{zv}{m}(24t+1)}\le\frac{1}{\frac{zv}{n}(24t+1)}
\]
In this case it can be stated:
\[
\frac{1}{z}<\frac{4}{24t+1}<\frac{1}{z}+\frac{1}{24t+1}+\frac{1}{24t+1}
\]
Hence we will have:
\[
\frac{1}{z}<\frac{4}{24t+1} \quad\rightarrow\quad 24t+1<4z \quad\rightarrow\quad \left \lceil \frac{24t+1}{4} \right \rceil\le z
\]
And also:
\[
\frac{2}{24t+1}<\frac{1}{z}  \quad\rightarrow\quad 2z<24t+1 \quad\rightarrow\quad z\le \left \lfloor \frac{24t+1}{2} \right \rfloor
\]
Therefore the range $z$ is determined:
\begin{equation}
6t+1 \le z \le 12t
\end{equation}
So to show that there are partitions of rational numbers $\frac{4}{24t+1}$, the following equation can be expressed: 
\[
\forall t\in\mathbb{N}\quad,\quad\exists y\in\mathbb{N} \quad4(6t+1+y)=(24t+1)+(4y+3)=(24t+1)+\frac{m}{v}+\frac{n}{v}
\]
Where in:
\[
m\mid v(24t+1)(6t+1+y) \quad,\quad n\mid v(24t+1)(6t+1+y)  
\]
And we will have:
\[
\frac{4}{24t+1}=\frac{1}{6t+1+y}+\frac{1}{\frac{v}{m}(6t+1+y)(24t+1)}+\frac{1}{\frac{v}{n}(6t+1+y)(24t+1)}.
\]
In theorem \ref{sara 098}, all possible states of this equation can be obtained from it.
\section{L.C.M and partition of positive rational numbers}
In the previous sections, we showed the relations between L.C.M and rational partitions. In this section, using the L.C.M relation, we express another conjecture equivalent to the desired conjecture.\\
Consider the partitions of the positive rational number $\frac{4}{7}$:
 \begin{align*}\label{ar aaaaa}
 	\frac{4}{7}&=\frac{1}{2}+\frac{1}{28}+\frac{1}{28}\\
 	\frac{4}{7}&=\frac{1}{4}+\frac{1}{4}+\frac{1}{14}\\
 	\frac{4}{7}&=\frac{1}{2}+\frac{1}{21}+\frac{1}{42}\\
 	\frac{4}{7}&=\frac{1}{3}+\frac{1}{6}+\frac{1}{14}\\
 	\frac{4}{7}&=\frac{1}{2}+\frac{1}{16}+\frac{1}{112}\\
 	\frac{4}{7}&=\frac{1}{2}+\frac{1}{18}+\frac{1}{63}\\
 	\frac{4}{7}&=\frac{1}{2}+\frac{1}{15}+\frac{1}{210}
 \end{align*}
Using the L.C.M relation, we will get the following positive integer partitions, which we have also specified the relation between them:
  \begin{align*}
	\left(4\times2=7+\frac{1}{2}+\frac{1}{2}\right)&\equiv A\\
	\left(4\times4=7+7+2\right)&\equiv 2A\\
		\left(4\times6=14+7+3\right)&\equiv 3A\\
		\left(4\times6=21+2+1\right)&\equiv 3A\\
	\left(4\times16=56+7+1\right)&\equiv 8A	\\
	\left(4\times18=63+7+2\right)& \equiv 9A	\\
	\left(4\times30=105+14+1\right) &\equiv 15A.	
\end{align*}
We express the partitions of the positive rational number $\frac{4}{19}$:
 \begin{align*}
	\frac{4}{19}&=\frac{1}{5}+\frac{1}{190}+\frac{1}{190}\\
	\frac{4}{19}&=\frac{1}{10}+\frac{1}{10}+\frac{1}{95}\\
	\frac{4}{19}&=\frac{1}{5}+\frac{1}{114}+\frac{1}{570}\\
	\frac{4}{19}&=\frac{1}{5}+\frac{1}{100}+\frac{1}{1900}\\
	\frac{4}{19}&=\frac{1}{5}+\frac{1}{120}+\frac{1}{456}\\
	\frac{4}{19}&=\frac{1}{5}+\frac{1}{96}+\frac{1}{9120}\\
	\frac{4}{19}&=\frac{1}{6}+\frac{1}{38}+\frac{1}{57}\\
	\frac{4}{19}&=\frac{1}{6}+\frac{1}{24}+\frac{1}{456}\\
	\frac{4}{19}&=\frac{1}{8}+\frac{1}{12}+\frac{1}{456}\\
	\frac{4}{19}&=\frac{1}{6}+\frac{1}{30}+\frac{1}{95}\\
	\frac{4}{19}&=\frac{1}{6}+\frac{1}{23}+\frac{1}{2622}
\end{align*}
Also, using the L.C.M relation, we will reach the following positive integer partitions, which we have also specified the relation between them:
  \begin{align*}
	\left(4\times5=19+\frac{1}{2}+\frac{1}{2}\right)&\equiv A\\
	\left(4\times10=19+19+2\right)&\equiv 2A\\
	\left(4\times30=114+5+1\right)&\equiv 6A\\
	\left(4\times100=380+19+1\right)&\equiv 20A\\
	\left(4\times120=456+19+5\right)&\equiv 24A	\\
	\left(4\times480=1824+95+1\right)& \equiv 1824A	\\
	\left(4\times6=19+3+2\right) &\equiv B\\
		\left(4\times24=76+19+1\right) &\equiv 4B\\
			\left(4\times6=57+38+1\right) &\equiv 4B\\
				\left(4\times30=95+19+6\right) &\equiv 5B\\
					\left(4\times138=437+114+1\right) &\equiv 23B	
\end{align*}
According to the expressed partitions, the general opinion of the author is as follows:\\
The partition of each number in form $\frac{4}{n}$, according to the definition provided for the partition of positive rational numbers, requires the expression of three necessary and sufficient conditions.\\
Suppose $n$ is a natural number, in which case if there are $a$, $b$, $c$, $d$ and $q$ belong to positive rational numbers, then we must have:

\begin{enumerate}
\item
\[\gcd(a,b,c)=1\]
\item
\[a+b+c=4q\]
\item
\[\lcm(a,b,c)=\frac{nq}{d}.\]
\end{enumerate}
Given the definition of the least common multiple\cite{sameh}:
\[
\lcm(a,b,c)=\frac{abc}{\gcd(ab,ac,bc)}
\]
Hence it can be stated that:
\[
\lcm(a,b,c)=\frac{abc}{\gcd(ab,ac,bc)}=\frac{nq(a+b+c)}{4d}
\]
In this case, the following well-defined equation can be inferred from it:
\begin{equation}\label{sara.aras}
\frac{4d}{n\gcd(ab,ac,bc)}=\frac{a+b+c}{abc}=\frac{1}{bc}+\frac{1}{ac}+\frac{1}{ab}
\end{equation}
Or in another form:
\[
\frac{4}{n\gcd(ab,ac,bc)}=\frac{1}{bcd}+\frac{1}{acd}+\frac{1}{abd}
\]
If we assume that:
\[
\gcd(ab,ac,bc)=1
\]
In this case:
\[
\frac{4}{n}=\frac{1}{bcd}+\frac{1}{acd}+\frac{1}{abd}.
\]
The numbers that equation\eqref{sara.aras} covers are as follows:
\[
n=\frac{\lcm(a,b,c)4d}{a+b+c}
\]
Therefore, the Erd\"{o}s-Straus equation can be considered equivalent to equation\eqref{sara.aras}.\\
By generalizing equation\eqref{sara.aras}, we can obtain an equation equivalent to equation\eqref{saraismylife}:
\begin{enumerate}
\item
\[\gcd(a_{1},a_{2},\cdots,a_{v})=1\]
\item
\[a_{1}+a_{2}+\cdots +a_{v}=kq\]
\item
\[\lcm(a_{1},a_{2},\cdots,a_{v})=\frac{nq}{d}\]
\end{enumerate}
So according to the definition of the least common multiple\cite{sameh}:
\[
\lcm(a_{1},a_{2},\cdots,a_{v})=\frac{(a_{1}a_{2}\cdots a_{v})=l}{\gcd\left(\frac{l}{a_{1}},\frac{l}{a_{2}},\cdots,\frac{l}{a_{v}}\right)}=\frac{n(a_{1}+a_{2}+\cdots +a_{v})}{kd}\]
Hence it can be stated that:
\[
\frac{k}{n\gcd\left(\frac{l}{a_{1}},\frac{l}{a_{2}},\cdots,\frac{l}{a_{v}}\right)}=\frac{1}{d}\left(\sum_{i=1}^{v} \frac{a_{i}}{l}\right).
\]

\section{Functions L.C.M}

\begin{definition}
	The function $ Q_{a}[b]$  is a function of counting the number of positive integers whose least common multiple with the positive integer $b $  is equal to the positive integer $ a$.
\end{definition}  
\begin{definition}
	The function $ SQ_{a}[b]$ is the function of the sum of positive integers whose least common multiple with a positive integer $b $ is equal to the positive integer  $ a$.
\end{definition}  
\begin{definition}
	The function $ MQ_{a}[b]$ is the function of the product of positive integers whose least common multiple with the positive integer $b $  is equal to the positive integer $ a$. 
\end{definition}  
\begin{definition}
	The function $\sigma(n)$ is the function of the sum of the divisors of the number $n$. We consider the decomposition of $n$ into the prime factors\cite{SARAISMYLIFE}:
	\[
	n=p_{1}^{\alpha_{1}}p_{2}^{\alpha_{2}}\cdots p_{s}^{\alpha_{s}}
	\]
	\[
	\sigma(n)=\prod_{1\le i\le s}\frac{p_{i}^{\alpha_{i}+1}}{p_{i}-1}.
	\]
\end{definition}

Now in the following theorem we express and prove these three functions.
\begin{theorem}
	\[
	\forall a\in\mathbb{N}\quad,\quad \exists b\in\mathbb{N}\quad,\quad b\mid a\quad,\quad \exists d,k\in\mathbb{N} \quad,\quad a=dk=bc
	\]
	We consider the decompositions of $b$ and $c$ into the prime factors as follows:
	\[
	b=p_{1}^{x_{1}}p_{2}^{x_{2}}\cdots p_{s}^{x_{s}}
	\]
	\[
	c=p_{1}^{y_{1}}\cdots p_{s}^{y_{s}}p_{s+1}^{y_{s+1}}\cdots p_{l}^{y_{l}}
	\]
	Hence we have:
\begin{enumerate}
	\item
	\[
	Q_{a}[b]=\tau(1)=1 \quad,\quad Q_{a}[c]=\tau\left(p_{s+1}^{y_{s+1}}\cdots p_{l}^{y_{l}}\right)
	\]
	\item
	\[
	SQ_{a}[b]=a \quad,\quad SQ_{a}[c]=\frac{a}{p_{s+1}^{y_{s+1}}\cdots p_{l}^{y_{l}}}\times\sigma\left(p_{s+1}^{y_{s+1}}\cdots p_{l}^{y_{l}}\right)
	\]
	\item
	\[
	MQ_{a}[b]=a^{\frac{\tau(a)}{2}}
	\]
	\[
	 MQ_{a}[c]=\left(p_{s+1}^{y_{s+1}}\cdots p_{l}^{y_{l}}\right)^{\tau\left(p_{1}^{x_{1}+y_{1}}p_{2}^{x_{2}+y_{2}}\cdots p_{s}^{x_{s}+y_{s}}\right)} \times k^{\tau\left(p_{1}^{x_{1}+y_{1}}p_{2}^{x_{2}+y_{2}}\cdots p_{s}^{x_{s}+y_{s}}\right)^{\frac{1}{2}}}.
	\]
\end{enumerate}

\end{theorem}
\begin{proof}
	According to the definition of the least common multiple\cite{sameh}, if $a = b$ then:
	\[
	Q_{a}[b]=\tau(b)
	\]
	And if $b = 1$, then for every $a$ belonging to natural numbers, we have:
	\[
	Q_{a}[1]=1
	\]
	If we consider $a$ as a product of two numbers, then there are two cases:
	\[
	a=bc \quad,\quad \gcd(b,c)=1 \quad or \quad \gcd(b,c)\ne1
	\] 
The first case if we have:
\[
\gcd(b,c)=1
\]
Given that the arithmetic function $\tau(n)$ is multiplicative, we can say that:
\[
\tau(a)=\tau(b)\tau(c)
\]
And it can also be said:
\[
\lcm(b,c)=bc=a
\]
Consider the following set as the set of divisors of $b$:
	\[
	b=\{1,b_{1},b_{2},\cdots,b_{v-1}\} \quad,\quad \tau(b)=v
	\]
	We can say that:
	\[
	\lcm(b,c)=a
	\]
	\[
	\lcm(b,b_{1}c)=a
	\]
	\begin{equation} \label{sara 2}
		\lcm(b,b_{2}c)=a
	\end{equation}
	\[
	\vdots
	\]
	\[
	\lcm(b,b_{v-1}c)=a
	\]
	Hence:
\[
Q_{a}[b]=\tau(b)
\]
And in this way it can be achieved that:
\[
Q_{a}[c]=\tau(c)
\]
Therefore:
	\[
	Q_{a}[a]=Q_{a}[b]Q_{a}[c].
	\]
	In this case, the function $Q_{a}[b]$ will be obtained in the following form:
	\[
	Q_{a}[b]=\frac{Q_{a}[a]}{Q_{a}[c]}=\frac{\tau(a)}{\tau(c)}=\tau(b).
	\]
	Let us now consider the second case:
	\[
	a=bc \quad,\quad \gcd(b,c)\ne1
	\]
In this case there are two other numbers such as $d$ and $k$ belong to natural numbers such that:
	\[
	a=kd \quad,\quad \gcd(k,d)=1
	\]
According to the first case, it can be stated that:
	\[
	Q_{a}[k]=\frac{\tau(a)}{\tau(d)}=\tau(k)
	\]
	According to the definition of the least common multiple\cite{sameh}, if $k$ is equal to multiplying the number $c$ by the prime factors common to $b$ with their power, then:
	\begin{equation}\label{sara,}
	Q_{a}[d]=Q_{a}[b]
	\end{equation}
	So in this case, the function $Q_{a}[b]$ will be equal to:
	\[
	Q_{a}[b]=\frac{\tau(a)}{\tau(d)}=\tau(k).
	\]
	To prove the other two functions, we do the same as we did for the function $Q_{a}[b]$:
	according to \ref{sara 2}, if we have:
		\[
	a=bc \quad,\quad \gcd(b,c)=1
	\]
	In this case, we can say:
		\[
	SQ_{a}[b]=c\sigma(b)
	\]
	\[
	MQ_{a}[b]=c^{\tau(b)} b^\frac{\tau(b)}{2}
	\]
	And according to equation \ref{sara,} can be expressed:
		\[
	SQ_{a}[b]=d\sigma(k)
	\]
		\[
	MQ_{a}[b]=d^{\tau(k)} k^\frac{\tau(k)}{2}.
	\]
\end{proof}
\begin{cor}
	Some corollaries of the function $ Q_{a}[b]$ are in the following form:
	\begin{enumerate}
		\item
		\[
		\forall m \in \mathbb{N} \quad,\quad \exists n \in \mathbb{N} \quad,\quad \gcd(m,n)=1
		\]
		\[
		Q_{mn}[mn]=Q_{mn}[m]\times Q_{mn}[n]
		\]
		\item
		\[
		\forall m, n \in \mathbb{P}
		\]
		\[
		Q_{mn}[mn]=Q_{mn}[m]+Q_{mn}[n]=Q_{mn}[m]\times Q_{mn}[n]
		\]
		\item
		\[
		\forall m \in \mathbb{P} \quad,\quad \forall n\in \mathbb{N} \quad,\quad \exists t\in \mathbb{N} \quad,\quad \gcd(t,m)=1
		\]
		\[
		Q_{m^n}[m^k]=\begin{cases}
			1 & \mbox{if } k\in\mathbb{N}_{n-1} \cup \{0\}
			\\
			n+1 & \mbox{if } k=n
		\end{cases}
		\]
		\[
		,\quad	Q_{t\times m^{n}}[m^{n}]=n+1
		\]
		\item
		\[
		\forall p\in\mathbb{P}
		\]
		\[	Q_{pk}[p]=\begin{cases}
			1 & \mbox{if } k\equiv  0 \pmod{p}
			\\
			2 & \mbox{if } k\not\equiv 0 \pmod{p}
		\end{cases}\]
	
	\item
	\[
	\forall a\in\mathbb{N}\quad,\quad\exists b\in\mathbb{N}\quad,\quad b\mid a\quad,\quad\exists x_{1},x_{2},\cdots,x_{s}\in\mathbb{N}
	\]
	\[
	\begin{cases}
		SQ_{a}[b]=x_{1}+x_{2}+\cdots+x_{s}\\
		MQ_{a}[b]=x_{1}\times x_{2}\times\cdots\times x_{s}
	\end{cases}
	\]
	\[
	\frac{SQ_{a}[b]}{MQ_{a}[b]}=\frac{1}{x_{2}\times x_{3}\times\cdots\times x_{s}}+\frac{1}{x_{1}\times x_{3}\times\cdots\times x_{s}}+\cdots+\frac{1}{x_{1}\times\cdots\times x_{s-1}}
	\]
	\end{enumerate}
\end{cor}
In this case, according to the stated corollaries, we express the series and the multiplications deduced from them.
\begin{cor}
	Series and multiplications based on the function $Q_{a}[b]$:
	\begin{enumerate}
		\item
		\[
		\forall s \in \mathbb{N} \quad,\quad \forall p \in \mathbb{P}
		\]
		\[
		\sum_{k=1}^{s} Q_{pk}[p] =2s-\left \lfloor \frac{s}
		{p} \right \rfloor
		\]
		\[
		\prod_{k=1}^{s} Q_{pk}[p] = 2^{s-\left \lfloor \frac{s}
			{p} \right \rfloor}
		\]
		
		\item
		\[
		\forall m \in \mathbb{P} \quad,\quad \forall k \in \mathbb{N}
		\]
		\[
		\sum_{x|m^k} Q_{m^k}[x]=2k+1
		\]
		\[
		\prod_{x|m^k} Q_{m^k}[x]=k+1
		\]
		\item
		\[
		\forall k\in \mathbb{N}\cup\{0\} \quad,\quad \forall p \in \mathbb{P} \quad,\quad \exists t\in \mathbb{N} \quad,\quad \gcd(t,p)=1
		\]
		\[
		\sum_{r=0}^{k}Q_{t\times p^{r}}[p^{r}]=\sum_{n=1}^{k+1}n=\binom{k+2}{2}
		\]
		\[
		\prod_{r=0}^{k}Q_{t\times p^{r}}[p^{r}]=\prod_{n=1}^{k+1}n=(k+1)!
		\]
		\item
		\[
		\forall p \in \mathbb{P} \quad,\quad \forall k \in \mathbb{N}\cup\{0\}
		\]
		\[
		\sum_{b \mid p^{k}}\frac{\tau(b)}{Q_{p^{k}[b]}}=\binom{k+1}{2}+1
		\]
		\[
		\prod_{b \mid p^{k}}\frac{\tau(b)}{Q_{p^{k}[b]}}=k!.
		\]
		
	\end{enumerate}
\end{cor}
\begin{definition}
The Euler function is the number of non-negative integers less than $n$ that are relatively prime to $n$. We consider the decomposition of $n$ into the prime factors\cite{SARAISMYLIFE}:
\[
n=p_{1}^{\alpha_{1}}p_{2}^{\alpha_{2}}\cdots p_{s}^{\alpha_{s}}
\]
\[
\varphi(n)=\prod_{i=1}^{s}\left(p_{i}^{\alpha_{i}}-p_{i}^{\alpha_{i}-1}\right).
\]
\end{definition}
\begin{cor}
	The relation between Euler function and function $Q_{a}[b]$:
	\[
	\forall p \in \mathbb{P}
	\]
	\[
	\sum_{k=1}^{p} Q_{pk}[p] = p+\varphi(p)=2p-1
	\]
	\[
	\prod_{k=1}^{p} Q_{pk}[p] = 2^{\varphi(p)}=2^{p-1}.
	\]
\end{cor}
\begin{definition}
	The Fermat quotient of an integer $k$ with respect to an odd prime $p$ is defined as:
	\[
	q_{p}(k)=\frac{k^{p-1}-1}{p}.
	\]
\end{definition}
\begin{theorem}
	The relation between function $SQ_{k-1}[p]$ and function $q_{p}(k) $:
	\[
	\forall k\in \mathbb{N}-\{1\} \quad,\quad \exists d,p \in \mathbb{N} \quad,\quad k-1=dp
	\]
	\[
	q_{p}(k) \equiv -SQ_{k-1}[p] \pmod{dp}.
	\]
\end{theorem}
\begin{proof}
	We consider the theorem in the following form:
\[
dp\mid q_{p}(k)+SQ_{k-1}[p]
\]
We now consider the functions as follows:
\[
SQ_{k-1}[p]=d+k-1
\]
\[
q_{p}(k)=\frac{k^{p-1} -1}{p}=\frac{dk^{p-1} -d}{k-1}
\]
Therefore:
\[
q_{p}(k)+SQ_{k-1}[p]=d\left(\frac{k^{p-1}-1}{k-1}\right)+d+k-1
\]
We rewrite the right side of the relation:
\[
d\left(k^{p-2}+k^{p-3}+\cdots+1\right)+d+k-1
\]
And we place the value of $k$:
\[
d\left((dp+1)^{p-2}+(dp+1)^{p-3}+\cdots+1\right)+d+dp
\]
The value of the geometric series will be as follows:
\[
S=\frac{(dp+1)^{p-1}-1}{dp}
\]
In this case we will have:
\begin{equation}\label{sarasm}
dS+d+dp=\frac{(dp+1)^{p-1}-1}{p}+d+dp
\end{equation}
Now consider the following relation:
\begin{equation}\label{sarasm2}
(dp+1)^{p-1}=(dp)^{p-1}+(p-1)(dp)^{p-2}+\cdots+(p-1)(dp)+1
\end{equation}
We add the value $dp-1$ to the sides of equation\eqref{sarasm2}:
\[
(dp+1)^{p-1}-1+dp=(dp)^{p-1}+(p-1)(dp)^{p-2}+\cdots+(p-1)(dp)+1-1+dp
\]
And finally we will have:
\[
 q_{p}(k)+SQ_{k-1}[p]=(dp)^{p-1}+(p-1)(dp)^{p-2}+\cdots+\binom{p-1}{2}(dp)^{2}+2dp^{2}
\]
So it is clear that:
\[
dp\bigg{|} (dp)^{p-1}+(p-1)(dp)^{p-2}+\cdots+\binom{p-1}{2}(dp)^{2}+2dp^{2}.
\]
\end{proof}
\begin{theorem}\label{sara 099}
     For all $n\in\mathbb{N} $, we consider the decomposition of $n$ into the prime factors:
     \[
     n=p_{1}^{\alpha_{1}}p_{2}^{\alpha_{2}}\cdots p_{s}^{\alpha_{s}}
     \]
     The following series can be expressed:
      \begin{align*}
      	&\sum_{x\mid n}Q_{n}[x]=\tau(n)+1+\sum_{i=1}^{s}(\alpha_{i}-1)&+\\
      	&+\sum_{m=1}^{s}\left(\sum_{k_{1}=1}^{s}\cdots\sum_{k_{m}=1}^{s}\prod_{i=1}^{m}\left(\alpha_{\sum_{v=1}^{i}k_{v}}-1\right)\sum_{j=1+\sum_{v=1}^{i}k_{v}}^{s}(\alpha_{j}-1)\right)&+\\
      		&+\sum_{m=1}^{s}\tau(p_{m}^{\alpha_{m}})+\sum_{m=1}^{s}\left(\tau(p_{m}^{\alpha_{m}})\sum_{m\ne i, i=1}^{s}(\alpha_{i}-1)\right)&+\\
      			&+\sum_{m=1}^{s}\tau(p_{m}^{\alpha_{m}})&\times\\&\times\left(\sum_{r=1}^{s}\left(\sum_{ k_{1}=1}^{s}\cdots\sum_{ k_{r}=1}^{s}\prod_{i=1}^{r}\left(\alpha_{m\ne l,l=\sum_{v=1}^{i}k_{v}}-1\right)\sum_{m\ne j , j=1+\sum_{v=1}^{i}k_{v}}^{s}(\alpha_{j}-1)\right)\right)&+\\
&+\sum_{m=1}^{s}\left(\sum_{k_{1}=1}^{s}\cdots\sum_{k_{m}=1}^{s}\prod_{i=1}^{m}\tau\left(p_{f_{i}=\sum_{v=1}^{i}k_{v}}^{\alpha_{f_{i}=\sum_{v=1}^{i}k_{v}}}\right)\sum_{u=1+\sum_{v=1}^{i}k_{v}}^{s}\tau\left(p_{u}^{\alpha_{u}}\right)\right)&+\\
&+\sum_{m=1}^{s}\left(\sum_{k_{1}=1}^{s}\cdots\sum_{k_{m}=1}^{s}\prod_{i=1}^{m}\tau\left(p_{f_{i}=\sum_{v=1}^{i}k_{v}}^{\alpha_{f_{i}=\sum_{v=1}^{i}k_{v}}}\right)\sum_{u=1+\sum_{v=1}^{i}k_{v}}^{s}\tau\left(p_{u}^{\alpha_{u}}\right)\right)\times\sum_{f_{i}\ne b,u\ne b, b=1}^{s}(\alpha_{b}-1)&+\\
&+\sum_{m=1}^{s}\left(\sum_{k_{1}=1}^{s}\cdots\sum_{k_{m}=1}^{s}\prod_{i=1}^{m}\tau\left(p_{f_{i}=\sum_{v=1}^{i}k_{v}}^{\alpha_{f_{i}=\sum_{v=1}^{i}k_{v}}}\right)\sum_{u=1+\sum_{v=1}^{i}k_{v}}^{s}\tau\left(p_{u}^{\alpha_{u}}\right)\right)&\times\\&\times
\left(\sum_{r=1}^{s}\left(\sum_{ k_{1}=1}^{s}\cdots\sum_{ k_{r}=1}^{s}\prod_{i=1}^{r}\left(\alpha_{f_{i}\ne l,u\ne l,l=\sum_{v=1}^{i}k_{v}}-1\right)\sum_{f_{i}\ne j,u\ne j , j=1+\sum_{v=1}^{i}k_{v}}^{s}(\alpha_{j}-1)\right)\right).
      	 \end{align*}
\end{theorem}
\begin{proof}
	According to the definition of the function $Q_{n}[x]$, in the following conventional decomposition:
	\[
	n=p_{1}^{\alpha_{1}}p_{2}^{\alpha_{2}}\cdots p_{s}^{\alpha_{s}}
	\]
	We now consider the expansion of the $n$ decomposition in the form of the following matrix.
	\[
	p_{1}^{0} \quad p_{1}^{1} \quad p_{1}^{2} \quad\cdots\quad p_{1}^{\alpha_{1}}
	\]
	\[
	p_{2}^{0} \quad p_{2}^{1} \quad p_{2}^{2} \quad\cdots\quad p_{2}^{\alpha_{2}}
	\]
	\[\vdots\]
	\[
	p_{s}^{0} \quad p_{s}^{1} \quad p_{s}^{2} \quad\cdots\quad p_{s}^{\alpha_{s}}
	\]
	For every $p_{i}^{\alpha_{i}}$ that $1\le i\le s$, there are $(\alpha_{i}-1)$ choices whose result in the function will be equal to one. In this case, according to the principle of AP\footnote{Addition principle\cite{T}} we will have: 
	\[
	\alpha_{1}+\alpha_{2}+\cdots+\alpha_{s}-s
	\]
	If we consider the divisors of the number $n$ based on permutations, then for Two-combination, we have:
	\[
	(\alpha_{1}-1)(\alpha_{2}-1)+(\alpha_{1}-1)(\alpha_{3}-1)+\cdots+(\alpha_{1}-1)(\alpha_{s}-1)
	\]
	\[
	(\alpha_{2}-1)(\alpha_{3}-1)+(\alpha_{2}-1)(\alpha_{4}-1)+\cdots+(\alpha_{2}-1)(\alpha_{s}-1)
	\]
	\[
	\vdots
	\]
	\[
	(\alpha_{s-1}-1)(\alpha_{s}-1)
	\]
	As well as for Three-combination:
	\[
	(\alpha_{1}-1)(\alpha_{2}-1)(\alpha_{3}-1)+(\alpha_{1}-1)(\alpha_{2}-1)(\alpha_{4}-1)+\cdots+(\alpha_{1}-1)(\alpha_{2}-1)(\alpha_{s}-1)
	\]
	\[
	(\alpha_{1}-1)(\alpha_{3}-1)(\alpha_{4}-1)+(\alpha_{1}-1)(\alpha_{3}-1)(\alpha_{5}-1)+\cdots+(\alpha_{1}-1)(\alpha_{3}-1)(\alpha_{s}-1)
	\]
	\[
	\vdots
	\]
	\[
	(\alpha_{1}-1)(\alpha_{s-1}-1)(\alpha_{s}-1)
	\]
	\[
	(\alpha_{2}-1)(\alpha_{3}-1)(\alpha_{4}-1)+(\alpha_{2}-1)(\alpha_{3}-1)(\alpha_{5}-1)+\cdots+(\alpha_{2}-1)(\alpha_{3}-1)(\alpha_{s}-1)
	\]
	\[
	(\alpha_{2}-1)(\alpha_{4}-1)(\alpha_{5}-1)+(\alpha_{2}-1)(\alpha_{4}-1)(\alpha_{6}-1)+\cdots+(\alpha_{2}-1)(\alpha_{4}-1)(\alpha_{s}-1)
	\]
	\[
	\vdots
	\]
	\[
	(\alpha_{2}-1)(\alpha_{s-1}-1)(\alpha_{s}-1)
	\]
	\[
	\vdots
	\]
	\[
(\alpha_{3}-1)(\alpha_{s-1}-1)(\alpha_{s}-1)
\]
\[
\vdots
\]
\[
(\alpha_{s-2}-1)(\alpha_{s-1}-1)(\alpha_{s}-1)
\]
And this is how we get all the $r$-combinations that $1\le r\le s$, until finally for the $s$-combination:
\[
(\alpha_{1}-1)(\alpha_{2}-1)(\alpha_{3}-1)\cdots(\alpha_{s-3}-1)(\alpha_{s-2}-1)(\alpha_{s-1}-1)(\alpha_{s}-1).
\]
According to the definition of the function $Q_{n}[x] $, by selecting any $p_{i}^{\alpha_{i}}$ and according to the principle of AP we will have:
\[
\tau(p_{1}^{\alpha_{1}})+\tau(p_{2}^{\alpha_{2}})+\cdots+\tau(p_{s}^{\alpha_{s}})
\]
And now we consider all the $r$-combinations that will be obtained with $\tau(p_{i}^{\alpha_{i}})$:
\[
\left(\alpha_{2}+\alpha_{3}+\cdots+\alpha_{s}-(s-1)\right)\tau(p_{1}^{\alpha_{1}})
\]
\[
\left(\alpha_{1}+\alpha_{3}+\cdots+\alpha_{s}-(s-1)\right)\tau(p_{2}^{\alpha_{2}})
\]
\[
\vdots
\]
\[
\left(\alpha_{1}+\alpha_{2}+\cdots+\alpha_{s-1}-(s-1)\right)\tau(p_{s}^{\alpha_{s}})
\]
\[
\vdots
\]
\[
(\alpha_{1}-1)(\alpha_{2}-1)(\alpha_{3}-1)\cdots(\alpha_{s-3}-1)(\alpha_{s-2}-1)(\alpha_{s-1}-1)\tau\left(p_{s}^{\alpha_{s}}\right)
\]
Thus we perform the $\tau$ function for other $r$-combinations permutations, and the last divisor remains the number $n$ itself, which is the last combination of the $\tau$ function:
\[
\tau\left(p_{1}^{\alpha_{1}}p_{2}^{\alpha_{2}}\cdots p_{s}^{\alpha_{s}}\right).
\]
\end{proof}
\begin{remark}
	The two series $\sum_{x|n}SQ_{n}[x]$ and $\sum_{x|n}MQ_{n}[x]$ can be deduced from the $\sum_{x|n}Q_{n}[x]$ series.
\end{remark}
\begin{theorem}\label{sara 09}
	For all $n\in\mathbb{N} $, we consider the decomposition of $n$ into the prime factors:
	\[
	n=p_{1}^{\alpha_{1}}p_{2}^{\alpha_{2}}\cdots p_{s}^{\alpha_{s}}
	\]
	The following series can be expressed:
	\[\sum_{x\mid n}Q_{n}[x]=\tau\left(n^2\right)=\left(2\alpha_{1}+1\right)\left(2\alpha_{2}+1\right)\cdots\left(2\alpha_{s}+1\right).\]
	
\end{theorem}
\begin{proof}
Using induction, the correctness of the theorem can be obtained.
\end{proof}
According to the two theorems \ref{sara 099} and \ref{sara 09} in general, the following corollary can be expressed:
\begin{cor}
\[
\forall q\in\mathbb{N}
\]
\begin{align*}
&\tau\left(n^q\right)=\tau(n^{q-1})+1+\sum_{i=1}^{s}(\alpha_{i}-1)&+\\
      	&+\sum_{m=1}^{s}\left(\sum_{k_{1}=1}^{s}\cdots\sum_{k_{m}=1}^{s}\prod_{i=1}^{m}\left(\alpha_{\sum_{v=1}^{i}k_{v}}-1\right)\sum_{j=1+\sum_{v=1}^{i}k_{v}}^{s}(\alpha_{j}-1)\right)&+\\
      		&+\sum_{m=1}^{s}\tau(p_{m}^{(q-1)\alpha_{m}})+\sum_{m=1}^{s}\left(\tau(p_{m}^{(q-1)\alpha_{m}})\sum_{m\ne i, i=1}^{s}(\alpha_{i}-1)\right)&+\\
      	&+\sum_{m=1}^{s}\tau\left(p_{m}^{(q-1)\alpha_{m}}\right)&\times\\&\times\left(\sum_{r=1}^{s}\left(\sum_{ k_{1}=1}^{s}\cdots\sum_{ k_{r}=1}^{s}\prod_{i=1}^{r}\left(\alpha_{m\ne l,l=\sum_{v=1}^{i}k_{v}}-1\right)\sum_{m\ne j , j=1+\sum_{v=1}^{i}k_{v}}^{s}(\alpha_{j}-1)\right)\right)&+\\
&+\sum_{m=1}^{s}\left(\sum_{k_{1}=1}^{s}\cdots\sum_{k_{m}=1}^{s}\prod_{i=1}^{m}\tau\left(p_{f_{i}=\sum_{v=1}^{i}k_{v}}^{(q-1)\alpha_{f_{i}=\sum_{v=1}^{i}k_{v}}}\right)\sum_{u=1+\sum_{v=1}^{i}k_{v}}^{s}\tau\left(p_{u}^{(q-1)\alpha_{u}}\right)\right)&+\\
&+\sum_{m=1}^{s}\left(\sum_{k_{1}=1}^{s}\cdots\sum_{k_{m}=1}^{s}\prod_{i=1}^{m}\tau\left(p_{f_{i}=\sum_{v=1}^{i}k_{v}}^{(q-1)\alpha_{f_{i}=\sum_{v=1}^{i}k_{v}}}\right)\sum_{u=1+\sum_{v=1}^{i}k_{v}}^{s}\tau\left(p_{u}^{(q-1)\alpha_{u}}\right)\right)&\times\\&\times\sum_{f_{i}\ne b,u\ne b, b=1}^{s}\left(\alpha_{b}-1\right)&+\\
&+\sum_{m=1}^{s}\left(\sum_{k_{1}=1}^{s}\cdots\sum_{k_{m}=1}^{s}\prod_{i=1}^{m}\tau\left(p_{f_{i}=\sum_{v=1}^{i}k_{v}}^{(q-1)\alpha_{\sum_{v=1}^{i}k_{v}}}\right)\sum_{u=1+\sum_{v=1}^{i}k_{v}}^{s}\tau\left(p_{u}^{(q-1)\alpha_{u}}\right)\right)&\times\\&\times
\left(\sum_{r=1}^{s}\left(\sum_{ k_{1}=1}^{s}\cdots\sum_{ k_{r}=1}^{s}\prod_{i=1}^{r}\left(\alpha_{f_{i}\ne l,u\ne l,l=\sum_{v=1}^{i}k_{v}}-1\right)\sum_{f_{i}\ne j,u\ne j , j=1+\sum_{v=1}^{i}k_{v}}^{s}(\alpha_{j}-1)\right)\right).
\end{align*}
\end{cor}
\begin{example}
\[
\tau\left(p_{1}^{ak}p_{2}^{am}\right)=mk+m\tau\left(p_{1}^{(a-1)k}\right)+k\tau\left(p_{2}^{(a-1)m}\right)+\tau\left(p_{1}^{(a-1)k}p_{2}^{(a-1)m}\right).
\]
\end{example}

\begin{theorem}
	For all $n\in\mathbb{N} $, we consider the decomposition of $n,x$ into the prime factors:
	\[n=p_{1}^{\alpha_{1}}p_{2}^{\alpha_{2}}\cdots p_{s}^{\alpha_{s}}\quad,\quad x=p_{r}^{a_{r}}p_{l}^{a_{l}}\cdots p_{i}^{a_{i}}\quad,\quad x\mid n\]
		\begin{enumerate}
		\item
	\[\sum_{xt|n}Q_{n}[xt]=\tau\left(x\right)\tau\left(\frac{n^2}{p_{r}^{\alpha_{r}+a_{r}}p_{l}^{\alpha_{l}+a_{l}}\cdots p_{i}^{\alpha_{i}+a_{i}}}\right)
	\]
or
\[
\gcd(n,xy)=1\quad,\quad\sum_{xt|n}Q_{n}[xt]=\tau(x)\tau(y)\sum_{t|\frac{n}{xy}}Q_{\frac{n}{xy}}[t]
\]
	\item
	\[
	 \exists j,d\in\mathbb{N} \quad,\quad1\le j \le s \quad,\quad 1\le d \le \alpha_{j}
	\]
	\[
	\sum_{p_{j}^{d}t\big{|}n}Q_{n}\left[p_{j}^{d}t\right]=\sum_{t|n}Q_{n}[t]-d\sum_{t\big{|}n\left(p_{j}^{\alpha_{j}}\right)^{-1}}Q_{n\left(p_{j}^{\alpha_{j}}\right)^{-1}}[t]=\tau\left(n^2\right)-d\tau\left(\frac{n^2}{p_{j}^{2\alpha_{j}}}\right).
	\]
\end{enumerate}
\end{theorem}
\begin{proof}
	Considering the factoring action and the multiplicative nature of the $\tau$ function, the correctness of the above two series is obvious.
\end{proof}

\bibliographystyle{amsplain}

\end{document}